\documentclass{article}

\usepackage{times}
\usepackage{graphicx} 
\usepackage[FIGTOPCAP]{subfigure}

\usepackage{algorithm}\usepackage{algorithmic}

\usepackage{hyperref}

\usepackage{epsfig,fullpage}
\usepackage{epstopdf}
\usepackage{amsmath,amssymb,theorem,float,bm,enumerate,multirow}
\usepackage{rotating}
\usepackage{array}
\usepackage{color}

\allowdisplaybreaks

\newcommand{\beq}{\begin{equation}}
\newcommand{\eeq}{\end{equation}}

\newcommand\I{\mathbb{I}}

\newcommand\R{\mathbb{R}}

\renewcommand{\u}{\mathbf{u}}

\newcommand{\cL}{{\cal L}}

\newcommand{\vertiii}[1]{{\left\vert\kern-0.25ex\left\vert\kern-0.25ex\left\vert #1
    \right\vert\kern-0.25ex\right\vert\kern-0.25ex\right\vert}}

\DeclareMathOperator{\argmin}{argmin}

\newcounter{exampleI}
{\theorembodyfont{\rmfamily} \theoremstyle{plain} }

\newcounter{exampleII}
{\theorembodyfont{\rmfamily} \theoremstyle{plain} }

\newcounter{exampleIII}
{\theorembodyfont{\rmfamily} \theoremstyle{plain} }

{\theorembodyfont{\rmfamily} }
{\theorembodyfont{\rmfamily} \newtheorem{defn}{Definition}}
{\theorembodyfont{\rmfamily} }
\newtheorem{theo}{Theorem}
\newtheorem{rem}{Remark}

\newtheorem{lemm}{Lemma}

\newcommand{\proof}{\noindent{\itshape Proof:}\hspace*{1em}}

\newcommand{\qed}{\nolinebreak[1]~~~\hspace*{\fill} \rule{5pt}{5pt}\vspace*{\parskip}\vspace*{1ex}}

\newcommand {\commentout}[1] {}

\def\ints{{{\rm Z} \kern -.35em {\rm Z} }}  
\def\smallints{{{\rm Z} \kern -.3em {\rm Z} }}  
\def\pints{{{\rm I} \kern -.15em {\rm N} }}      
\newcommand{\reals}{\mathbb R}

\def\cplx{{{\rm I} \kern -.45em {\rm C} }}       
\def\l2{\rm {\mathcal L}^{2}(\reals)}            

\newtheorem{nad}{Notation and Definitions}[section]

\newtheorem{corollary}{Corollary}

\newcommand{\be}{\begin{eqnarray}}
\newcommand{\ee}{\end{eqnarray}}
\newcommand{\bea}{\begin{eqnarray}}
\newcommand{\eea}{\end{eqnarray}}
\newcommand{\beaa}{\begin{eqnarray*}}
\newcommand{\eeaa}{\end{eqnarray*}}
\newcommand{\bnad}{\begin{nad}}
\newcommand{\enad}{\end{nad}}

\newcommand{\sign}{{\mbox{\rm sign}}}

\renewcommand{\widehat}{\hat}

\title{Generalized Direct Change Estimation in Ising Model Structure}
\author{Farideh Fazayeli \qquad Arindam Banerjee  \vspace*{2mm}
\\
\{farideh,banerjee@cs.umn.edu\}\vspace*{2mm}\\
Department of Computer Science \& Engineering\\
University of Minnesota, Twin Cities}
\date{}
\begin{document}
\maketitle

\begin{abstract}
We consider the problem of estimating change in the dependency structure between two $p$-dimensional Ising models, based on respectively $n_1$ and $n_2$ samples drawn from the models. The change is assumed to be structured, e.g., sparse, block sparse, node-perturbed sparse, etc., such that it can be characterized by a suitable (atomic) norm. We present and analyze a norm-regularized estimator for directly estimating the change in structure, without having to estimate the structures of the individual Ising models. The estimator can work with any norm, and can be generalized to other graphical models under mild assumptions. We show that only one set of samples, say $n_2$, needs to satisfy the sample complexity requirement for the estimator to work, and the estimation error decreases as $\frac{c}{\sqrt{\min(n_1,n_2)}}$, where $c$ depends on the Gaussian width of the unit norm ball. For example, for $\ell_1$ norm applied to $s$-sparse change, the change can be accurately estimated with $\min(n_1,n_2)=O(s \log p)$ which is sharper than an existing result $n_1= O(s^2 \log p)$ and $n_2 = O(n_1^2)$. Experimental results illustrating the effectiveness of the proposed estimator are presented.
\end{abstract}

\section{Introduction}
\label{sec:intro}

Over the past decade, considerable progress has been made on estimating the statistical dependency structure of graphical models based on samples drawn from the model. In particular, such advances have been made for Gaussian graphical models, Ising models, Gaussian copulas, as well as certain multi-variate extensions of general exponential family distributions including multivariate Poisson models~\cite{baga08,kahs09,mebu06,rawl10,rwry11,yglr12}.

In this paper, we consider Ising models and focus on the problem of {\em estimating changes in Ising model structure}:
given two sets of samples $ \mathfrak{X}_1^{n_1} = \{{\bf{x}}^1_i\}_{i=1}^{n_1}$ and $\mathfrak{X}_2^{n_2}=\{{\bf{x}}^2_i\}_{i=1}^{n_2}$ respectively drawn from two $p$-dimensional Ising models with true parameters $\theta_1^*$ and $\theta_2^*$, where $\theta_1^*,\theta_2^* \in \R^{p \times p}$, the goal is to estimate the change $\delta \theta^* = (\theta_1^* - \theta_2^*)$. In particular, we focus on the situation when the change $\delta \theta^*$ has structure, such as sparsity, block sparsity, or node-perturbed sparsity, which can be characterized by a suitable (atomic) norm~\cite{crpw12,mlfw14}. However, the individual model parameters $\theta_1^*,\theta_2^*$ need not have any specific structure, and they may both correspond to dense matrices. The goal is to get an estimate $\delta \hat{\theta}$ of the change $\delta \theta^*$ such that the estimation error $\Delta = (\delta \hat{\theta} - \delta \theta^* )$ is small. Such change estimation has  potentially wide range of applications including identifying the changes in the neural connectivity networks, the difference between plant trait interactions at different climate conditions, and the changes in the stock market dependency structures.


One can consider two broad approaches for solving such change estimation problems: (i) {\em indirect change estimation}, where we estimate $\widehat{\theta}_1$ and $\widehat{\theta}_2$ from two sets of samples separately and obtain $\delta \widehat{\theta} = (\widehat{\theta}_1 - \widehat{\theta}_2)$, or (ii) {\em direct change estimation}, where we directly estimate $\delta \hat{\theta}$ using the two sets of samples, without estimating ${\theta}_1$ and ${\theta}_2$ individually. In a high dimensional setting, recent advances~\cite{call11,rawl10,rwry11} illustrate that accurate estimation of the parameter $\theta^*$ of an Ising model depends on how sparse or otherwise structured the true parameter $\theta^*$ is.
For example, if both $\theta_1^*$ and $\theta_2^*$ are sparse and the samples $n_1, n_2$ are sufficient to estimate them accurately~\cite{rawl10}, indirect estimation of $\delta \hat{\theta}$ should be accurate.
However, if the individual parameters $\theta_1^*$ and $\theta_2^*$ are somewhat dense, and the change $\delta \theta^*$ has considerably more structure, such as block sparsity (only a small block has changed) or node perturbation sparsity (only edges from a few nodes have changed)~\cite{mlfw14}, direct estimation may be considerably more efficient both in terms of the number of samples required as well as the computation time.

{\bf Related Work:} In recent work, Liu et al.~\cite{liss14} proposed a direct change estimator for graphical models based on the ratio of the probability density of the two models~\cite{gshs09,kahs09,snkp08,susk12,vaiz15}. They focused on the special case of $L_1$ norm, i.e., $\delta \theta^* \in \mathbb{R}^{p^2}$ is sparse, and provided non-asymptotic error bounds for the estimator along with a sample complexity of $n_1=O(s^2 \log p)$ and $n_2=O(n_1^2)$ for an unbounded density ratio model, where $s$ is the number of the changed edges with $p$ being the number of variables.
Liu et al.~\cite{liss15} improved the sample complexity to $\min(n1, n2)=O(s^2 \log p)$ when a bounded density ratio model is assumed. Zhao et~al.~\cite{zhcl14} considered estimating direct sparse changes in Gaussian graphical models (GGMs). Their estimator is specific to GGMs and can not be applied to Ising models.

{\bf Our Contributions:} We consider general structured direct change estimation, while allowing the change to have any structure which can be captured by a suitable (atomic) norm $R(\cdot)$. Our work is a considerable generalization of the existing literature which can only handle sparse changes, captured by the $L_1$ norm. In particular, our work now enables estimators for more general structures such as group/block sparsity, hierarchical group/block sparsity, node perturbation based sparsity, and so on \cite{bcfs14,crpw12,mlfw14,nrwy12}. Interestingly, for the unbounded density ratio model, our analysis yields sharper bounds for the special case of $\ell_1$ norm, considered by Liu et al.~\cite{liss14}. In particular, when $\delta \theta^*$ is sparse and our estimator is run with $L_1$ norm, we get a sample complexity of $n_1=n_2=O(s \log p)$ which is sharper than  $n_1=O(s^2 \log p)$ and $n_2=O(n_1^2)$ in \cite{liss14}.

The regularized estimator we analyze is broadly a Lasso-type estimator, with key important differences: the objective {\em does not} decompose additively over the samples,
and the objective depends on samples from two distributions. The estimator builds on the density ratio estimator in~\cite{liss14}, but works with general norm regularization \cite{bcfs14,crpw12,nrwy12} where the regularization parameter $\lambda_{n_1,n_2}$ depends on the sample size for both Ising models.
Our analysis is quite different from the existing literature in change estimation. Liu et.~\cite{liss14} build on the primal-dual witness approach of Wainwright~\cite{wain09}, which is effective for the special case of $L_1$ norm. Our analysis is largely geometric, where generic chaining~\cite{tala14} plays a key role, and our results are in terms of Gaussian widths of suitable sets associated with the norm~\cite{bcfs14,crpw12}.

The rest of the paper is organized as follows. In Section \ref{sec:directchange}, we introduce the direct change estimator based on the ratio of the probability density of the Ising models. In Section \ref{sec:theo}, we establish statistical consistency of the direct change estimator, and conclude in Section \ref{sec:concl}.

\section{Generalized Direct Change Estimation}
\label{sec:directchange}
\begin{figure}[h!]
\centering
\subfigure[$\theta_1$]{\includegraphics[trim = 10mm 0mm 15mm 10mm, clip, width = 0.32\columnwidth]{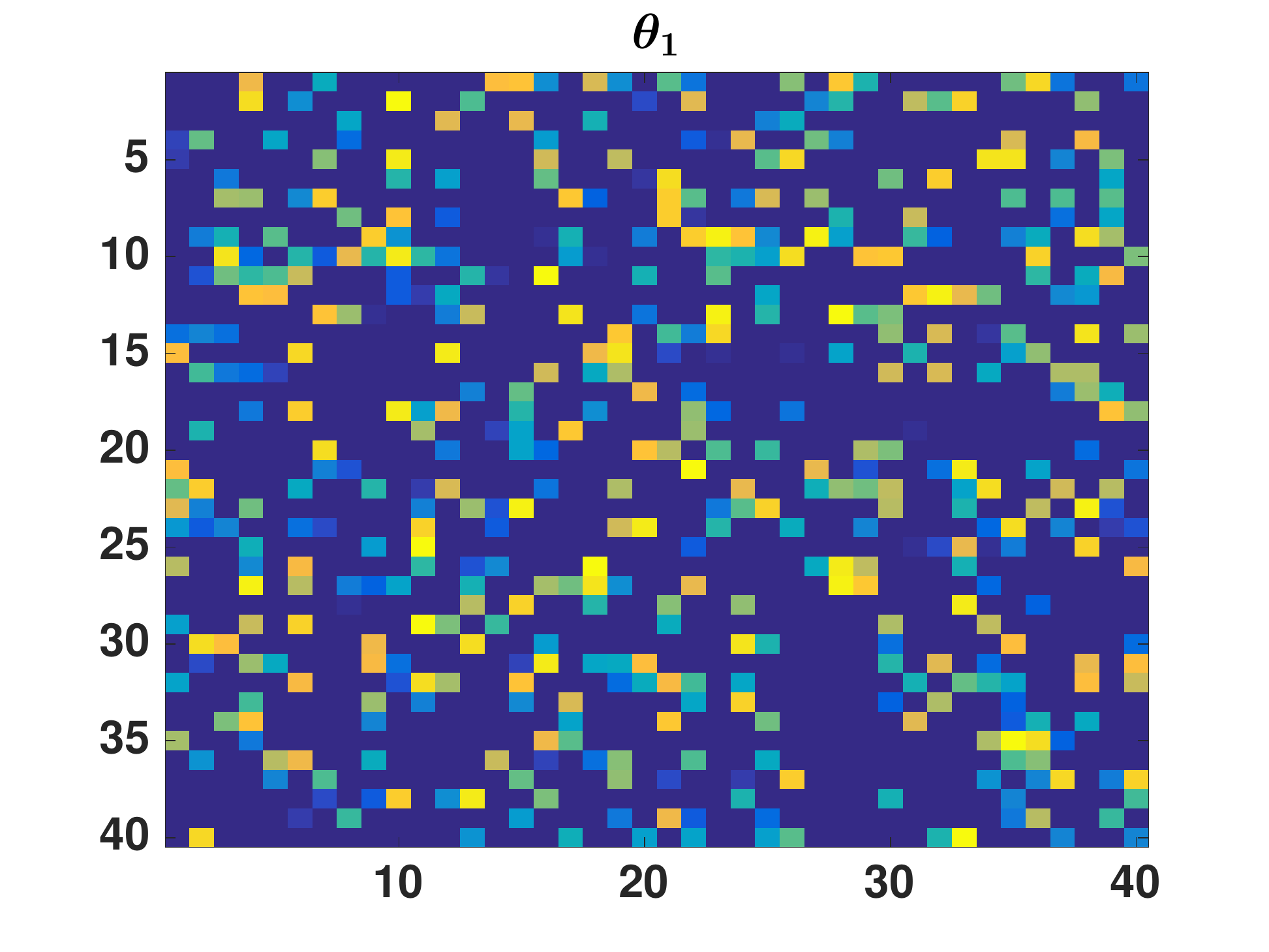}}
\subfigure[$\theta_2$]{\includegraphics[trim = 10mm 0mm 15mm 10mm, clip, width = 0.32\columnwidth]{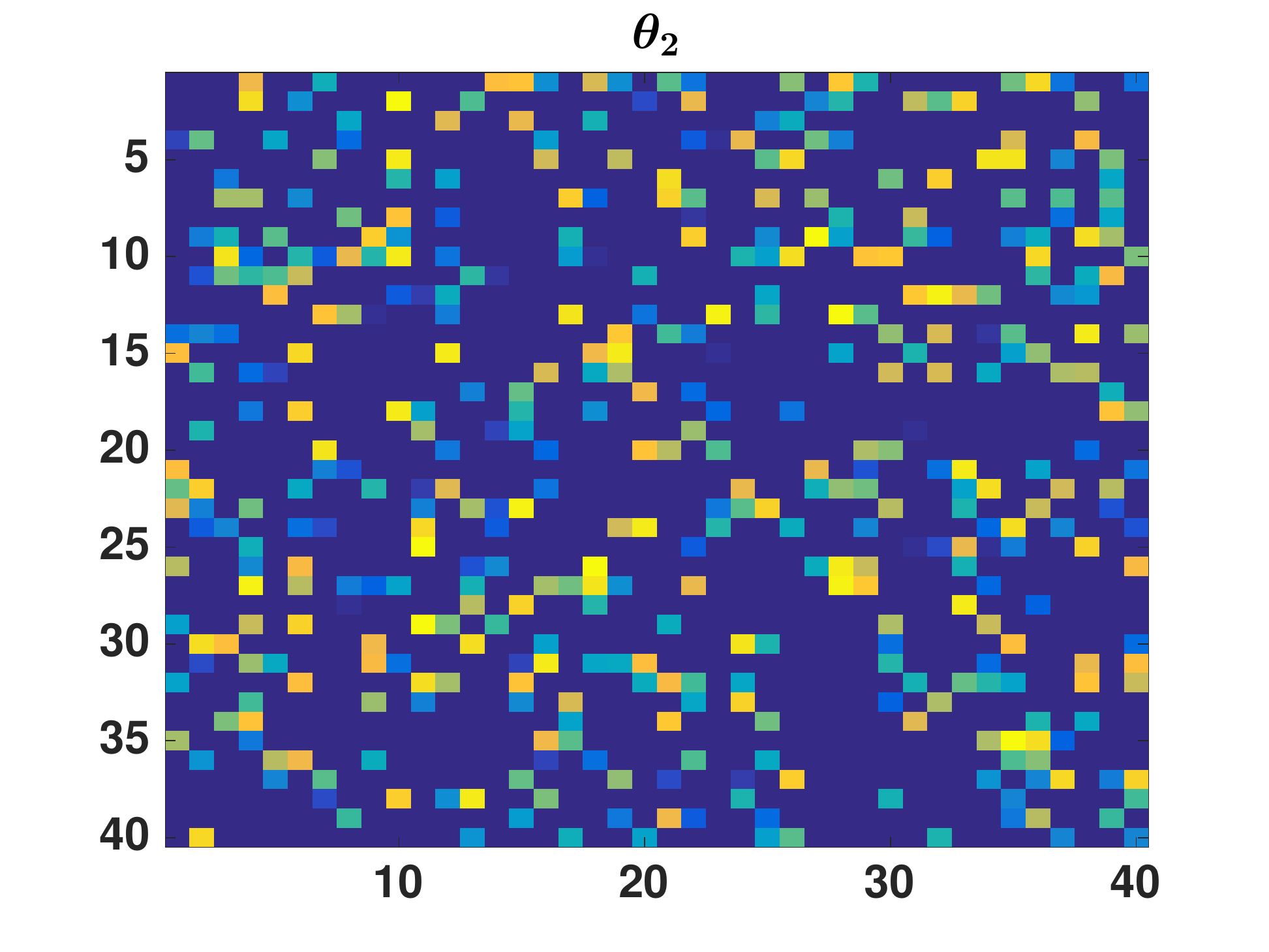}}
\subfigure[$\delta \theta=\theta_1-\theta_2$]{\includegraphics[trim = 10mm 0mm 15mm 10mm, clip, width = 0.32\columnwidth]{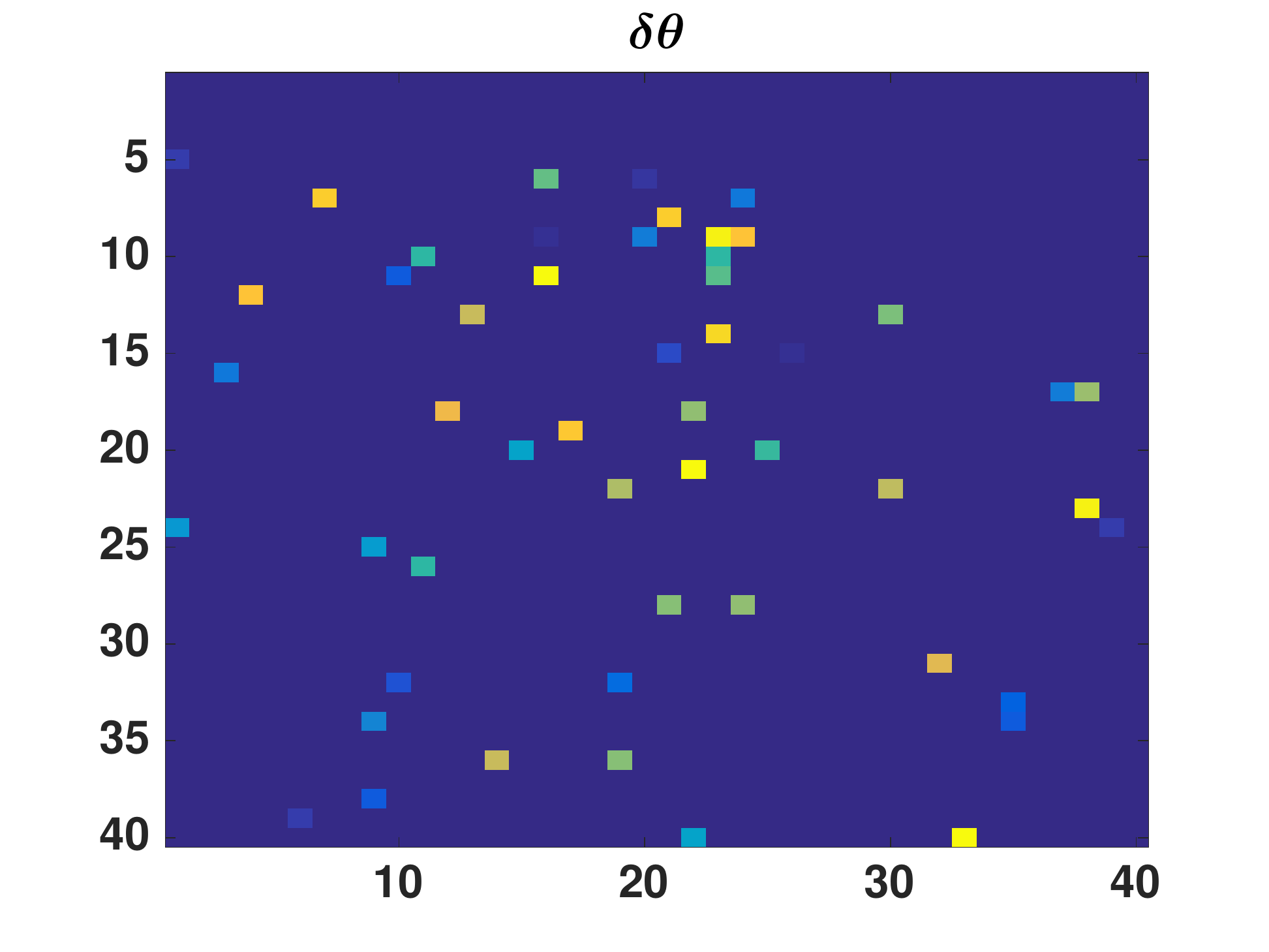}}
\\
(a) Sparse Structure

\subfigure[$\theta_1$]{\includegraphics[trim = 10mm 0mm 15mm 10mm, clip, width = 0.32\columnwidth]{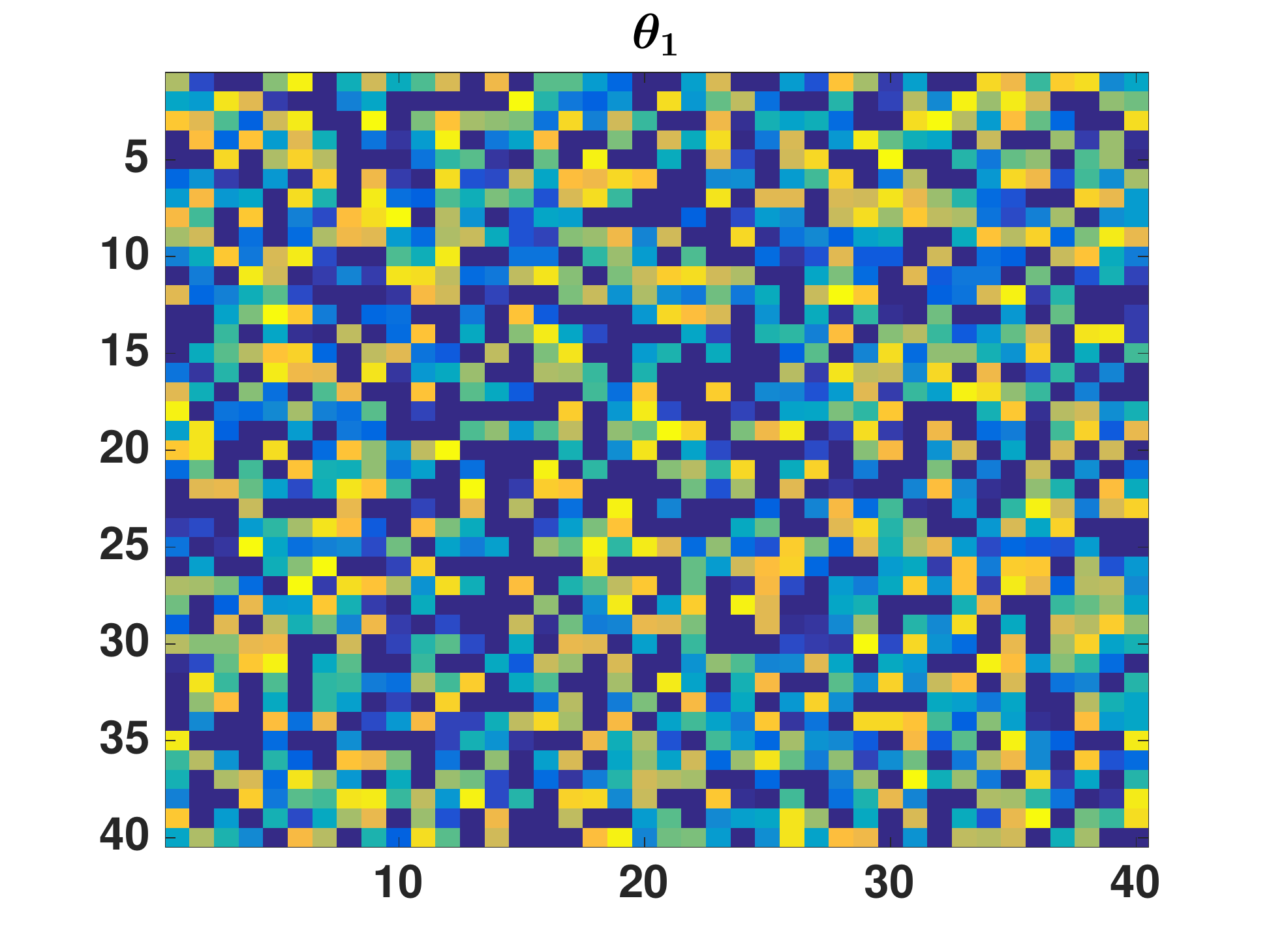}}
\subfigure[$\theta_2$]{\includegraphics[trim = 10mm 0mm 15mm 10mm, clip, width = 0.32\columnwidth]{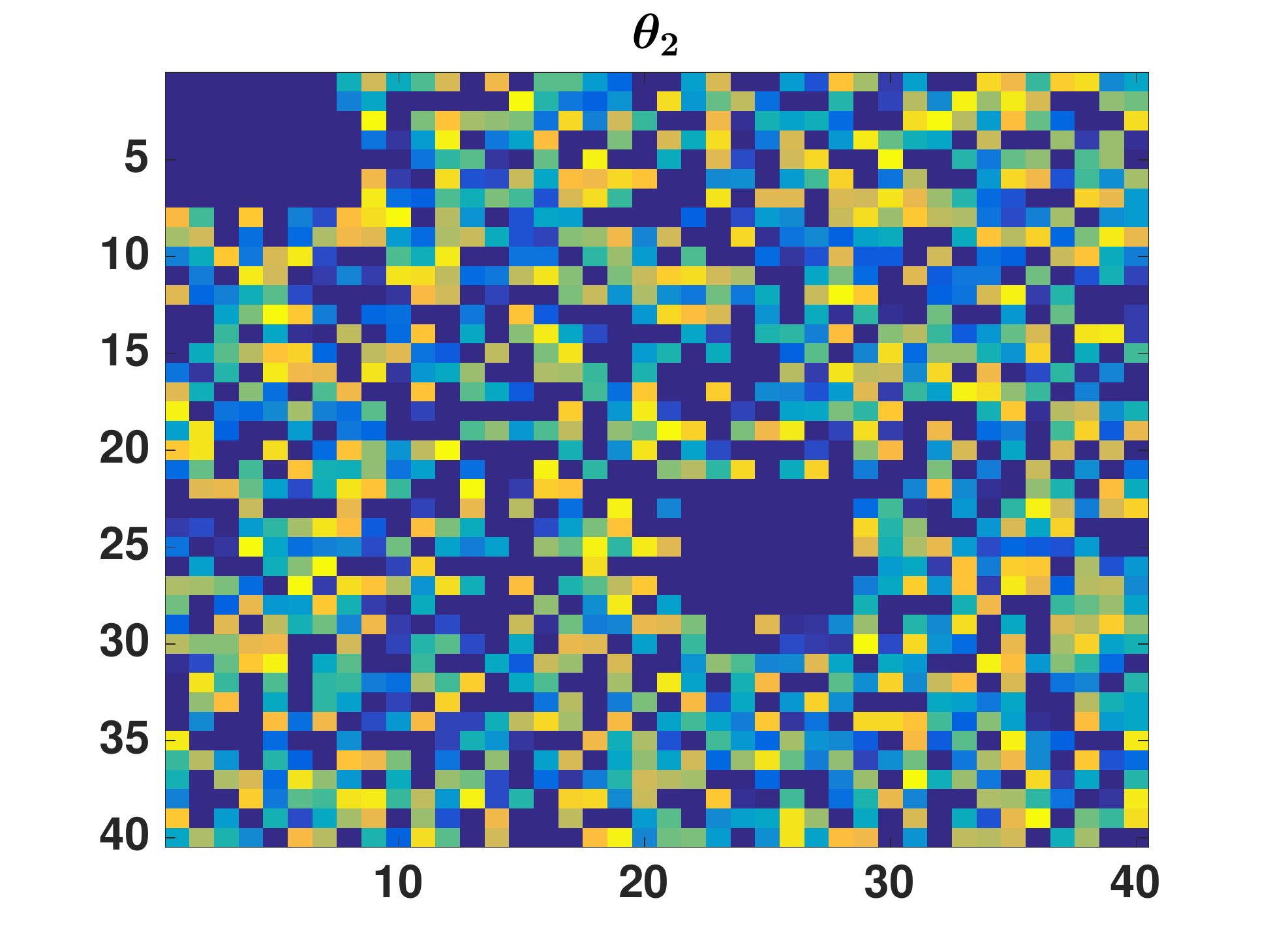}}
\subfigure[$\delta \theta=\theta_1-\theta_2$]{\includegraphics[trim = 10mm 0mm 15mm 10mm, clip, width = 0.32\columnwidth]{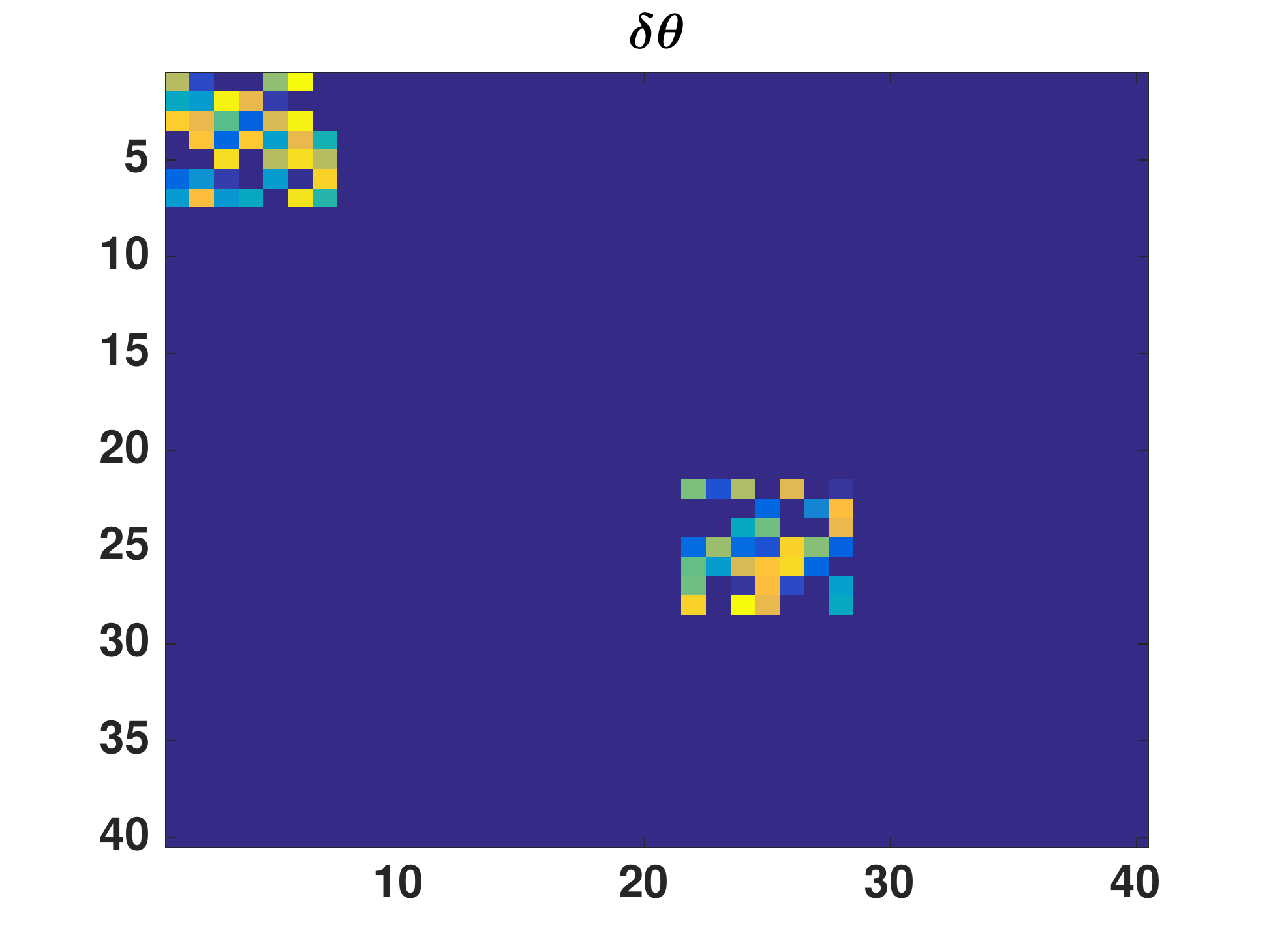}}
\\
(b) Group Sparsity Structure

\subfigure[$\theta_1$]{\includegraphics[trim = 10mm 0mm 15mm 10mm, clip, width = 0.32\columnwidth]{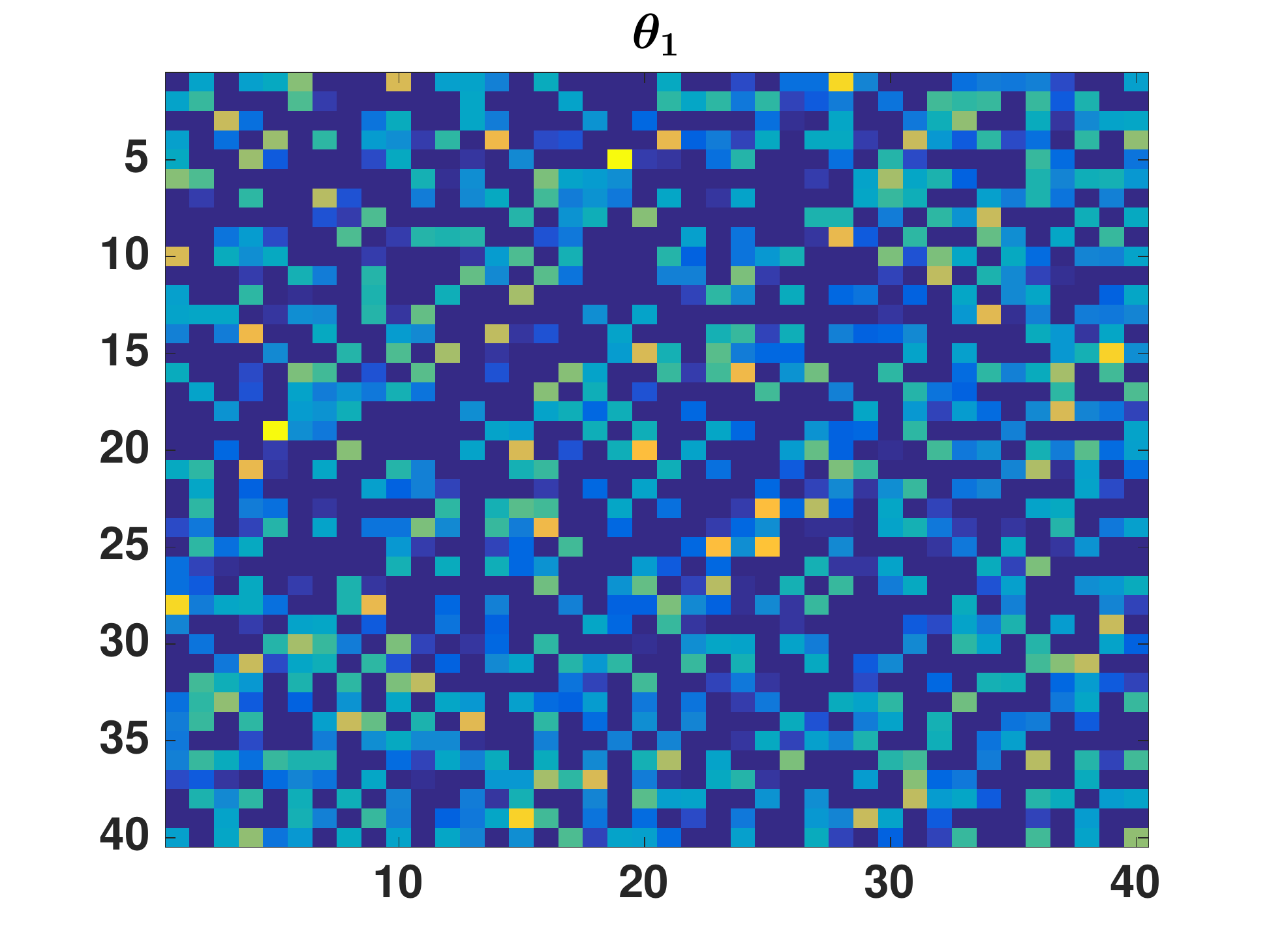}}
\subfigure[$\theta_2$]{\includegraphics[trim = 10mm 0mm 15mm 10mm, clip, width = 0.32\columnwidth]{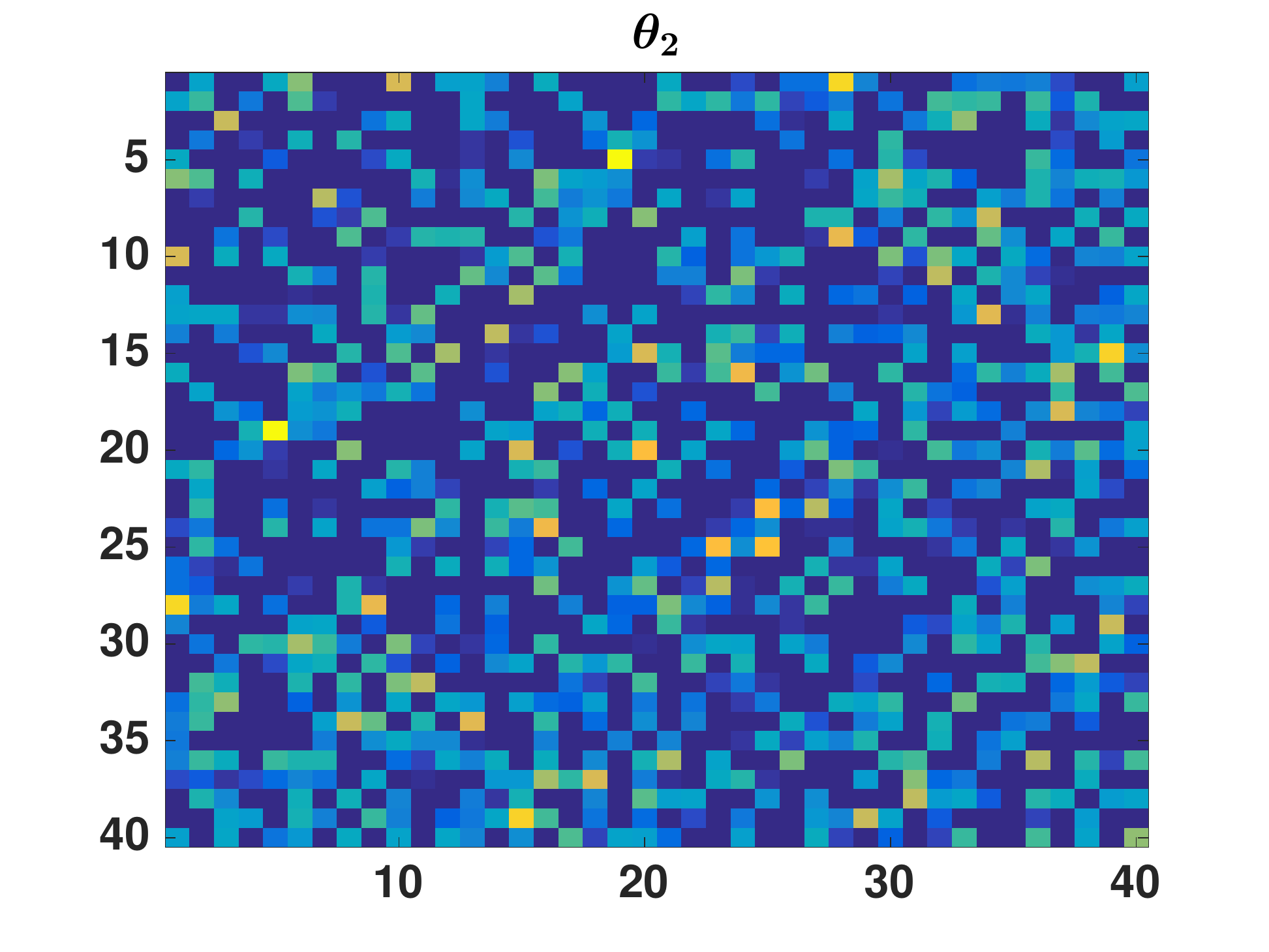}}
\subfigure[$\delta \theta=\theta_1-\theta_2$]{\includegraphics[trim = 10mm 0mm 15mm 10mm, clip, width = 0.32\columnwidth]{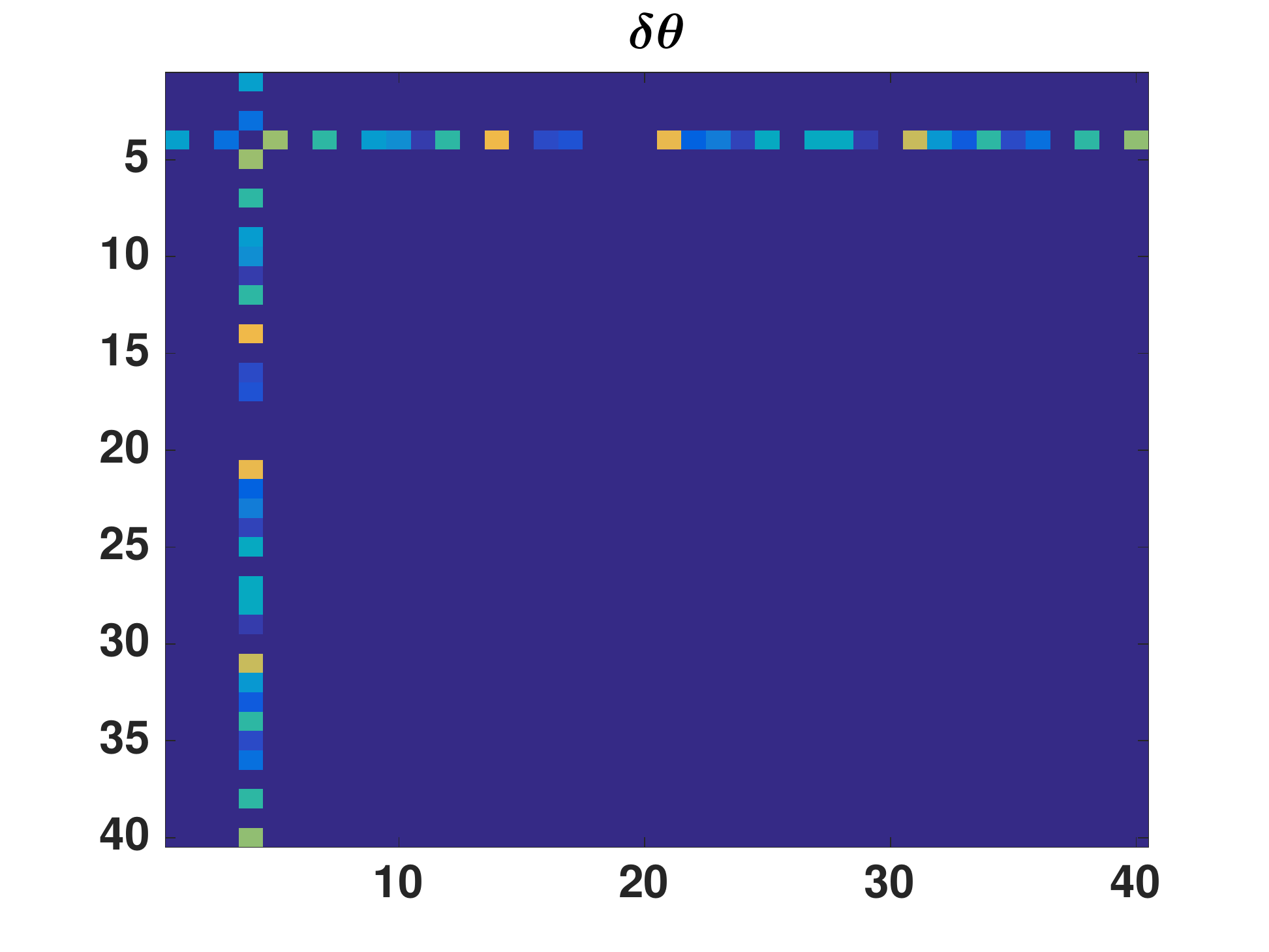}}
\\
(c) Node Perturbation Structure

\vspace{-0.3cm}
\caption{Examples of $\delta \theta$ with different structures. (a) This example shows the sparsity structure of $\delta \theta$, where a few edges has changed between $\theta_1$ and $\theta_2$. (b) This example presents the group sparsity structure, where a block of edges has changed between $\theta_1$ and $\theta_2$. (c) This example shows the node perturbation structure, where the connection of one node has changed. Blue represents zeros.}
\vspace{-0.3cm}
\label{fig:exmaples}
\end{figure}

We consider the following optimization problem
\begin{align}
\underset{{\delta \theta}}{\argmin} ~~\cL(\delta \theta; \mathfrak{X}_1^{n_1}, \mathfrak{X}_2^{n_2}) + \lambda_{n_1,n_2} R(\delta \theta),
\label{eq:changeestimateopt}
\end{align}
where $ \mathfrak{X}_1^{n_1} = \{{\bf{x}}^1_i\}_{i=1}^{n_1}$ and $\mathfrak{X}_2^{n_2}=\{{\bf{x}}^2_i\}_{i=1}^{n_2}$ are two sets of i.i.d~ binary samples drawn from from Ising graphical models with parameter $\theta_1^*$ and $\theta_2^*$, respectively, each ${\bf{x}}^1_i$ and ${\bf{x}}^2_i$ are $p-$dimensional vectors, and $n_1, n_2$ are the respective sample sizes.

In this Section, we first give a brief background on Ising model selection. Then, we explain how to develop the loss function $\cL(\delta \theta; \mathfrak{X}_1^{n_1}, \mathfrak{X}_2^{n_2})$ based on the density ratio \cite{gshs09,kahs09,snkp08,vaiz15} to directly estimate $\delta \theta = \theta_1 - \theta_2$, and finally we describe how to solve the optimization problem \eqref{eq:changeestimateopt} for any norm $R(\delta \theta)$.

\subsection{Ising Model}
Let $X = (X_1,X_2, \cdots ,X_p)$ denote a random vector in which each variable $X_s \in \{-1, 1\}$.
Let $G = (V, E)$ be an undirected graph with vertex set $V = \{1, \cdots, p\}$ and edge set $E$
whose elements are unordered pairs of distinct vertices. The pairwise Ising Markov random field associated with the graph $G$ over the random vector $X$ is 
\begin{align}
P(X= {\bf{x}} | \theta^*) &= \frac{1}{Z(\theta^*)} \exp  \{ \sum_{s,t \in E} \theta_{s,t}^* x_s x_t  \}\\
&= \frac{1}{Z(\theta^*)} \exp \{\langle \theta^*, T({\bf{x}}) \rangle \} 
\label{eq:isingvec}
\\
&= \frac{1}{Z(\Theta^*)} \exp \{ {\bf{x}} ^T \Theta^*{\bf{x}}  \}
\label{eq:isingMat}
\end{align}
where $T({\bf{x}}) = \{ x_s x_t\}_{s,t=1}^{p}$ is a vector of size $m=p^2$, $\theta^*= \{\theta_{s,t}^*\}_{s,t=1}^{p}\in \mathbb{R}^m$ and $\langle. ,. \rangle$ is the inner product operator, and $\Theta^* \in \mathbb{R}^{p \times p}$ where $\Theta^*_{s,t} = \theta_{s,t}^*$.
Note that basic Ising models also have non-interacting terms like $\alpha_s x_s$ and we are assuming these terms are zero, and they do not affect the dependency structure.

The parameter $\theta^*$ associated with the structure of the graph $G$ reveals the statistical conditional independence structure among the variables i.e., if $\theta^*_{s,t} = 0$, then feature $X_s$ is conditionally independent of $X_t$ given all other variables and there is no edge in the graph $G$.

The partition function, $Z(\theta^*)$, plays the role of a normalizing constant, ensuring that the probabilities add up to one which is defined as
\beq
Z(\theta^*) = \sum_{{\bf{x}} \in \mathcal{X}} \exp \{ \langle \theta^*, T({\bf{x}}) \rangle\} = \exp\{\Psi(\theta^*)\},
\eeq
where $\mathcal{X}$ be the set of all possible configurations of $X$.

\subsection{Loss Function} 
Here, we build the loss function based on equation \eqref{eq:isingvec}. Similarly, one can rewrite the loss function based on \eqref{eq:isingMat} if the regularization function is over matrices. 
Consider two Ising models with parameters $\theta_1^* \in \mathbb{R}^{p^2}$ and $\theta_2^* \in \mathbb{R}^{p^2}$.
Following Liu et. al \cite{lqgs14, liss14}, a direct estimate for the changes detection problem based on density ratio can be posed as follows
\begin{align}
r(X = {\bf{x}}|\delta \theta) =  \frac{p(X = {\bf{x}}|\theta_1)}{p(X = {\bf{x}}|\theta_2)} = 
\underbrace{\frac{\exp\{ \langle T({\bf{x}}), \theta_1 \rangle \}}{\exp\{ \langle T({\bf{x}}), \theta_2 \rangle \}}}_{r^*({\bf{x}}|\delta \theta)} 
\underbrace{\frac{Z(\theta_2)}{Z(\theta_1)}}_{1/Z(\delta \theta)}
 = \frac{\exp\{ \langle T({\bf{x}}), \delta \theta \rangle) \}}{Z(\delta \theta)},
\end{align}
where the parameter $\delta \theta = \theta_1 - \theta_2$ encodes the change between two graphical models $\theta_1$ and $\theta_2$. 

First, we show that $Z(\delta \theta) =  E_{X\sim q} [e ^{\langle T(X), \delta \theta \rangle}]$:
\begin{align}
Z(\delta \theta) = \frac{Z(\theta_1)}{Z(\theta_2)} = \frac{1}{Z(\theta_2)} \sum_{{\bf{x}} \in \mathcal{X}} e ^{\langle T({\bf{x}}), \theta_1 \rangle} 
& = \frac{1}{Z(\theta_2)}  \sum_{{\bf{x}} \in \mathcal{X}} e ^{\langle T({\bf{x}}), \theta_2 \rangle} \frac{e ^{\langle T({\bf{x}}), \theta_1 \rangle} }{e ^{\langle T({\bf{x}}), \theta_2 \rangle} } \nonumber \\
 & = \sum_{{\bf{x}} \in \mathcal{X}} \underbrace{\frac{e ^{\langle T({\bf{x}} ), \theta_2 \rangle}}{Z(\theta_2)}}_{p({\bf{x}}|\theta_2)} e ^{\langle T({\bf{x}} ), \delta \theta \rangle} 
= E_{X\sim p(X|\theta_2)} [e ^{\langle T(X), \delta \theta \rangle}].
\end{align}

Next, using the samples $\mathfrak{X}_2^{n_2}$ from $p(X|\theta_2)$, we estimate $Z(\delta \theta)$ empirically as
\begin{align}
\hat{Z}(\delta \theta)  
= \frac{1}{n_2} \sum_{i=1}^{n_2} \exp\{ \langle T({\bf{x}}_i^2), \delta \theta \rangle \},
\end{align}
and the sample approximation of $r(X|\delta \theta)$ is given as
\begin{align}
 \hat{r}(X = {\bf{x}}|\delta \theta) =  \frac{r^*({X=\bf{x}}|\theta_1)}{\hat{Z}(\delta \theta)}  = \frac{\exp\{ \langle T({\bf{x}}), \delta \theta \rangle\}}{\frac{1}{n_2} \sum_{i=1}^{n_2} \exp\{ \langle T({\bf{x}}_i^2), \delta \theta \rangle \}}.
 \label{eq:approxRatio}
\end{align}

Using the fact that $r(X|\delta \theta^*) q(X|\theta_2^*) = p(X|\theta_1^*)$, we approximate $ \hat{r}(X|\delta \theta)$, by minimizing the $KL$ divergence,
\begin{align}
KL \left( p(X|\theta_1^*) \| \hat{r}(X|\delta \theta) p(X|\theta_2^*) \right)  &=  \sum_{{\bf{x}} \in \mathcal{X}}  p({\bf{x}}|\theta_1^*) \log \frac{p({\bf{x}}|\theta_1^*)}{p({\bf{x}}|\theta_2^*) \hat{r}({\bf{x}}|\delta \theta)} \nonumber \\
&= \underbrace{KL\left(p(X|\theta_1^*) \| p(X|\theta_2^*) \right)}_{\text{Constant}} - E_{X \sim p(X|\theta_1^*)} \left[\log \hat{r}(X|\delta \theta) \right] \label{eq:KLdiv2}
\end{align}

Thus, using the samples $\mathfrak{X}_1^{n_1}$ and $\mathfrak{X}_2^{n_2}$, we define the empirical loss function
\begin{align}
\cL(\delta \theta; \mathfrak{X}_1^{n_1}, \mathfrak{X}_2^{n_2}) = \frac{-1}{n_1} \sum_{i=1}^{n_1}  \log \hat{r}({\bf{x}}_i^1 | \delta \theta)   =\frac{-1}{n_1} \sum_{i=1}^{n_1} \langle T({\bf{x}}_i^1), \delta \theta \rangle  + \underbrace{\log \frac{1}{n_2} \sum_{i=1}^{n_2} \exp\{ \langle T({\bf{x}}_i^2), \delta \theta \rangle \}}_{\hat{\Psi}(\delta \theta)} \label{eq:KLloss} 
\end{align}

\begin{rem}
\normalfont
Note that the loss function \eqref{eq:KLloss} does not additively decompose over the samples. The second term in \eqref{eq:KLloss} is the logarithm over sum of a function of samples.
\end{rem}
\subsection{Optimization}

The optimization problem \eqref{eq:changeestimateopt} has a composite objective with a smooth convex term corresponding to the loss function \eqref{eq:KLloss} and a a potentially non-smooth convex term corresponding to the regularizer. 
In this section, we present an algorithm in the class of Fast Iterative Shrinkage-Thresholding Algorithms (FISTA) for efficiently solving the problem \eqref{eq:changeestimateopt} \cite{bete09}.
For convenience, we refer the loss function $\cL(\delta \theta; \mathfrak{X}_1^{n_1}, \mathfrak{X}_2^{n_2})$ as $\cL(\delta \theta)$ and we drop the subscript $\{n_1, n_2\}$ of $\lambda_{n_1, n_2}$.

%
%
%
One of the most popular methods for composite objective functions is in the class of FISTA where at each iteration we linearize the smooth term and minimize 
the quadratic approximation of the form 
\begin{align}
Q_{L}(\delta \theta, \delta \theta_t) :=  \cL(\delta \theta) + \left \langle \delta\theta - \delta\theta_t, \nabla \cL(\delta \theta_t)  \right \rangle 
 + \frac{L}{2} \|\delta \theta - \delta \theta_t\|_2^2 + \lambda R(\delta\theta),
\label{eq:quadApprox}
\end{align}
where $L$ denotes the Lipschitz constant of the loss function $\cL(\delta \theta)$.
Ignoring constant terms in $\delta\theta_t$, the unique minimizer of the above expression \eqref{eq:quadApprox} can be written as
\begin{align}
p_{L} (\delta \theta_t) = \underset{\delta\theta}{\argmin}~~Q_{\cL}(\delta \theta, \delta \theta_t)  
& ~=\underset{\delta\theta}{\argmin}~~\lambda R(\delta \theta) + \frac{L}{2}  
 \left \| \delta \theta- \left( \delta \theta_t - \frac{1}{L} \nabla  \cL(\delta \theta_t)  \right) \right\|_2^2 \nonumber \\
&~= \underset{\delta\theta}{\argmin}~~\frac{\lambda}{L} R(\delta \theta) + \frac{1}{2}  
 \left \| \delta \theta- \left( \delta \theta_t - \frac{1}{L} \nabla  \cL(\delta \theta_t)  \right) \right\|_2^2. \label{eq:minimizer} 
\end{align}
In fact, the updates of $\delta \theta$ is to compute certain proximal operators of the non-smooth term $R(.)$. In general, the proximal operator $\text{prox}_h({\bf{x}})$ of a closed proper convex function $h : \mathbb{R}^d \mapsto \mathbb{R} \cup \{+ \infty \}$ \cite{pabo14}  is defined as
\begin{align}
\text{prox}_h({\bf{x}}) = \underset{{\bf{u}}}{\argmin} \left( h({\bf{u}}) + \frac{1}{2} \| {\bf{u}}-{\bf{x}} \|_2^2  \right).
\label{eq:prox}
\end{align}
Thus, the unique minimizer \eqref{eq:minimizer} correspond to 
$\text{prox}_{\frac{\lambda}{L}R}\left( \delta \theta_t - \frac{1}{L} \nabla  \cL(\delta \theta_t)  \right)$ which has  rate of convergence of $O(1/t)$~\cite{nest05,pabo14}.



To improve the rate of convergence, we adapt the idea of FISTA algorithm~\cite{bete09}. 
The main idea is to iteratively consider the proximal operator $\text{prox}(.)$ at a specific linear combination of the previous two iterates $\{\delta\theta_t, \delta\theta_{t-1}\}$
\begin{align}
\xi_{t+1} = \delta\theta_t + \alpha_{t+1} \left(\delta\theta_t - \delta\theta_{t-1}\right),
\end{align}
instead of just the previous iterate $\delta\theta_t$. The choice of $\alpha_{t+1}$ follows Nesterov's accelerated gradient descent \cite{nest05,pabo14} and is detailed in Algorithm \ref{alg:ratio}. The iterative algorithm simply updates
\begin{align}
\delta\theta_{t+1} = \text{prox}_{\frac{\lambda}{L}R}\left( \xi_{t+1} - \frac{1}{L} \nabla  \cL(\xi_{t+1})  \right).
\end{align}
The algorithm has a rate of convergence of $O(1/t^2)$~\cite{bete09}.

\begin{algorithm}[tb]
   \caption{Generalized Direct Change Estimator}
   \label{alg:ratio}
\begin{algorithmic}
   \STATE {\bfseries Input:} $L_0 >0$,  $\mathfrak{X}_1^{n_1}, \mathfrak{X}_2^{n_2}$
   \STATE {\bfseries Step 0.} Set $\xi_1 = \delta \theta_0$, $t = 1$
   \STATE {\bfseries Step $t$.} ($t\geq 1$)  Find the smallest non-negative integers $i_t$ such 
   that with $\tilde{L} = 2^{i_t} L_{t-1}$
   \begin{align}
   \cL\left(p_{\tilde{L}}(\xi_t)\right) + R\left(p_{\tilde{L}}(\xi_t)\right) \leq Q_{\tilde{L}}\left(p_{\tilde{L}}(\xi_t), \xi_t \right).
   \end{align}
   \STATE Set $L_t = 2^{i_t} L_{t-1}$ and Compute
	\begin{align}
	\delta\theta_t &= \text{prox}_{\frac{\lambda}{L}R}\left( \xi_{t} - \frac{1}{L} \nabla  \cL(\xi_{t})  \right) \\
	\beta_{t+1} &= \frac{1+\sqrt{1+4 \beta_t^2}}{2} \\
	\xi_{t+1} &= \delta \theta_t + \left(\frac{\beta_{t}-1}{\beta_{t+1}} \right) (\delta\theta_t - \delta\theta_{t-1})
	\end{align}
\end{algorithmic}
\end{algorithm}

\subsection{Regularization Function} 
We assume that the optimal $\delta \theta^*$ is sparse or suitably `structured' where such structure can be characterized by having a low value according to a suitable norm $R(\delta \theta^*)$. In below, we provide a few examples of such a norm.

{\bf $L_1$ norm:} One example for $R(.)$ we will consider throughout the paper is the $L_1$
norm regularization. We use $L_1$ norm if only a few edges has changed (1st row in Figure \ref{fig:exmaples}).~In particular, we consider $R(\delta \theta) = \|\delta \theta\|_1$ if number of non-zeros entries in $\delta \theta^*$ is $s < p^2$. The $\text{prox}_{\frac{\lambda}{L}\|.\|_1}(.)$  is given by the elementwise soft-thresholding operation~\cite{sidu09} as
\begin{align}
\left[\text{prox}_{\frac{\lambda}{L}\|.\|_1}\right]_i  ({\bf{z}})= \sign({\bf{z}}_i).\max(0, {\bf{z}}_i-\frac{\lambda}{L}).
\end{align}

{\bf Group-sparse norm:} Another popular example we consider is the group-sparse norm. 
We use group lasso norm if a group of edges has changed (2nd row in Figure \ref{fig:exmaples}).
For some kinds of data, it is reasonable to assume that the variables can be clustered (or grouped) into types, which share similar connectivity or correlation patterns.
Let $\mathcal{G} = \{\mathcal{G}_1, \mathcal{G}_2, \cdots, \mathcal{G}_{N_G} \}$ denote a
collection of groups, which are subsets of variables. 
We assume that $\delta\Theta^*(s,t)=0$ for any variable $s \in G_g$ and for any variable $t \in G_h$.
In the group sparse setting for any subset $S_\mathcal{G} \subseteq \{1, 2, \cdots N_G \} $ with cardinality $|S_\mathcal{G}| = s_\mathcal{G}$, we assume that the parameter $\delta \Theta^*$ satisfies $\{\delta \Theta^*_{s,t} = 0 : s,t \in G_g ~\&~ g \not \in S_\mathcal{G}\}$. 
We will focus on the case when $R(\delta\Theta) = \sum_{g=1}^{N_G} \| \delta \Theta(s,t) : s,t \in G_g \|_F$~\cite{masm09}. Let $\delta\Theta_{G_g}$ bd the sub-matrix of $\delta\Theta$ covering nodes in $G_g$. Proximal operator is given by the group specific soft-thresholding operation. 
\begin{align}
\left[\text{prox}_{\frac{\lambda}{L}R} \right]_{g} (\delta\Theta)=
\frac{\max(\|\delta\Theta_{G_g}\|_F - \frac{\lambda}{L} ,0)}{\|\delta\Theta_{G_g}\|_F}.
\end{align}

{\bf Node perturbation:} Another example is the row-column overlap norm (RCON)~\cite{mlfw14} to capture perturbed nodes i.e., nodes that have a completely different connectivity pattern to other nodes among two networks (3rd row in Figure \ref{fig:exmaples}). A special case of RCON we are interested is $\sum_{i=1}^p \|V_i\|_q$  where $\delta \Theta = V+V^T$, and $V_i$ is the $i-$th column of matrix $V$. This norm can be viewed as overlapping group lasso~\cite{mlfw14} and thus can be solved by applying Algorithm \ref{alg:ratio} with proximal operator for overlapping group lasso~\cite{yuly11}. Also, we can write problem \eqref{eq:changeestimateopt} as a constrained optimization
\begin{align}
& \underset{{\delta \Theta, V}}{\argmin} ~~\cL(\delta \Theta; \mathfrak{X}_1^{n_1}, \mathfrak{X}_2^{n_2}) + \lambda_1 \delta \Theta \|_1 + \lambda_{n_1,n_2} \sum_{i=1}^p \|V_i\|_q \nonumber \\
& \text{s.t.} \qquad \delta \Theta = V+V^T,
\end{align}
and solve it by applying in-exact ADMM techniques~\cite{mlfw14}.

\section{Statistical Recovery for Generalized Direct Change Estimation}
\label{sec:theo}
\label{sec:theoanal}
Our goal is to provide non-asymptotic bounds on $\|\Delta\|_2 = \|\delta \theta^* - \delta \hat{\theta}\|_2$ between the true parameter $\delta \theta^*$ and the minimizer  $\delta \hat{\theta}$ of \eqref{eq:changeestimateopt}. In this section, we describe various aspects of the problem, introducing notations along the way, and 
highlight our main result.


\subsection{Background and Assumption}
\label{sec:anal-back}

{\bf Gaussian Width:} In several of our proofs, we use the concept of Gaussian width~\cite{crpw12,gord88}, which is defined as follows.

\begin{defn}
For any set $A \in \mathbb{R}^p$, the \emph{Gaussian width} of the set $A$ is defined as:
\beq
w(A) = E_g\left[\sup_{u \in A} \langle g, u \rangle\right]~.
\eeq
where the expectation is over $g \sim N(0,\I_{p \times p})$, a vector of independent zero-mean unit-variance Gaussian random variable.
\label{def:gw}
\end{defn}
The Gaussian width $w(A)$ provides a geometric characterization of the size of the set $A$. 
Consider the Gaussian process $\{ Z_u \}$  where the constituent Gaussian random variables $Z_u = \langle u, g \rangle$ are indexed by $u \in A$, and $g \sim N(0,\I_{p \times p})$. Then the Gaussian width $w(A)$ can be viewed as the expectation of the supremum of the Gaussian process $\{ Z_u \}$. Bounds on the expectations of Gaussian and other empirical processes have been widely studied in the literature, and we will make use of generic chaining for some of our analysis \cite{bolm13,ledo13,tala05,tala14}.

{\bf The Error Set:} Consider solving the problem \eqref{eq:changeestimateopt}, under assumption $\lambda_{n_1,n_2} > \beta R^*\left(\nabla \cL(\delta \theta^*; \mathfrak{X}_1^{n_1}, \mathfrak{X}_2^{n_2})\right)$,  where $\beta > 1$ and $R^*(.)$ is the dual norm of $R(.)$. Banerjee et al. \cite{bcfs14} show that for any convex loss function the error vector $\Delta = (\delta \theta^* - \delta \hat{\theta})$ lies in a restricted set that is characterized as
\begin{align}
E_r = E_r(\delta \theta^*,\beta) 
= \left\{ \Delta \in \R^p ~ \left| ~ R(\delta \theta^* + \Delta) \leq R(\delta \theta^*) + \frac{1}{\beta} R(\Delta) \right. \right\}~. \label{eq:error1}
\end{align}

{\bf Restricted Strong Convexity (RSC) Condition:} The sample complexity of the problem \eqref{eq:changeestimateopt} depends on the RSC condition \cite{nrwy12}, which ensures that the estimation problem is strongly convex in
the neighborhood of the optimal parameter \cite{bcfs14, nrwy12}. A convex loss function satisfies the RSC condition in $C_r = cone(E_r)$, i.e., $\forall \Delta \in C_r$, if there exists a suitable constant $\kappa$ such that
\begin{align}
\delta \cL(\delta \theta^{*}, u) :=  \cL(\delta \theta^{*} + u) - \cL(\delta \theta^{*}) - \langle \nabla \cL(\delta \theta^{*}), u \rangle  \geq  \kappa \|u\|_{2}^{2}
 \label{eq:RSC}
\end{align}

{\bf Deterministic Recovery Bounds:} If the RSC condition is satisfied on the error set $C_r$ and $\lambda_{n_1,n_2}$ satisfies the assumptions stated earlier, for any norm $R(.)$, Banerjee et al. \cite{bcfs14} show a deterministic upper bound for $\| \Delta \|_2$ in terms of $\lambda_{n_1,n_2}$, $\kappa$, and the norm compatibility constant $\Psi(C_r) = \sup_{\u \in C_r} \frac{R(\u)}{\| \u \|_2}$, as
\beq
\| \Delta \|_2 \leq \frac{1+\beta}{\beta} \frac{\lambda_{n_1,n_2}}{\kappa}  \Psi(C_r)~.
\label{eq:recover}
\eeq

{ \bf Smooth Density Ratio Model Assumption:} For any vector ${\bf{u}}$ such that $\|{\bf{u}}\|_2 \leq \| \delta \theta^*\|_2$ and every $\epsilon \in R$, the following inequality holds:
\begin{align}
E_{X \sim p(X|\theta_2)} [\exp \{ \epsilon ~ r(X | \delta \theta^*+ {\bf{u}}) - 1 \}] \leq \exp \{ \epsilon^2 \}. \nonumber
\end{align}

A similar assumption is used in the analysis of Liu et al. \cite{liss14}.

\begin{rem}
\normalfont
 Bounded density ratio is a special case satisfying the smooth density ratio assumption. Lemma \ref{lem:bndratiocond} shows a sufficient condition under which the density ratio is bounded.
\end{rem}
  \begin{lemm}
Consider two Ising Model with true parameters $\theta_1^*$ and $\theta_2^*$. Let $d_1, d_2 \gg s$ where $\|\theta_1^*\|_0=d_1$, $\|\theta_2^*\|_0=d_2$, and $\|\delta\theta^*\|_0 = s$.
Assume
 \begin{align}
 \min_{i,j=1\cdots p}( | \theta_1^*(i,j) |) & \geq \frac{1}{d_1-1} - \frac{c_1}{(d_1-1)s} 
  \label{eq:alpha1cond}
 \\ 
 \min_{i,j=1\cdots p}( | \theta_2^*(i,j) |)  & \geq \frac{1}{d_2-1} -\frac{c_2}{(d_2-1)s},
  \label{eq:alpha2cond}
 \end{align}
 where $c_1$ and $c_2$ are positive constants.
 Then the density ratio $r(X = {\bf{x}}|\delta \theta^*)$ is bounded.
 \label{lem:bndratiocond}
\end{lemm}
Note that if individual graphs are dense, then the conditions \eqref{eq:alpha1cond} and \eqref{eq:alpha2cond} are satisfied and as a result the smooth density ratio is satisfied.

\begin{rem}
\normalfont 
In this paper, we focus on the Ising graphical model. But, our statistical analysis holds for any graphical models that satisfy the above mentioned assumption. Through our analysis, no assumption is required on the individual graphical models.
\end{rem}

%

\subsection{ Bounds on the regularization parameter}
To get the recovery bound \eqref{eq:recover} above, one needs to have $\lambda_{n_1, n_2} \geq \beta R^*\left(\nabla \cL(\delta \theta^*; \mathfrak{X}_1^{n_1}, \mathfrak{X}_2^{n_2})\right)$.
However, the bound on $\lambda_{n_1, n_2}$ depends on unknown quantity $\delta \theta^*$ and the samples $\mathfrak{X}_1^{n_1}, \mathfrak{X}_2^{n_2}$ and is hence random.
To overcome the above challenges, one can bound the expectation $E[R^*\left(\nabla \cL(\delta \theta^*; \mathfrak{X}_1^{n_1}, \mathfrak{X}_2^{n_2})\right)]$ over all samples of size $n_1$ and $n_2$, and obtain high-probability deviation bounds. 
The goal is to provide a sharp bound on $\lambda_{n_1,n_2}$ since the error bound in \eqref{eq:recover} is directly proportional to $\lambda_{n_1,n_2}$.

In theorem \ref{theo:lambdabnd}, we characterize the expectation $E[R^*\left(\nabla \cL(\delta \theta^*; \mathfrak{X}_1^{n_1}, \mathfrak{X}_2^{n_2})\right)]$ in terms of the Gaussian width of the unit norm-ball of $R(.)$,  which leads to a sharp bound.
The upper bound on Gaussian width of the unit norm-ball of $R$ for atomic norms which covers a wide range of norms is provided in \cite{crpw12,chba15}.
\begin{theo}
Define $\Omega_R = \{u: R(u) \leq 1\}$. Let $\phi(R) = \sup_{{\bf{u}}}  \frac{\|{\bf{u}}\|_2}{R({\bf{u}})}$. Assume that for any ${\bf{u}}$ that $\|{\bf{u}}\| \leq \| \theta^*\|$
\begin{align}
\frac{1}{2}  \lambda_{\max} \left(\nabla^2 {\cL}(\delta \theta^*+{{\bf{u}}}) \right) \leq \eta_0,
\end{align}
where $\lambda_{\max}(.)$ is the maximum eigenvalue.
Then under the smooth density ratio assumption, we have
\begin{align*}
E\left[R^*(\nabla \cL(\delta \theta^*; \mathfrak{X}_1^{n_1}, \mathfrak{X}_2^{n_2}))\right] \leq \frac{2 \sqrt{\eta_0} (c_1 w\left(\Omega_R) + \phi(R) \right)}{\sqrt{\min(n_1,n_2)}}.
\end{align*}
and with probability at least $1-c_2 e^{-\epsilon^2}$
\begin{align*}
R^*\left(\nabla \cL(\delta \theta^*; \mathfrak{X}_1^{n_1}, \mathfrak{X}_2^{n_2}) \right) \leq  \frac{c_2(1+\epsilon) w(\Omega_R) + \tau_1}{\sqrt{\min(n_1, n_2)}} .
\end{align*}
where $c_1$ and $c_2$ are positive constants, $\tau_1 = 2 \sqrt{\eta_0} \phi(R)$, and $w(\Omega_R)$ is the Gaussian width of set $\Omega_R$.
 \label{theo:lambdabnd}
\end{theo}
Note, that our analysis hold for any norm and it is expressed in terms of the Gaussian width.
In the following, we give the bound on the regularization parameter for two examples of the regularization function $R(.)$.

\begin{corollary}
If $R(\delta \theta)$ is the $L_1$ norm, and $\delta \theta \in \mathbb{R}^{p^2}$ then with high probability we have the bound
\beq
  R^*\left(\nabla \cL(\delta \theta^*; \mathfrak{X}_1^{n_1}, \mathfrak{X}_2^{n_2})\right) \leq  \frac{\eta_2 \sqrt{\log p}}{\sqrt{\min(n_1,n_2)}}.
\label{eq:lambdaL1}
\eeq
\end{corollary}
%

\begin{corollary}
If $R(\delta \theta)$ is the group-sparse norm, and $\delta \theta \in \mathbb{R}^{p^2}$ then with high probability we  have the bound
\beq
  R^*\left(\nabla \cL(\delta \theta^*; \mathfrak{X}_1^{n_1}, \mathfrak{X}_2^{n_2})\right) \leq  \frac{\eta_2 \sqrt{m+\log N_G}}{\sqrt{\min(n_1,n_2)}},
\label{eq:lambdaGroup}
\eeq
where $\mathcal{G}=\{\mathcal{G}_1, \cdots, \mathcal{G}_{N_G} \}$ is a collection of groups, $m=\max_i |\mathcal{G}_i|$ is the maximum size of any group.
\end{corollary}

\subsection{ RSC Condition}
In this Section, we establish the RSC condition for direct change detection estimator \eqref{eq:changeestimateopt}.
Simplifying the expression and applying mean value theorem twice on the left side of RSC condition \eqref{eq:RSC}, for $\forall \gamma_i \in [0, 1]$, we have
\begin{align}
\delta \cL(\delta \theta^{*}, u) := \cL(\delta \theta^{*} + u) - \cL(\delta \theta^{*}) - \langle \nabla \cL(\delta \theta^{*}), u \rangle  \geq u^T \nabla^2 \cL(\delta \theta^*+\gamma_i u) u.
\end{align}
Thus, the RSC condition depends on the non-linear terms of loss function.
Recall that the nonlinear term, second term, in Loss function \eqref{eq:changeestimateopt} 
which is the approximation of the log-partition functions only depends on $n_2$. As a results, only samples of $\mathfrak{X}_2^{n_2}$ affect the RSC conditions. Our analysis is an extension of the results on \cite{bcfs14} using the generic chaining. We show that, with high probability the RSC condition is satisfied once samples $n_2$ crosses $w^2(C_r \cap S^{d-1})$ the Gaussian width of restricted error set. The bound on Gaussian width of the error set for atomic norms has been provided in \cite{chba15}.

Let $r_i = r(X={\bf{x}_i^2} | \delta \theta^*)$ and $\bar{\varepsilon}$ denote the probability that $r_i$ exceeds some constant $\eta_0$: $\bar{\varepsilon} = p(r_i > \eta_0) \geq 1-  e^{-\frac{\eta_0^2}{2}}$.

\begin{theo}
Let $X \in \R^{n \times p}$ be a design matrix with independent isotropic sub-Gaussian rows with $\vertiii{X_i}_{\Psi_2} \leq \kappa$.
Then, for any set $A \subseteq S^{p-1}$, for suitable constants $\eta$, $c_1$, $c_2 >0$ with probability at least $
1-\exp(-\eta w^2(A))$, we have
\begin{align}
\inf_{u \in A} \partial \cL(\theta^*;u,X) \geq c_1 \underline{\rho}^2 \left ( 1 - c_2\kappa_1^2 \frac{w(A)}{\sqrt{n_2}} \right) - \tau
\end{align}
where $\kappa_1=\frac{\kappa}{\bar{\varepsilon}}$, $\underline{\rho}^2 = \inf_{{\bf{u}} \in A} \rho_{\bf{u}}^2$ with
$\rho_{\bf{u}}^2 = E\left[\left\langle {\bf{u}}, T(X_i^2) \right\rangle^2 \mathbb{I} (r_i > \eta_0)\right]$, and $\tau$ is smaller than the first term in right hand side. Thus, for $n_2 \geq c_2 w^2(A)$, with probability at least $1- \exp(-\eta w^2(A))$, we have
$\inf_{u \in A} \partial \cL(\theta^*;u,X) > 0$.
\label{theo:rsc}
\end{theo}

\subsection{Statistical Recovery}
With the above results in place, from \eqref{eq:recover}, Theorem \ref{theo:recovBnd} provides the main recovery bound for generalized direct change estimator \eqref{eq:changeestimateopt}.

\begin{theo}
Consider two set of i.i.d samples $ \mathfrak{X}_1^{n_1} = \{{\bf{x}}^1_i\}_{i=1}^{n_1}$ and $\mathfrak{X}_2^{n_2}=\{{\bf{x}}^2_i\}_{i=1}^{n_2}$.
Define $\Omega_R = \{u : R(u) \leq 1\}$. Assume that $\delta \hat{\theta}$ is the minimizer of the
problem \eqref{eq:changeestimateopt}.
Then, with probability at least $1-\eta_0 e^{-\epsilon^2}$ the followings hold
\beq
 \lambda_{n_1,n_2} \geq \frac{\eta_1}{\sqrt{\min(n_1, n_2)}} (w(\Omega_R)+\epsilon)
\eeq
and for $n_2 \geq cw^2(C_r \cap S^{d-1})$,  with high probability, the estimate $\delta \hat{\theta}$ satisfies
\beq
\| \Delta \|_2 \leq O \left(  \frac{w(\Omega_R)}{\sqrt{\min(n_1, n_2)}}  \right) \Psi(C_r)~,
\eeq
where $w(.)$ is the Gaussian width of a set,and $c_2$, $\eta_0$, and $\eta_1$ are positive constants.
\label{theo:recovBnd}
\end{theo}

\proof
Proof of the Theorem can be directly obtain as the results of \eqref{eq:recover} and Theorem \ref{theo:lambdabnd} and Theorem \ref{theo:rsc}.

In the following, we provide the recovery bound for two special cases as an example.


\begin{corollary}
If $R(\delta\theta)$ is the $L_1$ norm, $\delta\theta^* \in \mathbb{R}^{p^2}$ s $s$-sparse., $\Psi(C_r) \leq  4\sqrt{s}$, and for $n_2 > cs \log p$, the recovery error is bounded by
\begin{align*}
\| \Delta \|_2 \leq c_3 \frac{\Psi(C_r) \lambda_{n_1,n_2} }{\kappa} = O \left ( \sqrt{\frac{s \log p}{\min(n_1,n_2)}} \right).
\end{align*}
\end{corollary}
\begin{corollary}
If $R(\delta\theta)$ is the group-sparse norm, $\delta\theta \in \mathbb{R}^{p^2}$, $\Psi(C_r) \leq 4\sqrt{s_G}$ and for $n_2 \geq c(ms_G + s_G \log N_G)$, the recovery error  is bounded by
\begin{align*}
\| {\Delta} \|_2 & \leq c_3 \frac{\Psi(C_r) \lambda_{n_1,n_2} }{\kappa} = O \left( \sqrt{\frac{s_G m + \log N_G}{\min(n_1,n_2)}} \right).
\end{align*}
\end{corollary}




\section{Experiments}
\label{sec:res}
%

In this Section, we evaluate generalized direct change estimator (direct) with three different norms.
and we compare our direct approach with indirect approach. For indirect approach, we first estimate Ising model structures $\hat{\theta}_1$ and $\hat{\theta}_2$ with $L_1$ norm regularizer, separately~\cite{rawl10}. Then, we obtain $\delta \hat{\theta} = \hat{\theta}_1 - \hat{\theta}_2$. In all experiments, we draw $n_1$ and $n_2$ i.i.d samples from each Ising model by running Gibbs sampling. Here we set $n=n_1=n_2=\{20, 50, 100\}$. 

\begin{figure*}
\centering
\subfigure[$\theta_1$]{\includegraphics[trim = 10mm 0mm 15mm 0mm, clip, width = 0.24\textwidth]{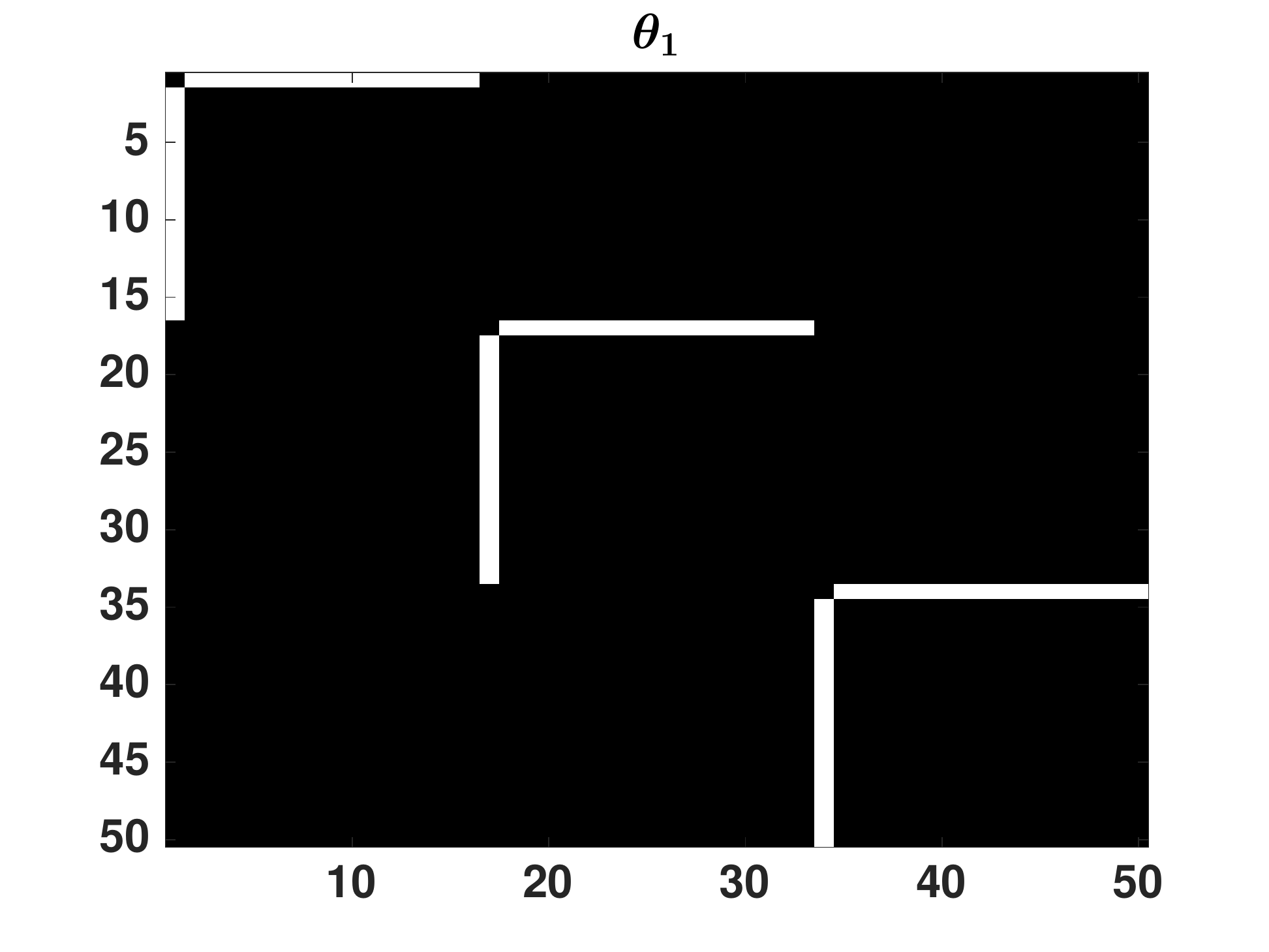}}
\subfigure[$\theta_2$]{\includegraphics[trim = 10mm 0mm 15mm 0mm, clip,  width = 0.24\textwidth]{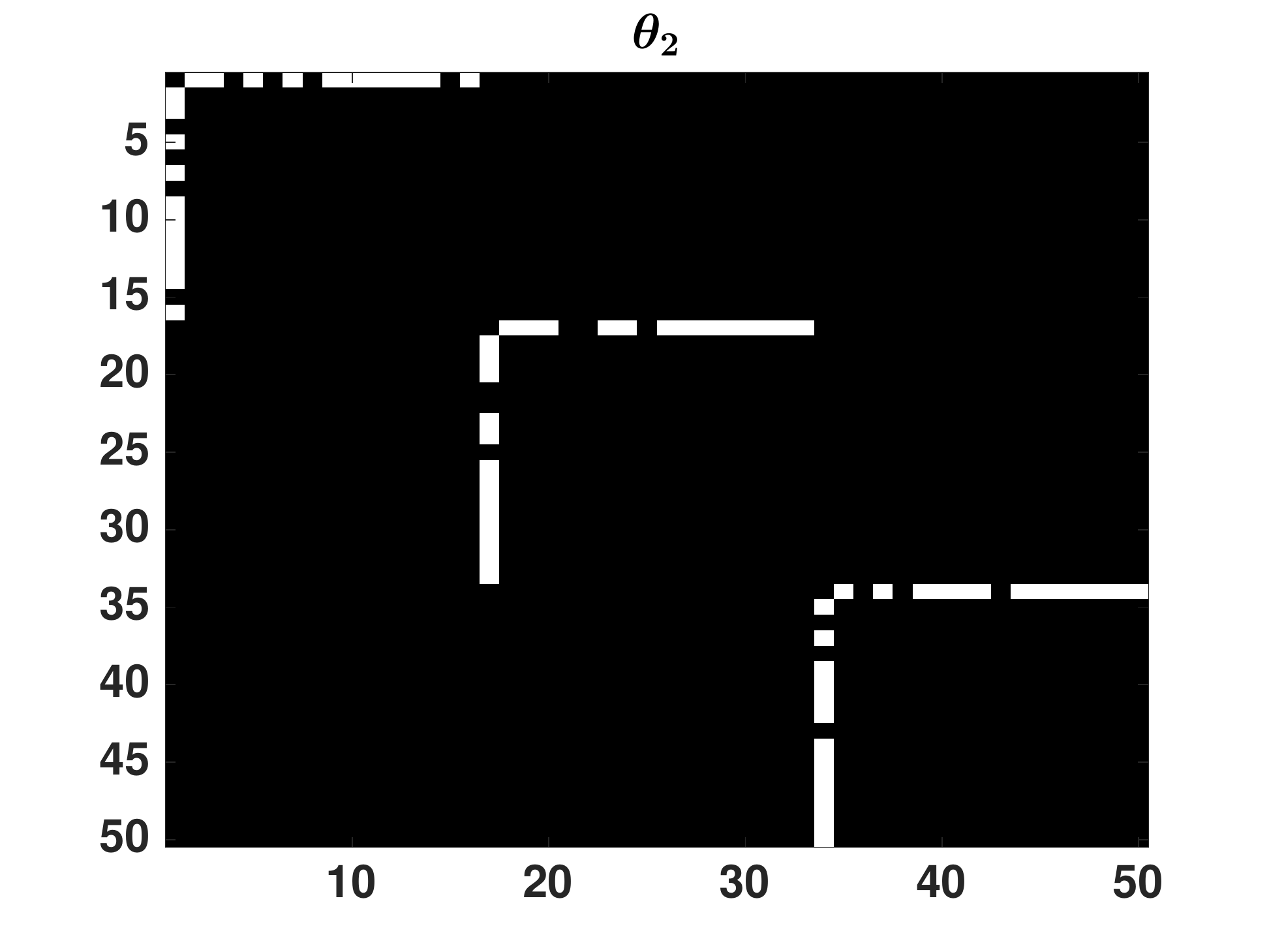}}
\subfigure[$\delta \theta = \theta_1-\theta_2$]{\includegraphics[trim = 10mm 0mm 15mm 0mm, clip,  width = 0.24\textwidth]{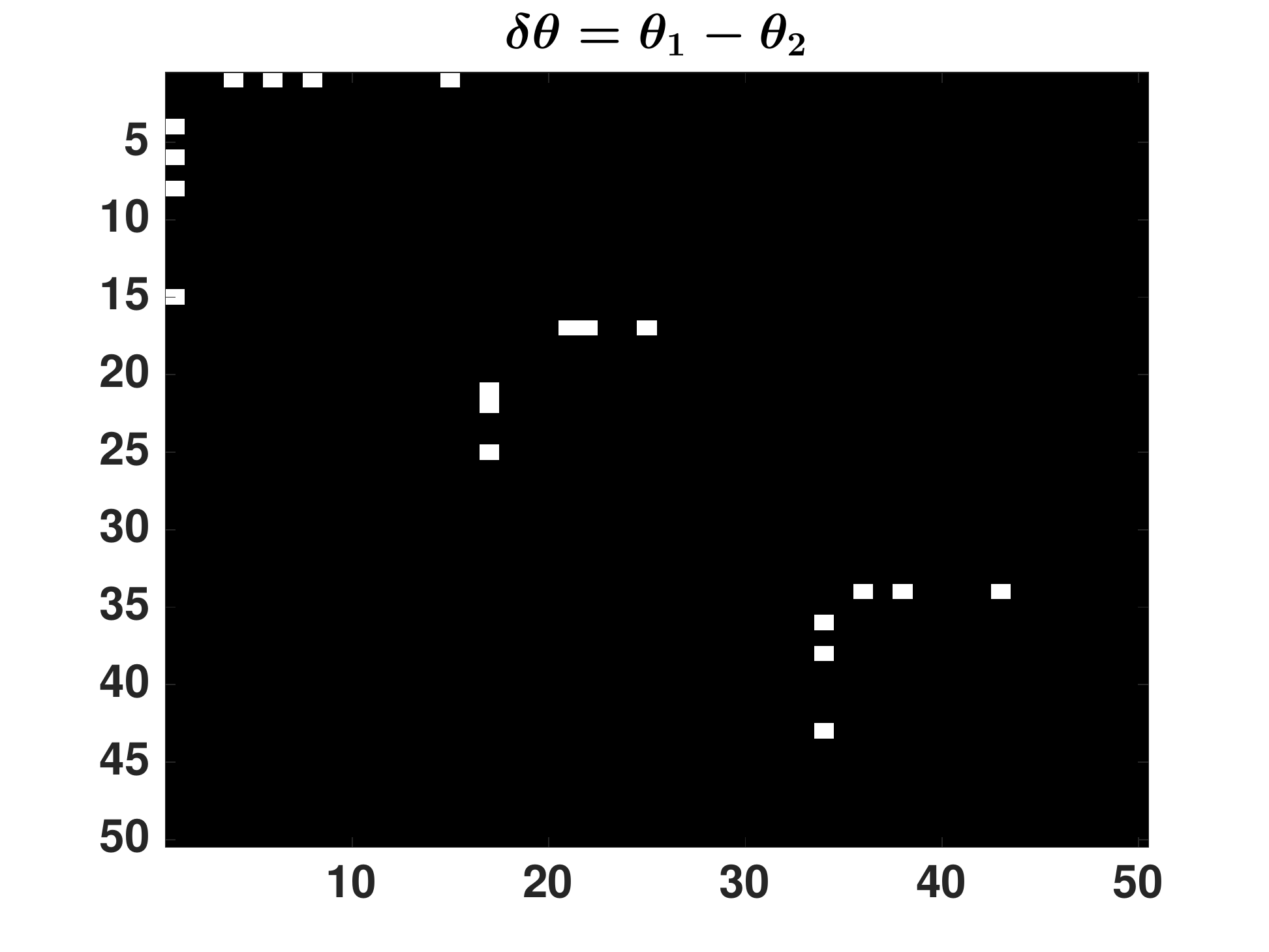}}
\subfigure[ROC]{\includegraphics[trim = 5mm 0mm 15mm 10mm, clip,  width = 0.245\textwidth]{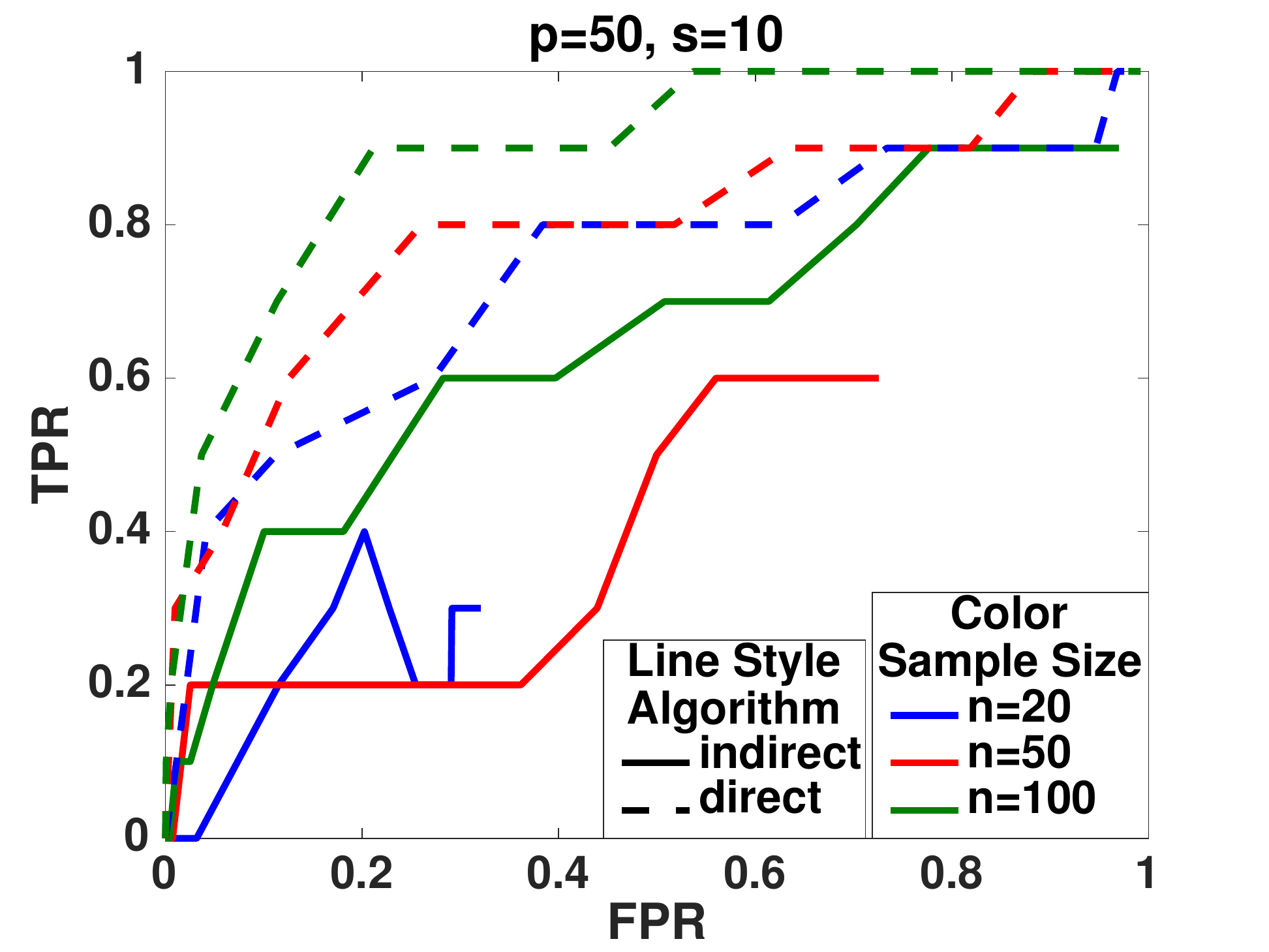}}

\vspace{-0.3cm}
\subfigure[$\theta_1$]{\includegraphics[trim = 10mm 0mm 15mm 0mm, clip, width = 0.24\textwidth]{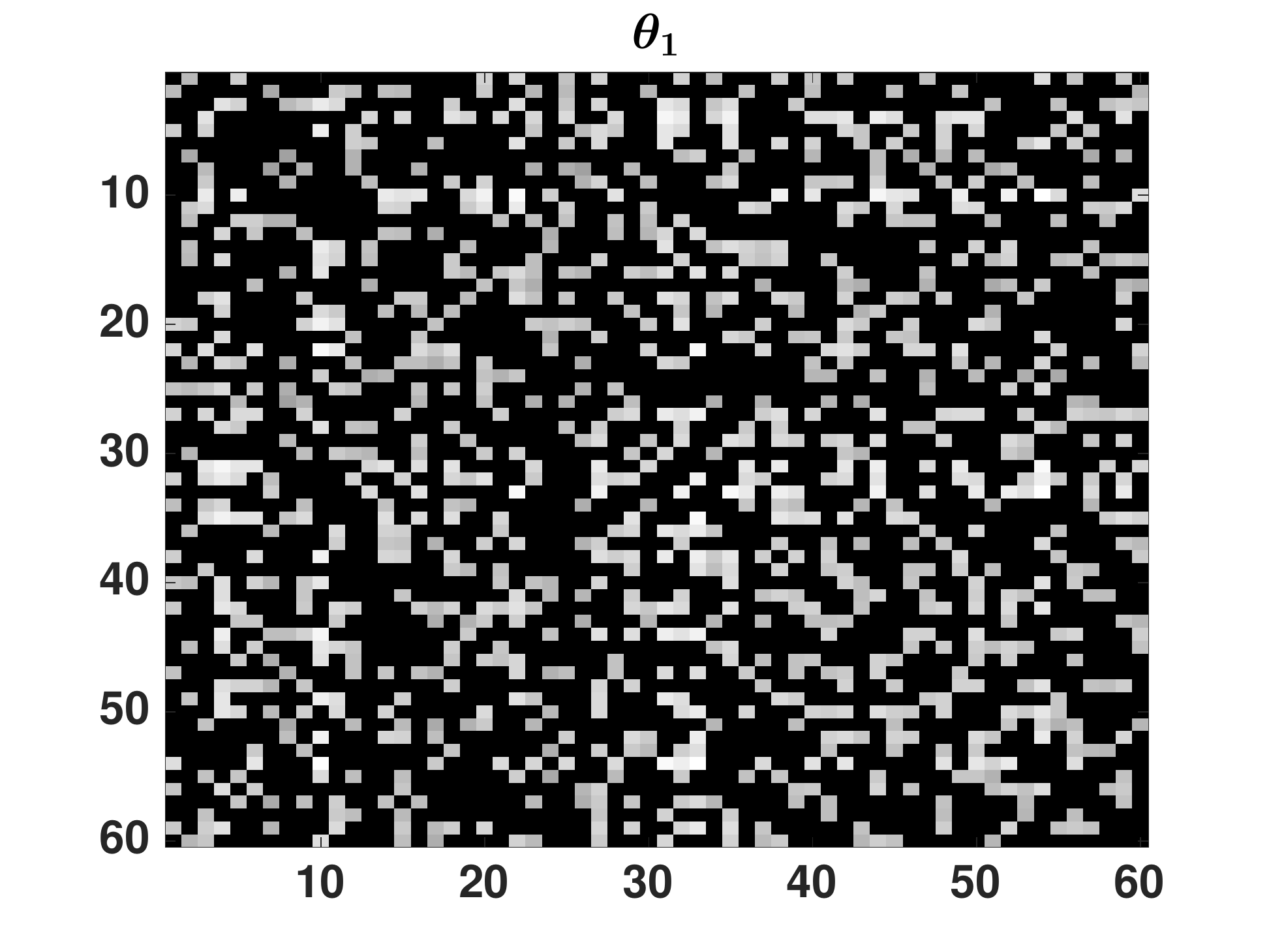}}
\subfigure[$\theta_2$]{\includegraphics[trim = 10mm 0mm 15mm 0mm, clip, width = 0.24\textwidth]{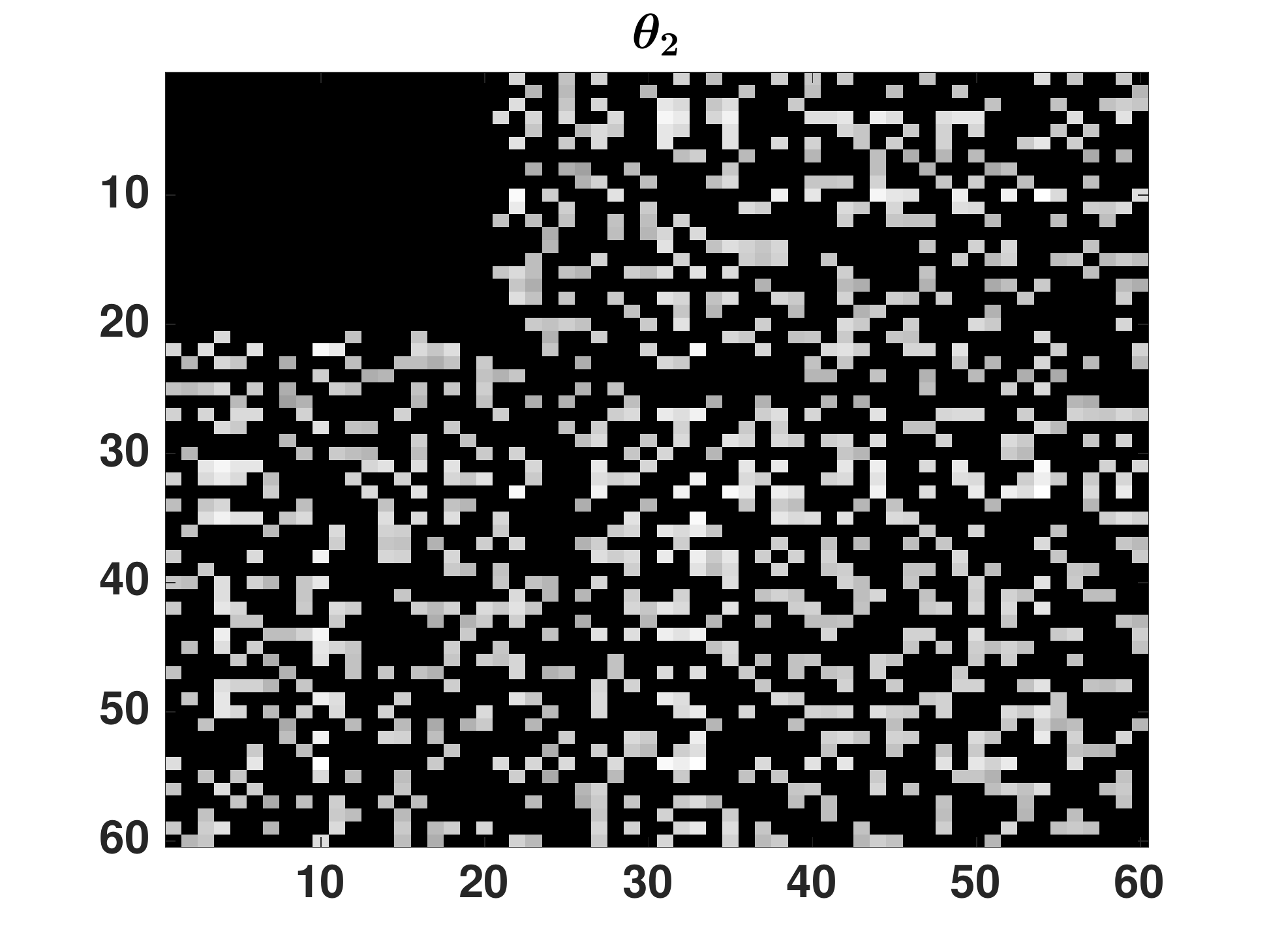}}
\subfigure[$\delta \theta = \theta_1-\theta_2$]{\includegraphics[trim = 10mm 0mm 15mm 0mm, clip, width = 0.24\textwidth]{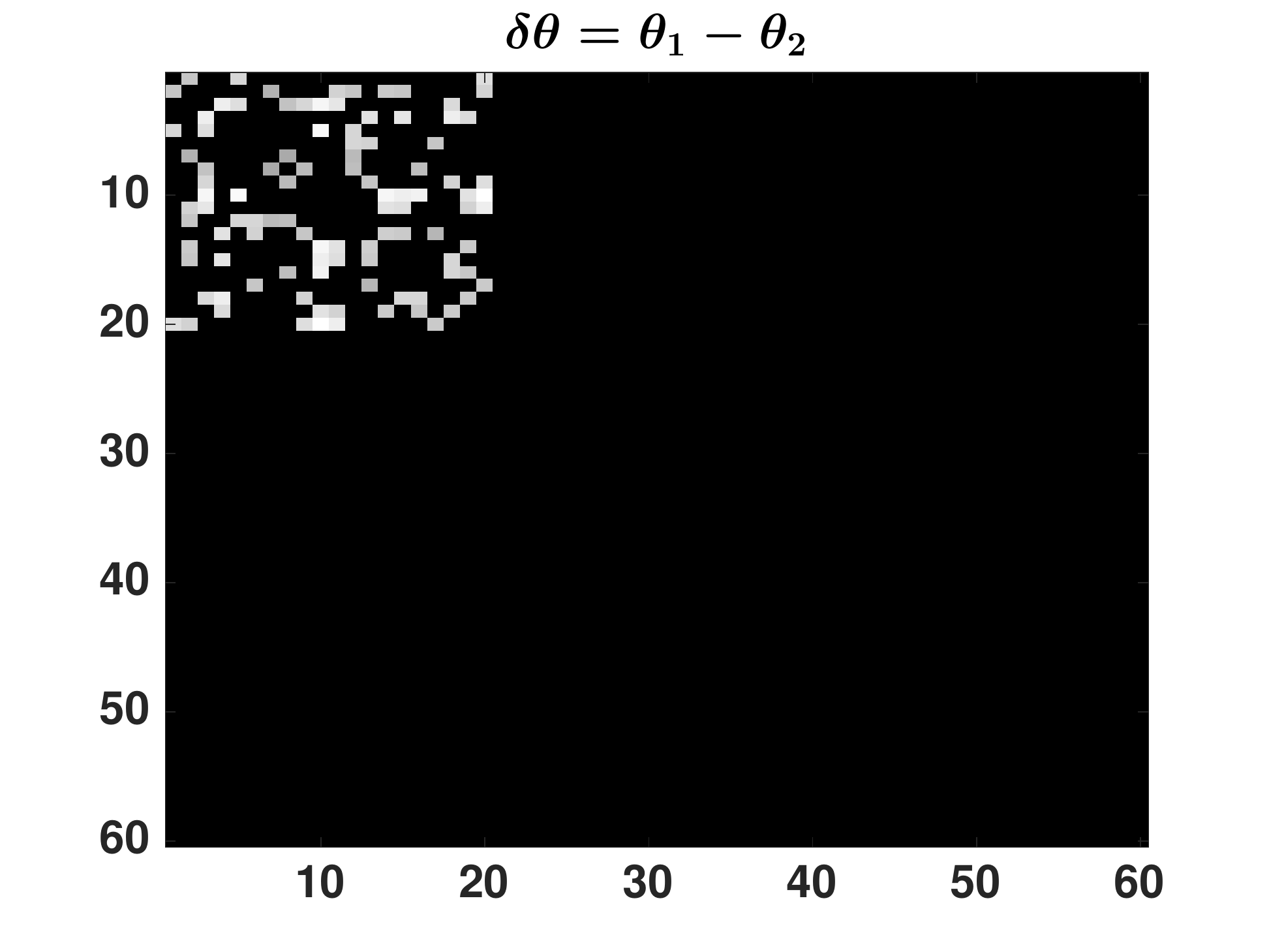}}
\subfigure[ROC]{\includegraphics[trim = 5mm 0mm 15mm 10mm, clip, width = 0.245\textwidth]{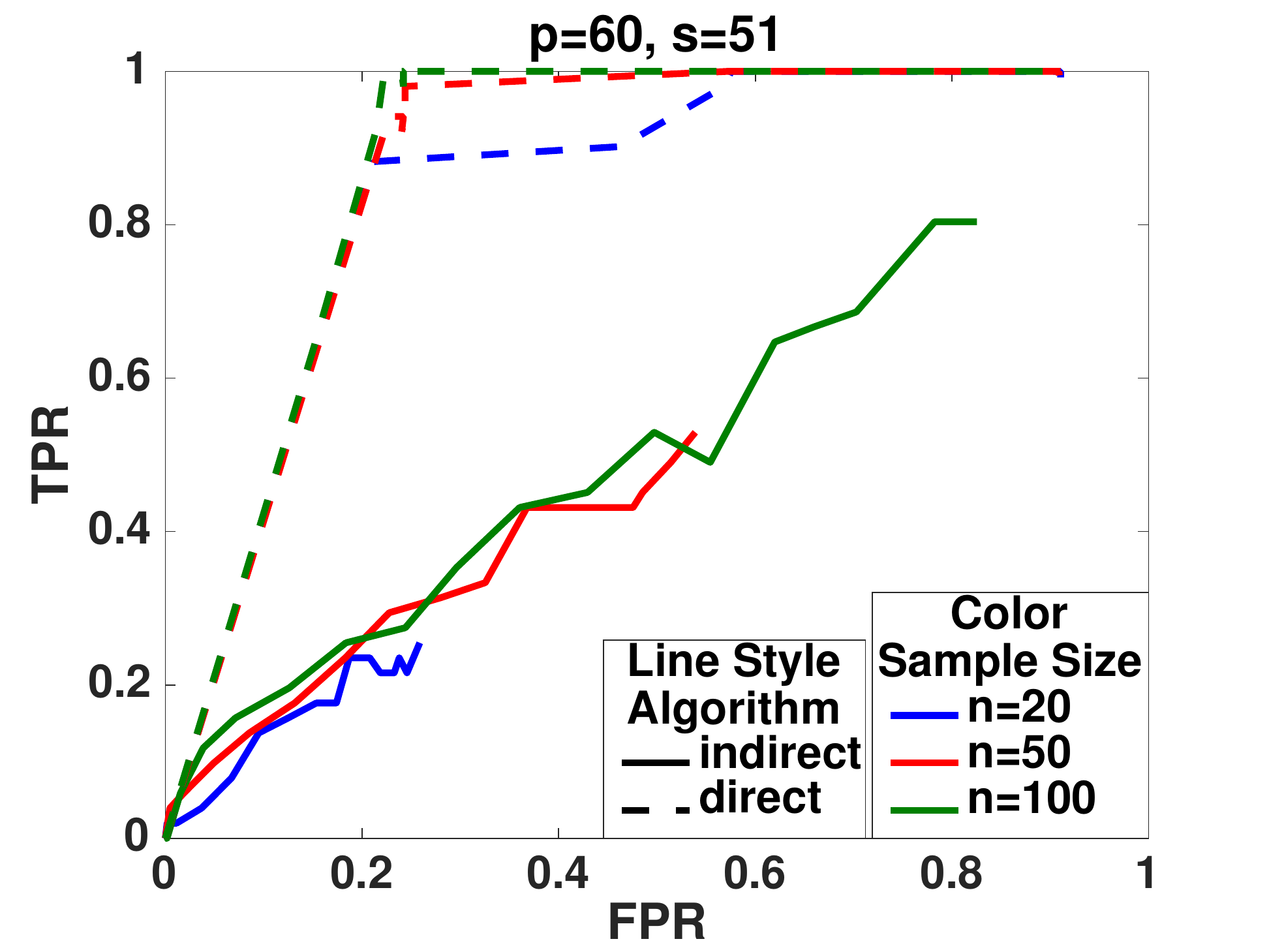}}

\vspace{-0.3cm}
\subfigure[$\theta_1$, scale-free]{\includegraphics[trim = 10mm 0mm 15mm 0mm, clip, width = 0.24\textwidth]{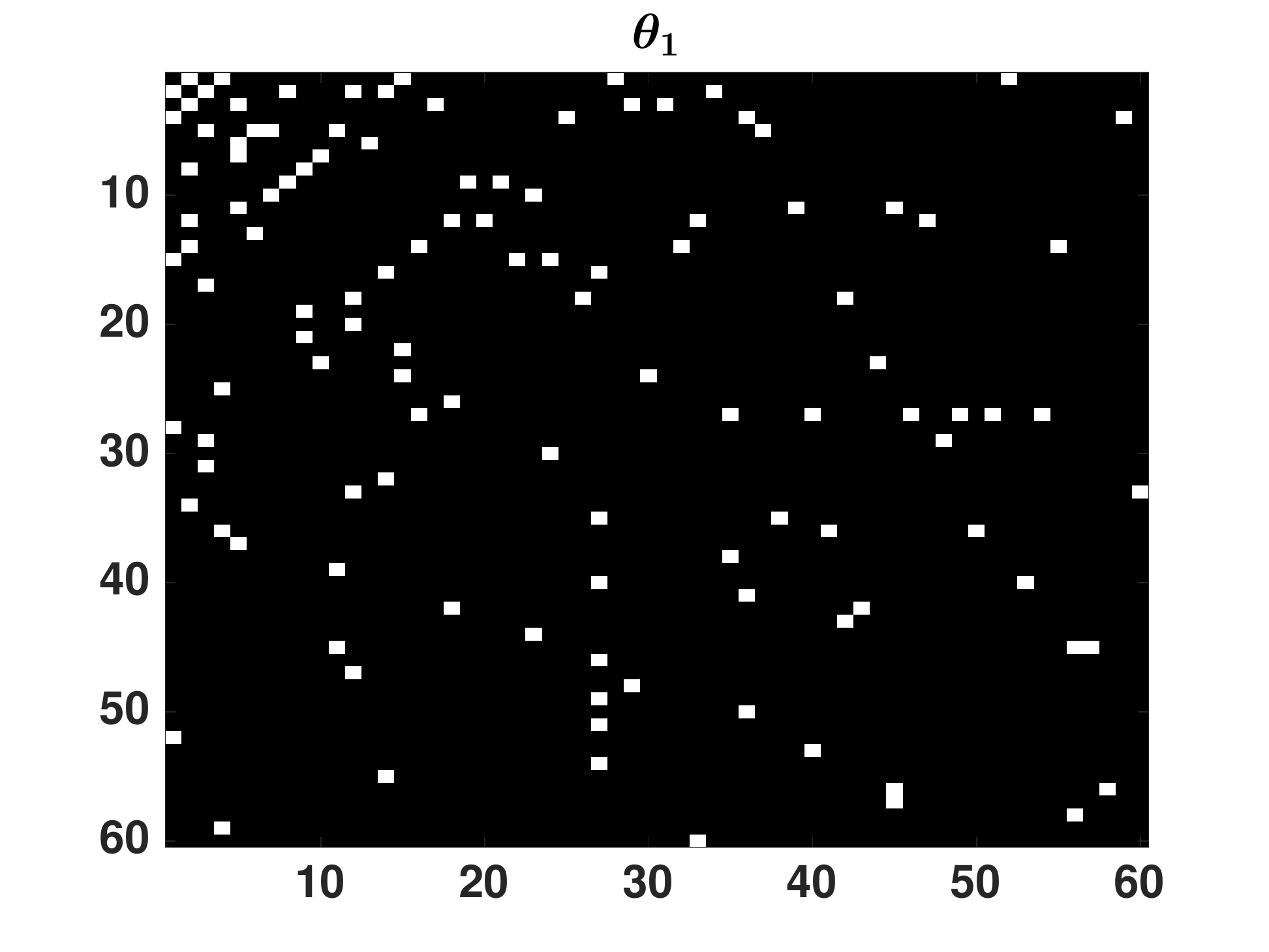}}
\subfigure[$\theta_2$]{\includegraphics[trim = 10mm 0mm 15mm 0mm, clip, width = 0.24\textwidth]{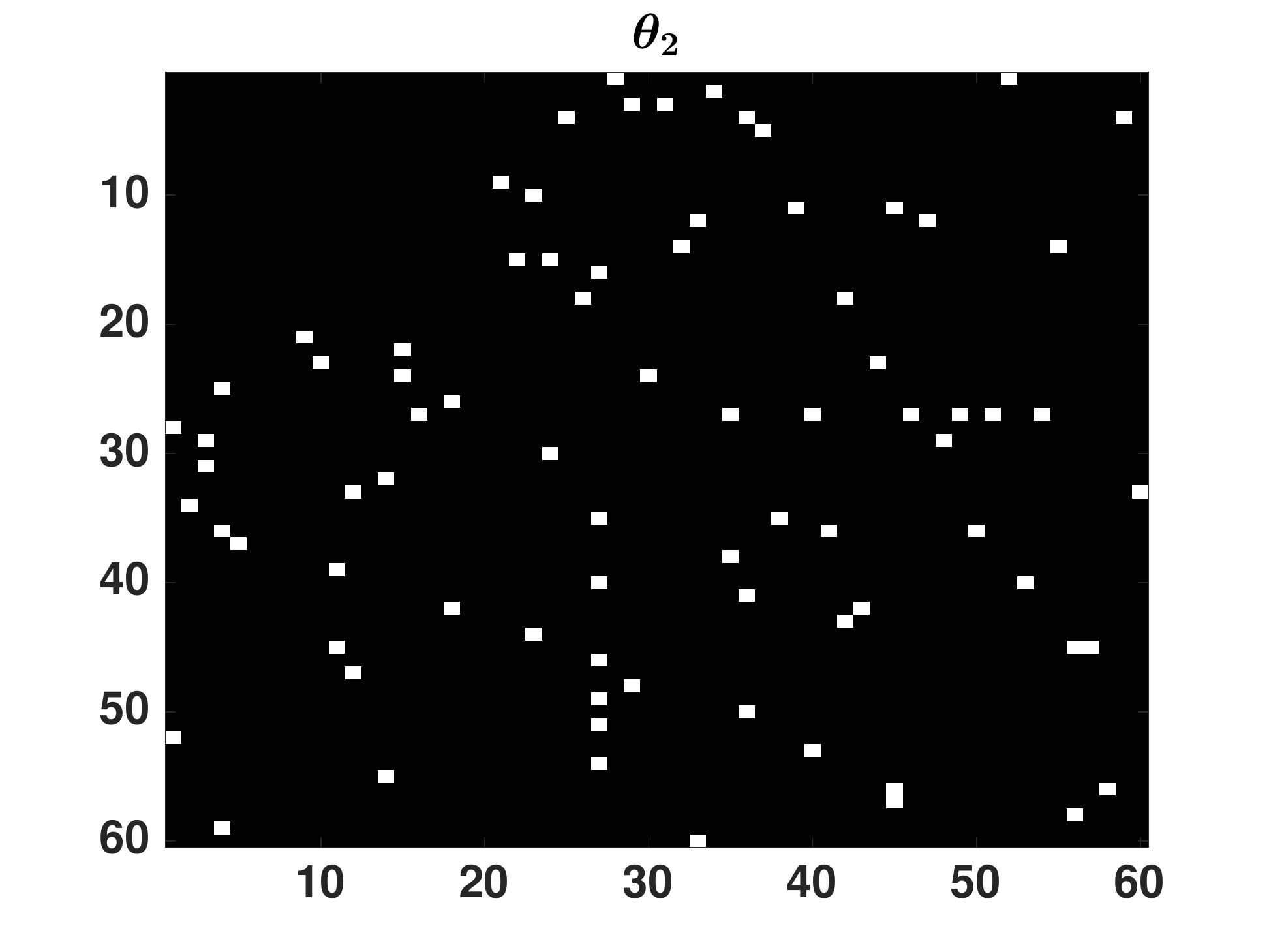}}
\subfigure[$\delta \theta = \theta_1-\theta_2$]{\includegraphics[trim = 10mm 0mm 15mm 0mm, clip, width = 0.24\textwidth]{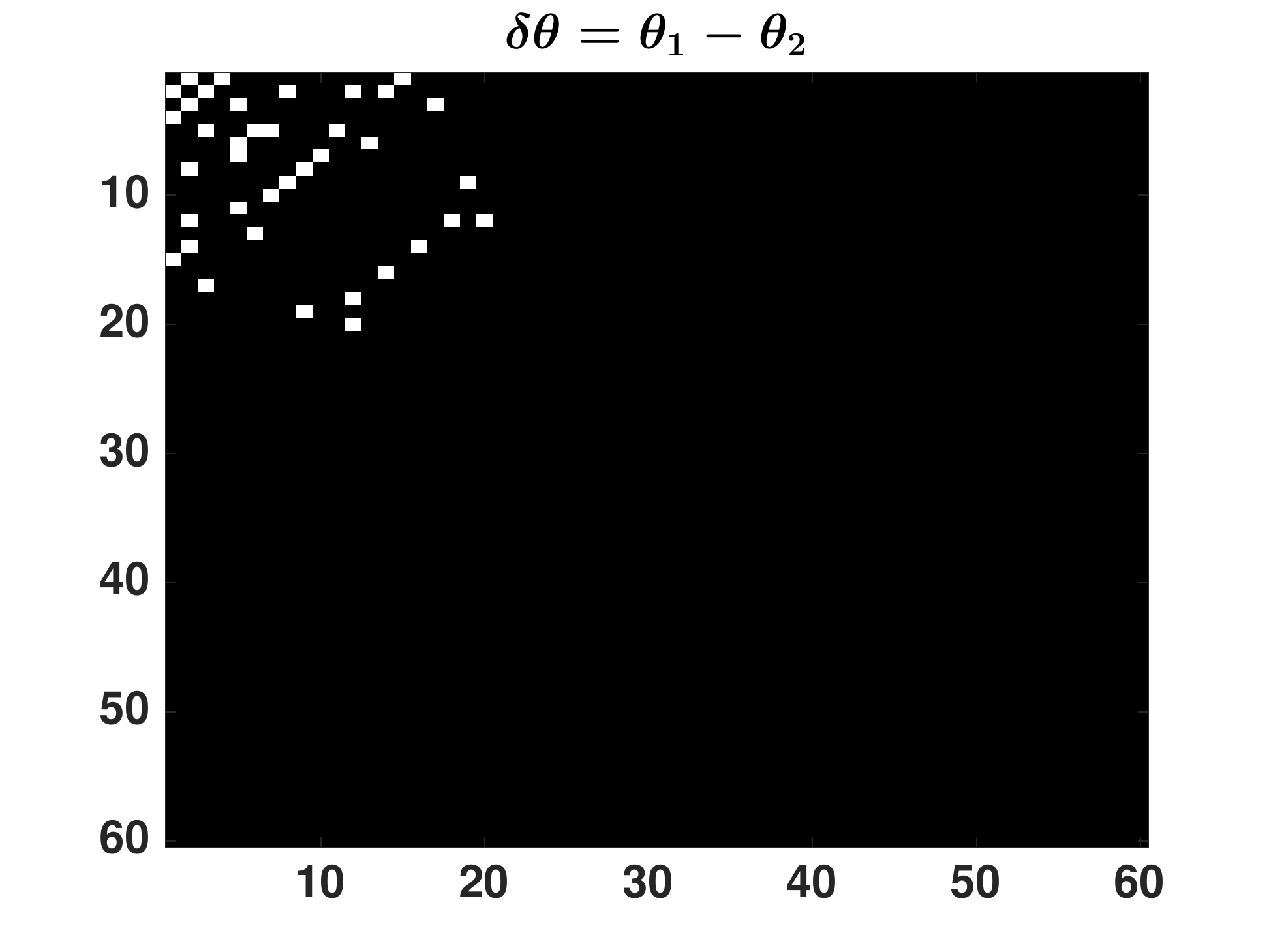}}
\subfigure[ROC]{\includegraphics[trim = 5mm 0mm 15mm 10mm, clip, width = 0.245\textwidth]{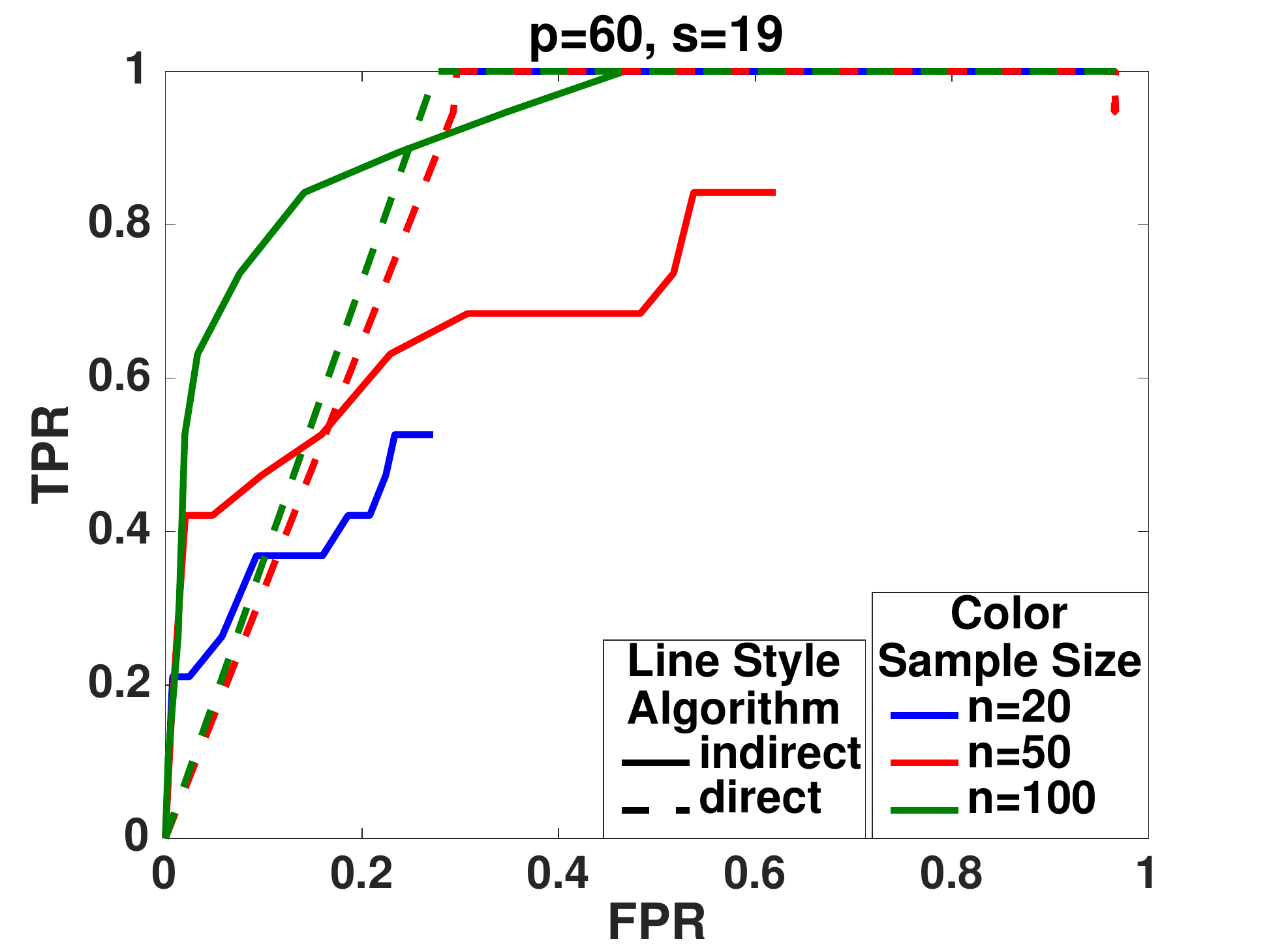}}

\vspace{-0.3cm}
\subfigure[$\theta_1$]{\includegraphics[trim = 10mm 0mm 15mm 0mm, clip, width = 0.24\textwidth]{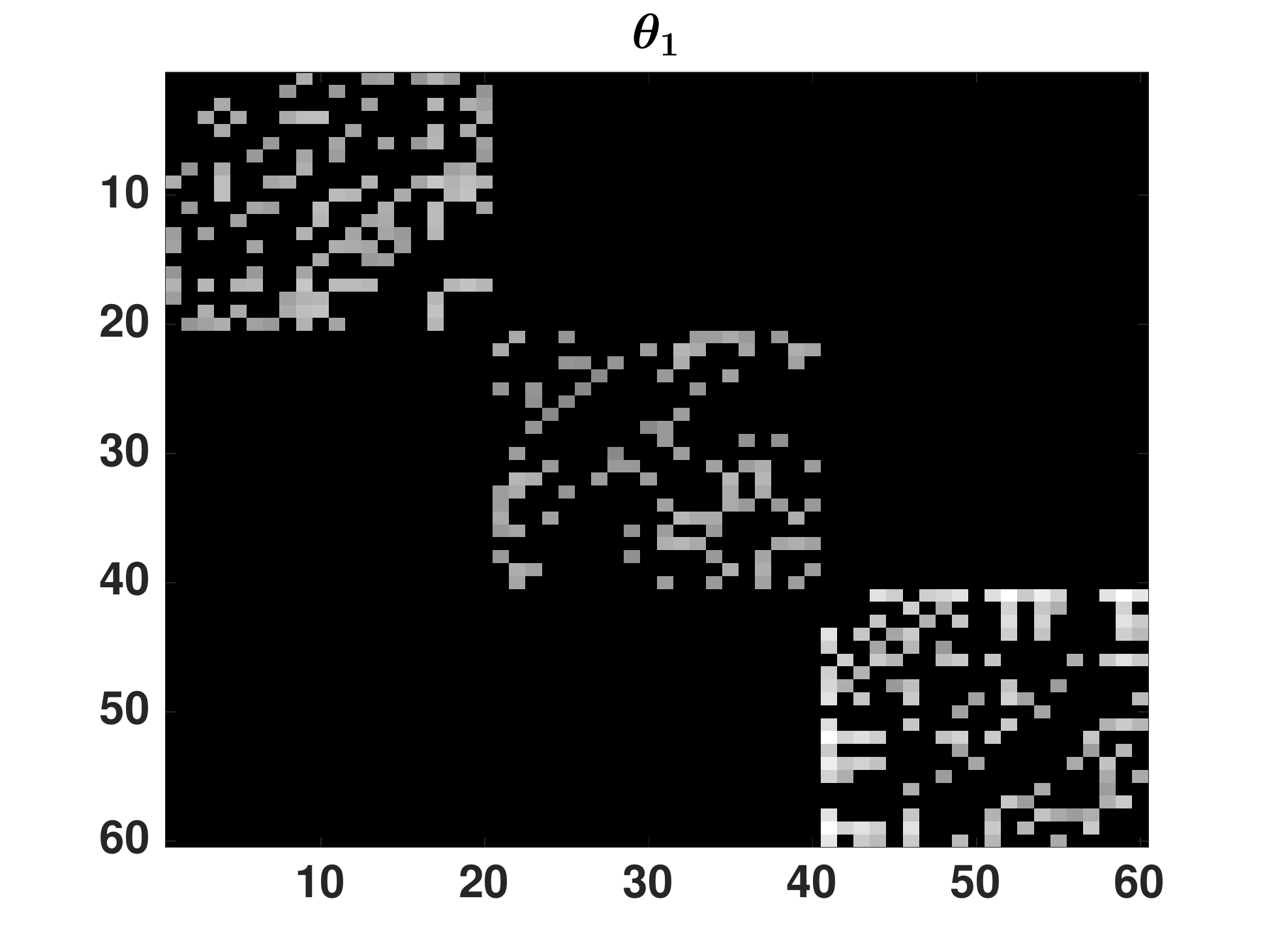}}
\subfigure[$\theta_2$]{\includegraphics[trim = 10mm 0mm 15mm 0mm, clip, width = 0.24\textwidth]{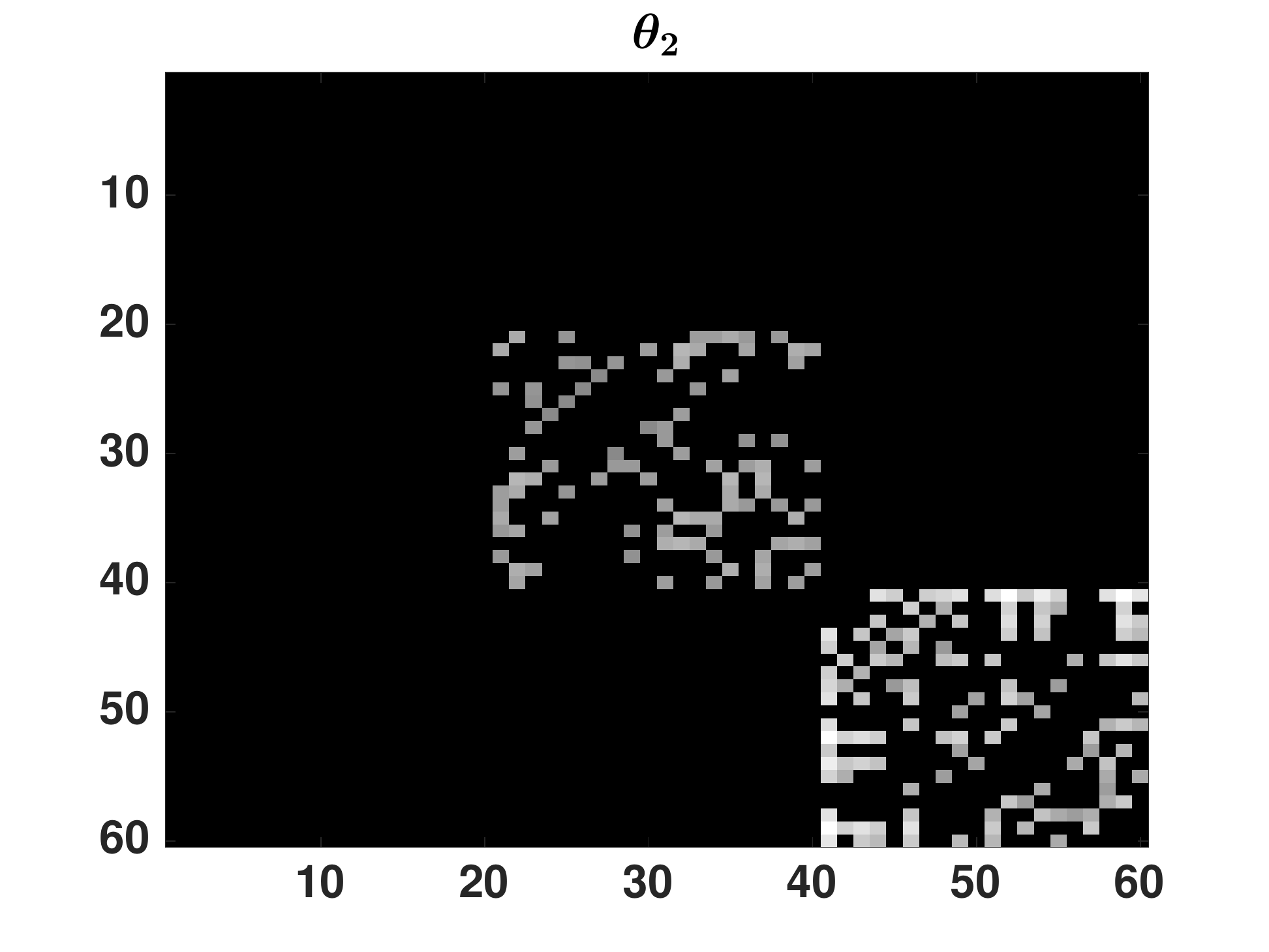}}
\subfigure[$\delta \theta = \theta_1-\theta_2$]{\includegraphics[trim = 10mm 0mm 15mm 0mm, clip, width = 0.24\textwidth]{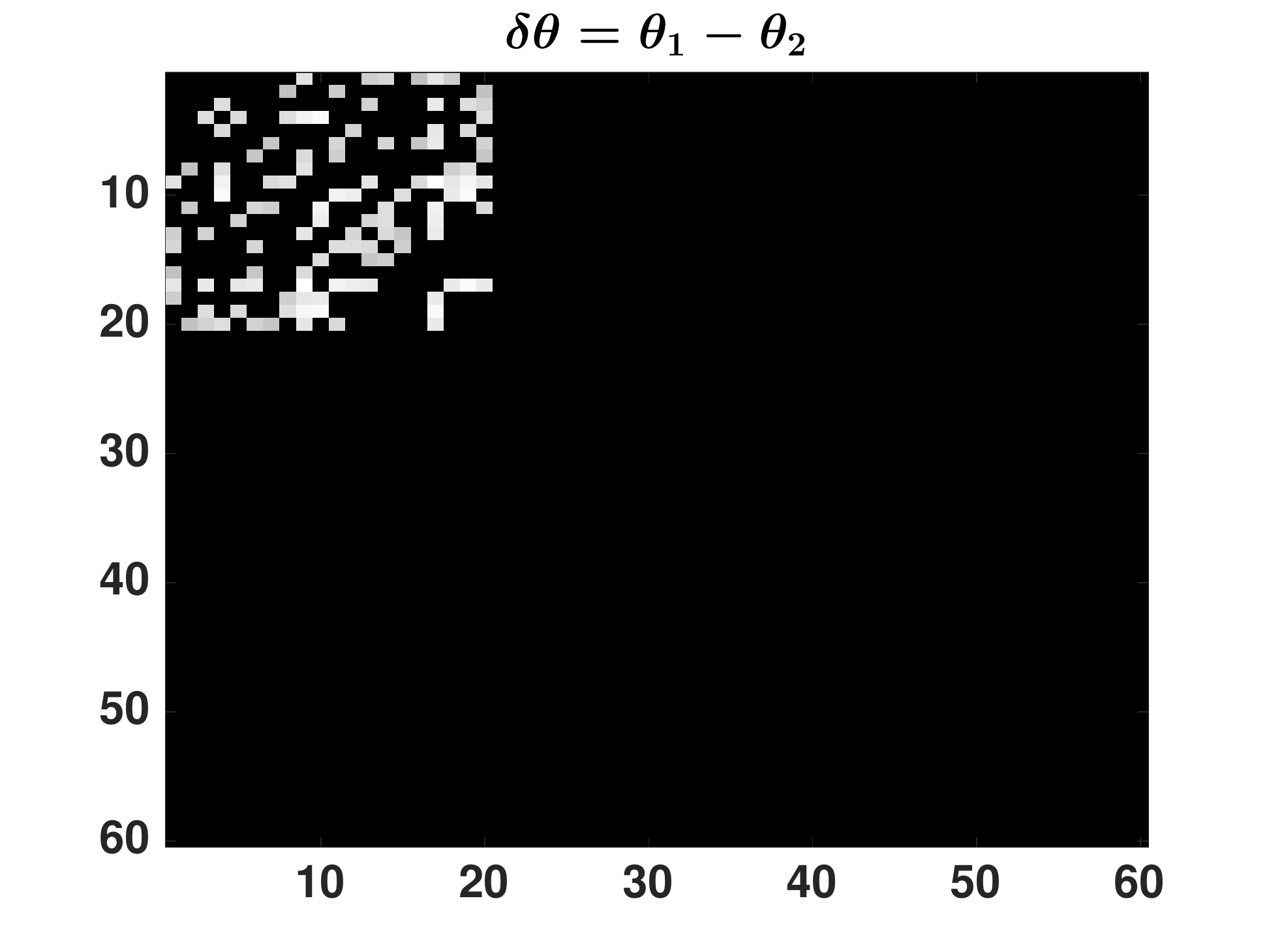}}
\subfigure[ROC]{\includegraphics[trim = 5mm 0mm 15mm 10mm, clip, width = 0.245\textwidth]{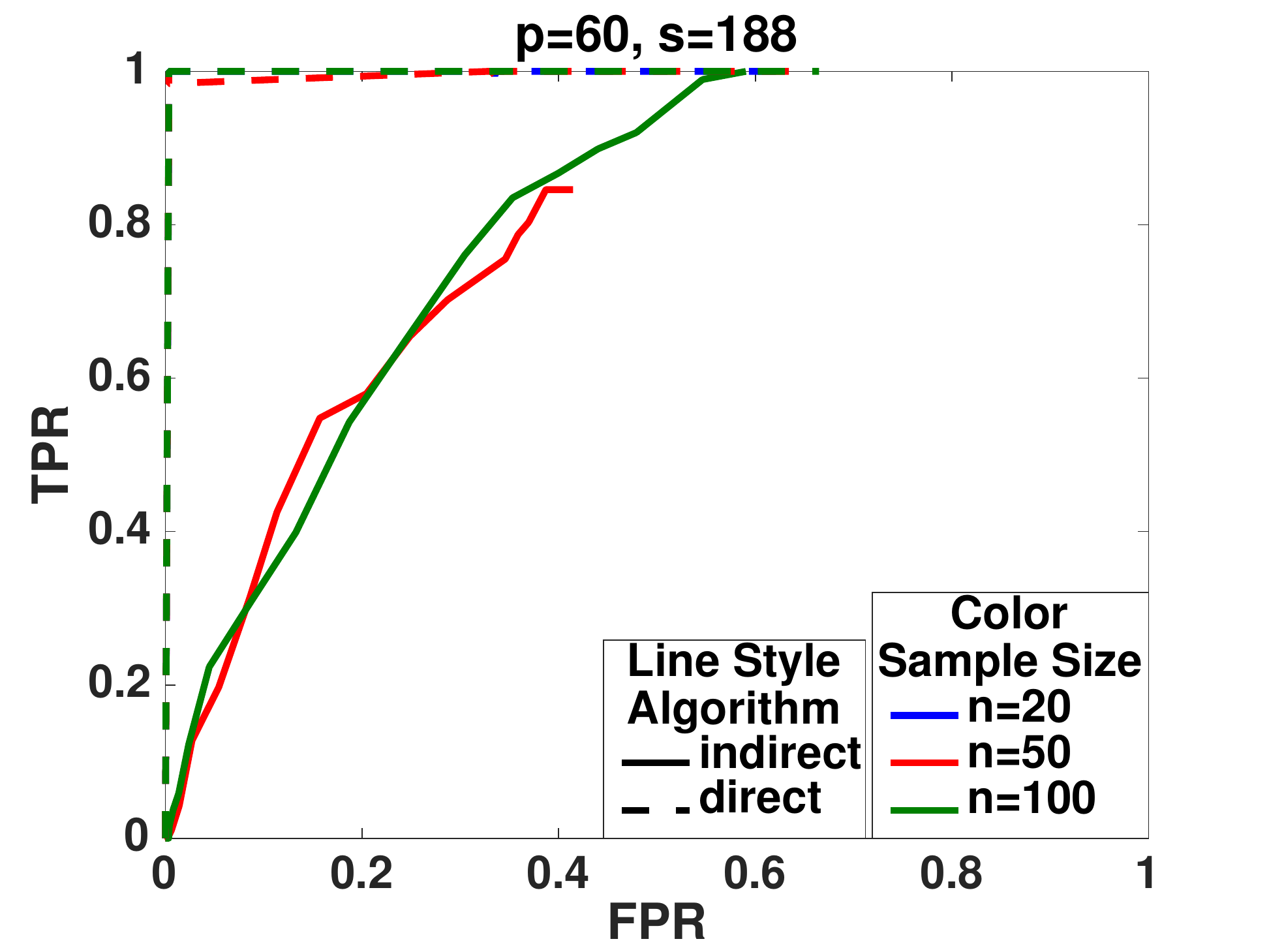}}

\vspace{-0.3cm}
\subfigure[$\theta_1$]{\includegraphics[trim = 10mm 0mm 15mm 0mm, clip, width = 0.24\textwidth]{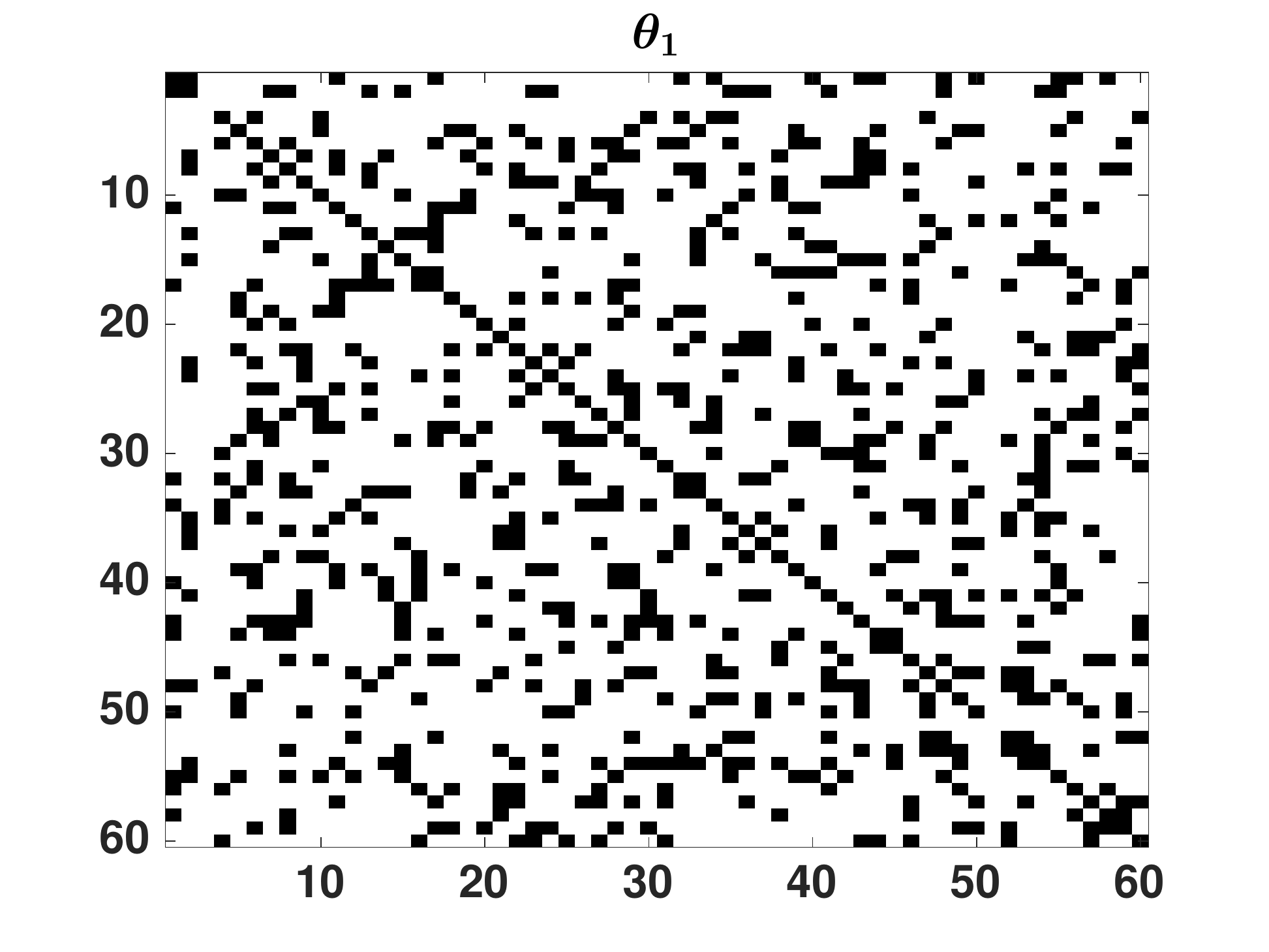}}
\subfigure[$\theta_2$]{\includegraphics[trim = 10mm 0mm 15mm 0mm, clip, width = 0.24\textwidth]{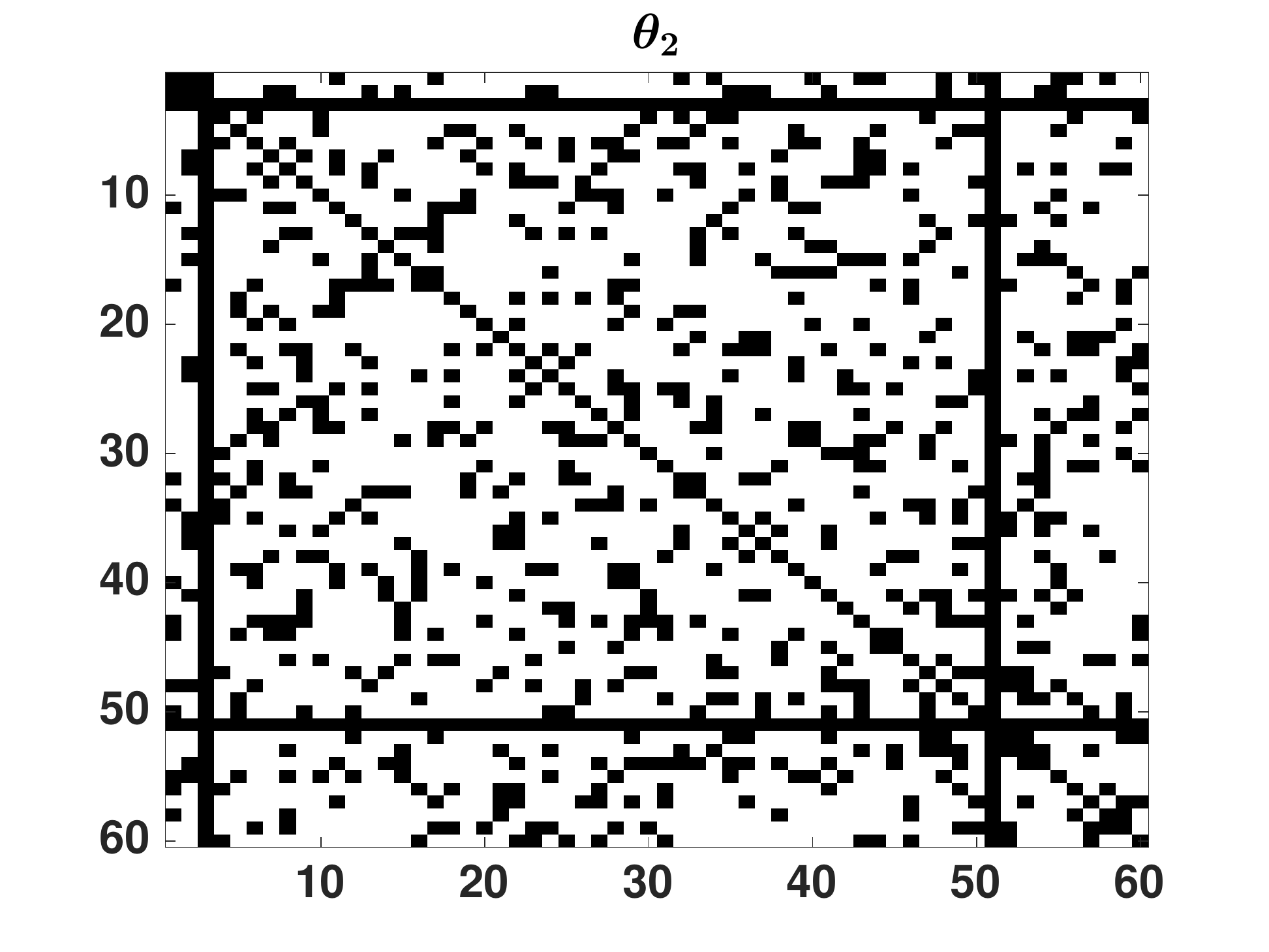}}
\subfigure[$\delta \theta = \theta_1-\theta_2$]{\includegraphics[trim = 10mm 0mm 15mm 0mm, clip, width = 0.24\textwidth]{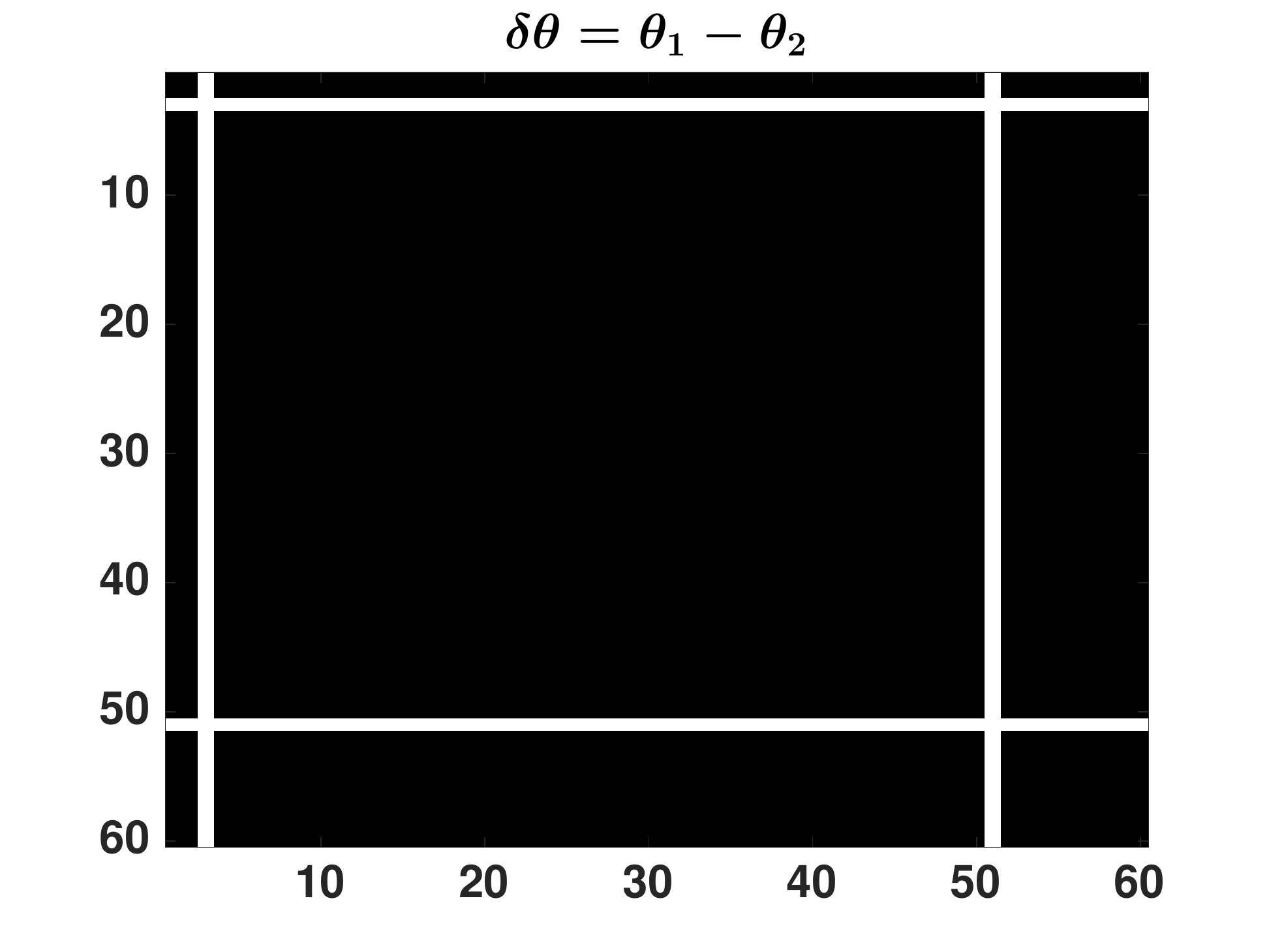}}
\subfigure[ROC]{\includegraphics[trim = 5mm 0mm 15mm 10mm, clip, width = 0.245\textwidth]{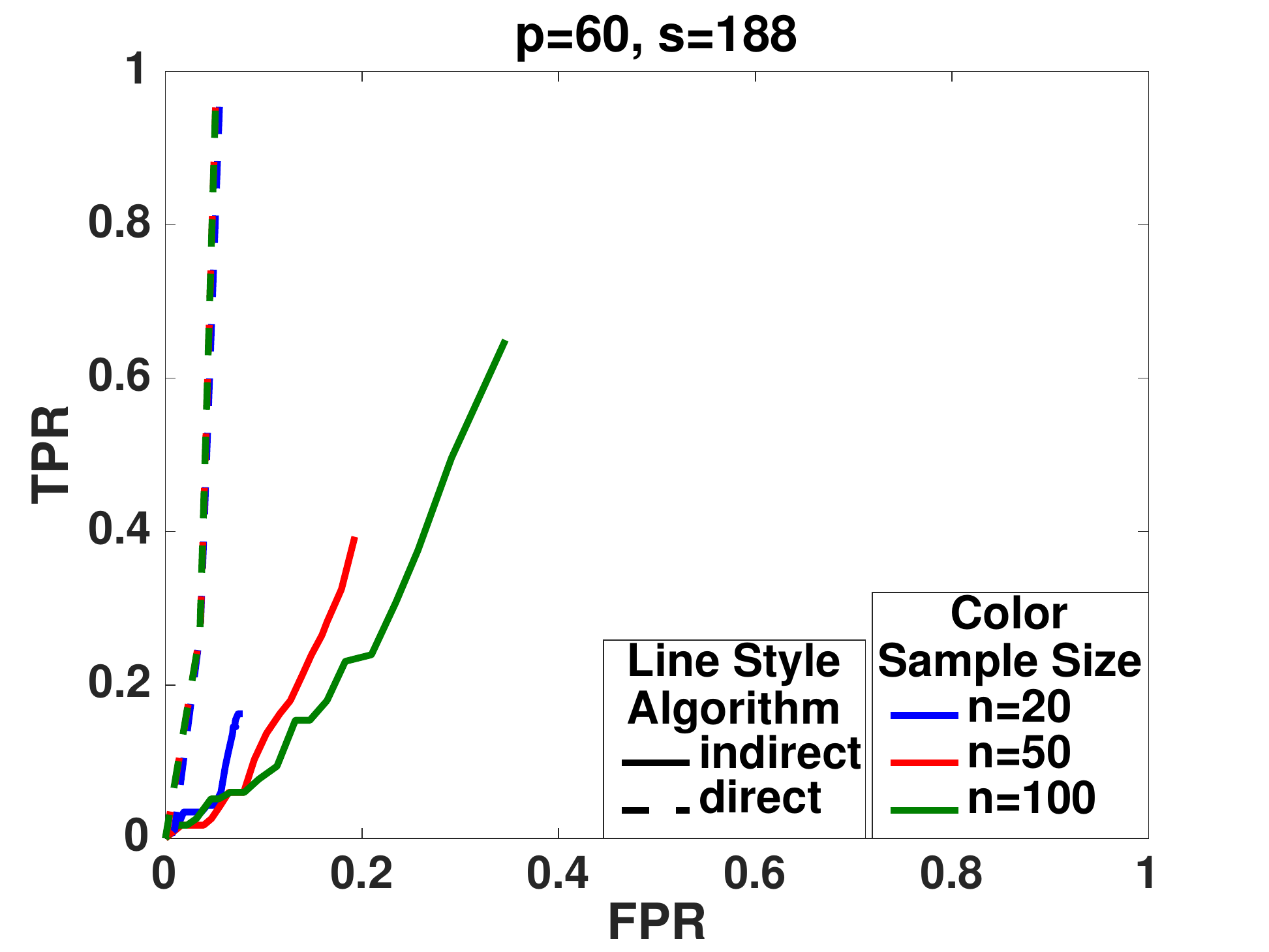}}
\label{fig:ROC}
\caption{First row $\delta\theta^*$ has a sparse structure ($L_1$ norm) and $\theta_1^*$ has 3 disconnected star graphs. Second, third, and forth rows $\delta\theta^*$ has group sparse structure (group sparse norm) where $\theta_1^*$ has a random graph structure in second row, scale-free structure in third row, and block structure in forth row. Last row  $\delta\theta^*$ has two perturbed norm (Node perturbation) and $\theta_1^*$ has a random graph structure. Blacks in heatmaps denotes zeros. ROC curve for different structures show in the last column. Direct approach has a better ROC curve for all structures except with scale-free structure of $\theta_1^*$.}
\end{figure*}

{\bf{$L_1$ norm:}} Here we first generate $\theta_1^*$ with three disconnected star sub-graphs (Figure \ref{fig:ROC}-a) with $p=50$. We generate the weights uniformly random between $\{0.3-0.5\}$. 
We then generate $\theta_2^*$ by removing 10 random edges from $\theta_1^*$ (Figure \ref{fig:ROC}-b).
It is interesting that although individual graphs are sparse, but direct approach has a better ROC curve for all values of $n$ (Figure \ref{fig:ROC}-d). Similar results obtained by with random graph structure of $\theta_1^*$ and $\theta_2^*$.

{\bf{Group-sparse norm:}} In this set of experiments, we evaluate direct method with three different structure for $\theta_1^*$: (i) a random graph structure (Figure \ref{fig:ROC}-e), (ii) scale free graph structure (Figure \ref{fig:ROC}-i), and (iii) a cluster graph structure (Figure \ref{fig:ROC}-m). In all settings, we set $p=60$ and generate $\theta_2^*$ by removing a block of edges from $\theta_1^*$ (Figure \ref{fig:ROC}-(f,j,n)). 
For random graph structure and block structure, direct method has a better ROC curve (Figure \ref{fig:ROC}-h,p). But, for scale-free structure, since the individual graphs are sparse, indirect method can estimate $\hat{\theta}_1$ and $\hat{\theta}_2$ correctly, and thus have a better ROC curve (Figure \ref{fig:ROC}-l).

{\bf{Node perturbation:}} Here, we first generate a random graph structure $\theta_1^*$, and then generate $\theta_2^*$ by perturbing two nodes in $\theta_1^*$. Here we set $p=60$ and 
generate $\theta_2^*$ by setting rows and columns $3, 51$ to zero in $\theta_1^*$ (Figure \ref{fig:ROC}-s). Although, the individual graphs are dense but direct approach can estimate edges in $\delta \theta$ with only $n=20$ samples (Figure \ref{fig:ROC}-t).

%
%

\section{Conclusion}
\label{sec:concl}
This paper presents the statistical analysis of direct change problem in Ising graphical models where any norm can be plugged in for characterizing the parameter structure. An optimization algorithm based on FISTA-style algorithms is proposed with the convergence rate of $O(1/T^2)$.
We provide the statistical analysis for any norm such as $L_1$ norm, group sparse norm, node perturbation, etc. Our analysis is based on generic chaining and illustrates the important role of Gaussian widths (a geometric measure of size of suitable sets) in such results.
For the special case of sparsity, we obtain a sharper result than previous results \cite{liss14} under the same smooth density ratio assumption. Liu et al. \cite{liss14} obtained the same result with a bounded density ratio assumption which is a more restrictive assumption. Although, we presented the results for Ising model, our analysis can be applied to any graphical model which satisfies the smooth density ratio assumption.
Further, we extensively compared our generalized direct change estimator with an indirect approach over a wide range of graph structures and norms. We show that our direct approach has a better ROC curve than indirect approach without any assumption on the structure of individual graphs.
We implemented indirect approach by estimating individual Ising model structures with $L_1$ norm regularizer. However, if individual graphs has a suitable structure such as group sparsity, one may apply a regularization that can characterize the graph structure and may improve performance of the indirect approach. We will investigate this possibility in our future research.

\section*{Appendix}

\appendix

\section{Background and Preliminaries}

\begin{defn}
{\bf Sub-Gaussian random variable:} We say that a random variable $x$ is sub-Gaussian if the moments satisfies
\beq
[E|x|^p]^{\frac{1}{p}} \leq K_2 \sqrt{p} 
\eeq
for any $p \geq 1$ with a constant $K_2$. The minimum value of $K_2$ is called sub-Gaussian norm of $x$, denoted by $\vertiii{x}_{\psi_2}$.
If $E[x]=0$, then
\beq
E[\exp\{tX\}] \leq \exp\{Ct^2 \vertiii{X}_{\Psi_2}^2\},
\eeq
where $C$ and $c$ are positive constant.
\label{def:subgvar}
\end{defn}

\begin{defn}
{\bf Sub-Gaussian random vector:} We say that a random vector $X$ in $\R^n$ is sub-Gaussian if the one-dimensional marginals $\langle X, {\bf{x}} \rangle$ are sub-Gaussian random variables for all ${\bf{x}} \in \R^n$. The sub-Gaussian norm of $X$ is defined as
\begin{align}
\vertiii{X}_{\psi_2} &= \sup\limits_{{\bf{x}} \in S^{n-1}} \|\langle X, {\bf{x}} \rangle \|_{\psi_2} 
\end{align}
\label{def:subgvec}
\end{defn}

\begin{lemm}
 Consider a sub-Gaussian vector $X \in \R^n$ with $\vertiii{X}_{\Psi_2} < K$, then for any vector $u$, $\langle X, u \rangle$ is a sub-Gaussian variable with $\vertiii{\langle X, u \rangle} < K \|u\|_2$.
 \label{lem:innprodsubgvec}
\end{lemm}
\proof
The argument is based on Definition \ref{def:subgvec} as follows,

\begin{align}
 \vertiii{\langle X, u \rangle}_{\Psi_2} =  \|u\|_2 \vertiii{\left \langle X, \frac{u}{\|u\|_2} \right \rangle}_{\Psi_2}  \leq \|u\|_2 \sup_{x\in S^{n-1}} \langle X, x \rangle = \|u\|_2 \vertiii{X}_{\Psi_2} \leq K \|u\|_2.
\end{align}
\qed

\begin{lemm}
Let $X_1$ and $X_2$ be centered sub-Gaussian random variables with $\vertiii{X_1}_{\Psi_2} = b_1$ and $\vertiii{X_2}_{\Psi_1} = b_2$. Then $X_1+X_2$ is centered sub-Gaussian with $\vertiii{X_1+X_2}_{\Psi_2} = b_1+b_2$.
\label{lem:sumsubgvar}
\end{lemm}
\proof
The argument is based on the definition of moment generating function of sub-Gaussian random variable:

Using Holder inequality for any $p,q>0$ where $\frac{1}{p}+\frac{1}{q}=1$, we have
\begin{align}
E[\exp\{t(X_1+X_2) \}] & \leq (E[\exp\{tX_1\}^p])^{1/p} (E[\exp\{tX_1\}^q])^{1/q} \nonumber \\
&\leq \exp\{Ct^2 (p b_1^2+q b_2^2)\} = \exp\{Ct^2 (p b_1^2+\frac{p}{1-p} b_2^2)\}.
\label{eq:holderMGFSumsubg}
\end{align}
The minimum of \eqref{eq:holderMGFSumsubg} occurs with $p=\frac{b_2}{b_1}$.
As a result, we have
\begin{align}
E[\exp\{t(X_1+X_2) \}] \leq \exp\{Ct^2 (b_1+ b_2)^2\}.
\end{align}
The proof is complete.
\qed

\subsection{Generic Chaining}


\begin{defn}[Majorizing measure \cite{tala96}]
Given $\alpha > 0$, and a metric space $(T,d)$ (that need not be finite), we define
\begin{align}
 \gamma_\alpha(T,d) = \inf \sup_t \sum_{n\geq 0} 2^{n/\alpha} \Delta(A_n(T)).
\end{align}
where the infimum is taken over all admissible sequences and $\Delta(A_n(T))$ is the diameter of $A_n(t)$.
\end{defn}
Note that $\gamma_2(T, \|. \|_2)$ coincides with the Gaussian width of $T$: $w(T)$.

\begin{lemm}
Given a metric space $(T,d)$, we have
\begin{align}
\gamma_1(T, \|.\|_\infty) \leq \gamma_2^2(T, \|.\|_2). 
\end{align}
\label{lem:bndgamma1}
 \end{lemm}
 
 \proof
Define $d_2(s,t) = \|s-t\|_2$ and $d_1(s,t) = \|s-t\|_\infty$.
We use the traditional definition of majorizing measure $\gamma_{\alpha,1}(T,d)$ from \cite{tala96} 
\begin{align}
\gamma_{\alpha,1}(T,d) = \inf \sup_t \left( \int_0^\infty \left( \log \frac{1}{\mu(B_d(t,\varepsilon))} \right)^{1/\alpha} d \varepsilon \right).
\end{align}
where $B_d(t,\varepsilon)$ is the closed ball of center $t$ and radius $\varepsilon$ based on the distance $d$ and the infimum is taken over all the probability measure $\mu$ on $T$.
 
Note that $\gamma_{\alpha,1}(T,d)$ coincides with the functional $\gamma_{\alpha}(T,d)$ \cite{tala01} as
\begin{align}
 K(\alpha)^{-1} \gamma_\alpha(T,d) \leq \gamma_{\alpha,1}(T,d) \leq K(\alpha) \gamma_\alpha(T,d),
\end{align}
where $K(\alpha)$ is a constant depending on $\alpha$ only.

As a result, it is enough to show that $\gamma_{1,1}(T, d_1) \leq \gamma_{2,1}^2(T, d_2)$.

Note that since for any vector $x$, we have $\| x \|_\infty \leq \|x \|_2$, therefore, for any probability measure $\mu$ and $t$, we have $\mu(B_{d_1}(t,\varepsilon)) \geq \mu(B_{d_2}(t,\varepsilon)) $.
As a result,
\begin{align}
\left( \int_0^\infty \left( \log \frac{1}{\mu(B_{d_1}(t,\varepsilon))} \right) d \varepsilon \right)   \leq \left( \int_0^\infty \left( \log \frac{1}{\mu(B_{d_2}(t,\varepsilon))} \right) d \varepsilon \right)
 \leq \left( \int_0^\infty \left( \log \frac{1}{\mu(B_{d_2}(t,\varepsilon))} \right)^{1/2} d \varepsilon \right)^2.
\label{eq:gammabndfixmu}
\end{align}

Since \eqref{eq:gammabndfixmu} holds for any $\mu$ and $t$, we have
\begin{align}
\gamma_{1,1}(T,d_1) = \inf \sup_t \left( \int_0^\infty \left( \log \frac{1}{\mu(B_{d_1}(t,\varepsilon))} \right) d \varepsilon \right) & \leq \inf \sup_t \left( \int_0^\infty \left( \log \frac{1}{\mu(B_{d_2}(t,\varepsilon))} \right)^{1/2} d \varepsilon \right)^2 \nonumber \\ 
&= \gamma^2_{2,1}(T,d_2).
\end{align}

This completes the proof.
\qed

\begin{theo} \textbf{[Theorem 1.2.7] in \cite{tala05}}
Consider a set $T$ provided with two distances $d_1$ and $d_2$. Consider a process $(X_t)_{t\in T}$ that satisfies $E[X_t] = 0$ and
\begin{align}
 P \left( |X_s - X_t| \geq u \right) \leq 2 \exp \left(- \min \left( \frac{u^2}{d_2(s,t)^2}, \frac{u}{d_1(s,t)} \right) \right).
\end{align}
Then
\begin{align}
 E [ \sup_{t,s \in T} |X_s - X_t| \leq L (\gamma_1(T,d_1) + \gamma_2(T,d_2)),
\end{align}
where 
$L$ is a constant.
\label{theo:talaGenExpbnd}
\end{theo}

\begin{theo} \textbf{[Theorem 1.2.9] in \cite{tala05}}
Under the conditions of Theorem \ref{theo:talaGenExpbnd}, for all values of $u_1, u_2 > 0$ we have
\begin{align}
 & P \left(|X_s - X_{t0}| \geq L (\gamma_1(T,d_1) + \gamma_2(T,d_2)) + u_1 D_1 + u_2 D_2 \right) \nonumber\\
 & \qquad \leq L \exp (- \min (u_2^2, u_1)),
\end{align}
where $D_j = 2 \sum_{n \geq 0} e_n(T, d_j)$. Note that $D_j \leq L \gamma_j(T, d_j)$.
\label{theo:talaConcSup}
\end{theo}

\begin{theo} \textbf{[Theorem  8.2 (Fernique-Talagrand's  comparison  theorem)] in \cite{vers14}}
Let $T$ be  an arbitrary  set. Consider a Gaussian  random  process $(G(t))_{t\in T}$ and a sub-Gaussian random process
$(H(t))_{t\in T}$. Assume that $E[G(t)] = E[H(t)] = 0$ for all $t\in T$. Assume also that for some $M > 0$, the following  increment comparison holds:
\begin{align}
 \vertiii{H(s)-H(t)}_{\psi_2} \leq M(E[ \| G(s) - G(t) \|_2^2)^{1/2} \qquad \forall s,t\in T.
\end{align}
Then
\begin{align}
 E[sup_{t\in T} H(t)] \leq C M E[ sup_{t\in T} G(t)].
\end{align}
\label{theo:FertalsubGau}
\end{theo}

\begin{theo}[Mendelson, Pajor, Tomczak-Jaegermann~\cite{mept07}]
There exist absolute constants $c_1$, $c_2$, $c_3$ for which the following
holds. Let $(\Omega,\mu)$ be a probability space, set $F$ be a subset of the unit
sphere of $L_2(\mu)$, i.e., $F \subseteq S_{L_2} = \{ f : \vertiii{f}_{L_2} = 1\}$, and assume that $ \sup_{f \in F}~\vertiii{f}_{\psi_2} \leq \kappa$. Then, for any $\theta > 0$ and $n \geq 1$ satisfying
\beq
c_1 \kappa \gamma_2(F, \vertiii{\cdot}_{\psi_2}) \leq \theta \sqrt{n}~,
\eeq
with probability at least $1- \exp(-c_2 \theta^2 n/\kappa^4)$,
\beq
\sup_{f \in F}~\left| \frac{1}{n} \sum_{i=1}^n f^2(X_i) - E\left[f^2\right] \right| \leq \theta~.
\eeq
Further, if $F$ is symmetric, then
\beq
E\left[ \sup_{f \in F} ~\left| \frac{1}{n} \sum_{i=1}^n f^2(X_i) - E\left[f^2\right] \right| \right] \leq c_3 \max \left\{ 2\kappa \frac{\gamma_2(F, \vertiii{\cdot}_{\psi_2})}{\sqrt{n}}, \frac{\gamma_2^2(F, \vertiii{\cdot}_{\psi_2})}{n} \right\}
\eeq
\label{thm:subGREmain}
\end{theo}

\section{Regularization Parameter}
  \begin{lemm}
Consider two Ising Model with true parameters $\theta_1^*$ and $\theta_2^*$. Let $d_1, d_2 \gg s$ where $\|\theta_1^*\|_0=d_1$, $\|\theta_2^*\|_0=d_2$, and $\|\delta\theta^*\|_0 = s$.
Assume
 \begin{align}
 \min_{i,j=1\cdots n_1}( | \theta_1^*(i,j) |) & \geq \frac{1}{d_1-1} - \frac{c_1}{(d_1-1)s} 
 \\ 
 \min_{i,j=1\cdots n_2}( | \theta_2^*(i,j) |)  & \geq \frac{1}{d_2-1} -\frac{c_2}{(d_2-1)s},
 \end{align}
 where $c_1$ and $c_2$ are positive constants.
 Then the density ratio $r(X = {\bf{x}}|\delta \theta^*)$ is bounded.
\end{lemm}

\proof
Let $\alpha_1 \leq | \theta_1^* | \leq \beta_1$ and $\alpha_2 \leq | \theta_2^* | \leq \beta_2$. 
Without loss of generality, assume that $\|\theta_1^*\|_2 = 1$ and $\|\theta_2^*\|_2 = 1$.

So,
\begin{align}
\beta_1 & \leq 1 - (d_1-1) \alpha_1 \\
\beta_2 &\leq 1 - (d_2-1) \alpha_2.
\end{align}

Based on triangle inequality of norms, we have
\begin{align}
\| \delta \theta^* \|_{\infty} = \| \theta_1^* - \theta_2^* \|_{\infty} \leq \| \theta_1^* \|_{\infty} + \|\theta_2^* \|_{\infty} \leq \beta_1 + \beta_2 \leq 2 - (d_1-1) \alpha_1 - (d_2-1) \alpha_2.
\end{align}
Let ${\bf{z}} = T({\bf{x}})$, then,
\begin{align}
| \langle {\bf{z}}, \delta \theta^* \rangle | & \leq \| {\bf{z}} \|_\infty \|\delta \theta^*\|_1 \\
& \leq s \|\delta \theta^*\|_\infty \\
& \leq 2s - [(d_1-1) \alpha_1 - (d_2-1) \alpha_2]s 
\end{align}
where the second inequality is the result of  $\| {\bf{z}} \|_\infty \leq 1$ since ${\bf{z}}$ comes from an Ising model.

Note that if 
\begin{align}
 \alpha_1 \geq \frac{s-c_1}{(d_1-1)s} = \frac{1}{d_1-1} - \frac{c_1}{(d_1-1)s},
\end{align}
then
\begin{align}
 s - (d_1-1) \alpha_1s \leq c_1.
\end{align}
Similarly, if
\begin{align}
 \alpha_2 \geq \frac{s-c_2}{(d_2-1)s} = \frac{1}{d_2-1} - \frac{c_2}{(d_2-1)s},
\end{align}
then
\begin{align}
 s - (d_2-1) \alpha_2s \leq c_2.
\end{align}
As a result, we have
\begin{align}
& | \langle {\bf{z}}, \delta \theta^* \rangle |  \leq c_1 + c_2. \\
\Rightarrow & \exp\{\langle {\bf{z}}, \delta \theta^* \rangle \}  \leq \exp{c_1 + c_2} \leq c_0.
\end{align}
For example, if $c_1=c_2=1$, then $c_0 = 10$.

Therefore, 
\begin{align}
r(X = {\bf{x}}|\delta \theta) =\frac{\exp\{ \langle {\bf{z}}, \delta \theta \rangle) \}}{Z(\delta \theta^*)} \leq \frac{c_0}{Z(\delta \theta^*)}.
\end{align}
This completes the proof. \qed
 
{ \bf{Assumption 1(Smooth Density Ratio Model Assumption)}} For any vector ${\bf{u}}$ such that $\|{\bf{u}}\|_2 \leq \| \delta \theta^*\|_2$ and every $t \in R$, the following inequality holds:
\begin{align}
E_{X \sim p(X|\theta_2)} [\exp \{ t r(X | \delta \theta^*+ {\bf{u}}) - 1 \}] \leq \exp \{ t^2 \}.
\end{align}

\begin{lemm}
 For any constant $\tau \leq 1$, define random event $M_\tau$ as follows,
 \begin{align}
 M_\tau = \{ \Psi(\delta \theta^*+ {\bf{u}}) - \Psi(\delta \theta^*) - \left[\hat{\Psi}(\delta \theta^*+{\bf{u}}) - \hat{\Psi}(\delta \theta^*)\right] \leq \tau \}.
 \end{align}
Then, for any vector ${\bf{u}}$ such that $\|{\bf{u}}\|_2 \leq \| \delta \theta^*\|_2$, under Assumption 1, we have
\begin{align}
P(M_\tau^c) = p \left( \Psi(\delta \theta^*+ {\bf{u}}) - \Psi(\delta \theta^*) -  \left[\hat{\Psi}(\delta \theta^*+{\bf{u}}) - \hat{\Psi}(\delta \theta^*)\right] > \tau \right) \leq 4 e^{- \frac{ n_2}{5} \tau^2}.
\end{align}
\label{lem:eventM}
\end{lemm}
 \proof
Recall that 
\begin{align}
& r(X={\bf{x}} | \delta \theta^*) = \frac{\exp\{\langle T({\bf{x}}), \delta \theta^* \rangle \}}{Z(\delta \theta^*)} \nonumber \\
\Rightarrow & \qquad \hat{Z}(\delta \theta^*) = \frac{1}{n_2} \sum_{i=1}^{n_2} \exp\{ \langle T({\bf{x}}_i^2), \delta \theta^* \rangle \} = \frac{1}{n_2} \sum_{i=1}^{n_2}  r(X={\bf{x}}_i^2 | \delta \theta) Z(\delta \theta^*) \nonumber \\
\Rightarrow & \qquad  \frac{\hat{Z}(\delta \theta^*)}{Z(\delta \theta^*)} = \frac{1}{n_2} \sum_{i=1}^{n_2}  r(X={\bf{x}}_i^2 | \delta \theta^*)
\end{align}
Note that $Z(\delta \theta) = E_{X \sim p(X|\theta_2)} [\exp\{\langle T({\bf{x}}), \delta \theta^* \rangle \}]$, therefore,
\begin{align}
E_{X \sim p(X|\theta_2)} [r(X | \delta \theta^*)] = 1.
\end{align}

Under the Assumption 1, we have
\begin{align}
p \left(\left|r(X={\bf{x}}_i^2 | \delta \theta^*) - 1\right| > \epsilon \right) \leq c_1 e^{- \epsilon^2}.
\end{align}
Applying Hoeffding inequality, we have
\begin{align}
& p \left( \left|\frac{1}{n_2} \sum_{i=1}^{n_2}  r(X={\bf{x}}_i^2 | \delta \theta^*) - 1 \right| \geq \epsilon \right) \leq 2 e^{-\epsilon^2} \\
\Rightarrow & \qquad p \left( \left |\frac{\hat{Z}(\delta \theta^*)}{Z(\delta \theta^*)} - 1 \right| \geq \epsilon \right) \leq 2 e^{-n_2 \epsilon^2}.
\end{align}

Taking the logarithm from both side, and using one side bound, we have
\begin{align}
& p \left( \log \frac{\hat{Z}(\delta \theta^*)}{Z(\delta \theta^*)} \geq \log(\epsilon + 1) \right) \leq e^{-n_2 \epsilon^2} \\
\Rightarrow & p \left( \hat{\Psi}(\delta \theta^*) - \Psi(\delta \theta^*) \geq \log(\epsilon + 1) \right) \leq e^{-n_2 \epsilon^2}.
\end{align}

Similarly, we have
\begin{align}
p \left( \Psi(\delta \theta^* + {\bf{u}}) - \hat{\Psi}(\delta \theta^* + {\bf{u}}) \geq - \log(1-\epsilon) \right) \leq e^{-n_2 \epsilon^2}.
\end{align}
 Applying the union bound, we have
 \begin{align}
p \left( \Psi(\delta \theta^* + {\bf{u}}) - \Psi(\delta \theta^*) - \left[\hat{\Psi}(\delta \theta^* + {\bf{u}}) - \hat{\Psi}(\delta \theta^*) \right]   \geq \log \frac{1+\epsilon}{1-\epsilon} \right) \leq 4 e^{-n_2 \epsilon^2}.
 \end{align}

Setting $\tau = \log \frac{1+\epsilon}{1-\epsilon} $, we have
 \begin{align}
p \left( \Psi(\delta \theta^* + {\bf{u}}) - \Psi(\delta \theta^*) - \left[\hat{\Psi}(\delta \theta^* + {\bf{u}}) - \hat{\Psi}(\delta \theta^*) \right]   \geq \tau \right) \leq 4 e^{- n_2 \left(\frac{e^\tau+1}{e^\tau-1}\right)^2} \leq 4 e^{- \frac{ n_2}{5} \tau^2},
 \end{align}
where the last inequality is obtained by using the fact that for any $\tau \leq 1$
\begin{align}
\left(\frac{e^\tau+1}{e^\tau-1}\right)^2 > \frac{\tau^2}{5}.
\end{align}

This completes the proof. \qed

 \begin{lemm}
Define random event $M_{\tilde{t}}$ as follows,
 \begin{align}
  M_{\tilde{t}} = \{ \Psi(\delta \theta^*+ t {\bf{u}}) - \Psi(\delta \theta^*) - \left[\hat{\Psi}(\delta \theta^*+ t {\bf{u}}) - \hat{\Psi}(\delta \theta^*) \right] \leq \tilde{t} \}, 
 \label{eq:eventMt}
 \end{align}
 where $\tilde{t} = \sqrt{5\eta_1}t+\frac{\sqrt{5}}{2}$. Let $Z  = T(X^1)$ and ${\bf{z}}  = T({\bf{x}}^1)$. Then, 
 \begin{align}
 P \left( \left| \left \langle {\bf{z}} - \nabla \hat{\Psi}(\theta^*), {\bf{u}} \right \rangle \right| \geq \epsilon \big|  M_{\tilde{t}} \right) \leq c\exp\{\frac{-\epsilon^2}{4\eta_0 \|{\bf{u}}\|_2^2}\},
 \end{align}
  where $\frac{1}{2} \lambda_{\max} \left(\nabla^2 \hat{\Psi}(\theta^*+\tilde{{\bf{u}}})\right)  \leq \eta_0$ and $c$ is a positive constant.
  \label{lem:momgenBnd}
 \end{lemm}
 \proof
 First, note that $p(X = {\bf{x}} | \theta_1^*) = p(X = {\bf{x}}|\theta_2^*) r(X = {\bf{x}}|\delta \theta^*)$. 
 Therefore,
 \begin{align}
 E_{X \sim p(X|\theta_1^*)} \left[e^{ \left \langle Z , t {\bf{u}} \right \rangle } \right] 
   & = \sum_{{\bf{x}} \in \mathcal{X}} e^{ \left \langle {\bf{z}} , t {\bf{u}} \right \rangle } p({\bf{x}} |\theta_1^*) \nonumber \\
   &=  \sum_{{\bf{x}} \in \mathcal{X}} e^{ \left \langle {\bf{z}} , t{\bf{u}} \right \rangle } p({\bf{x}}|\theta_2^*) r({\bf{x}}|\delta \theta^*) \nonumber \\
   &=  e^{-\Psi(\delta \theta^*)} \sum_{{\bf{x}} \in \mathcal{X}} e^{ \left \langle {\bf{z}} , t{\bf{u}}+\delta \theta^* \right \rangle } p({\bf{x}} |\theta_2^*) \nonumber \\
  &= e^{\Psi(\delta \theta^*+t{\bf{u}})-\Psi(\delta \theta^*)},
  \label{eq:expectPZt}
  \end{align}
  since $r({\bf{x}}|\delta \theta^*) = \exp\{ {\bf{x}} , \delta \theta^* \rangle - \Psi(\delta \theta^*)\}$, and $\Psi(\delta \theta^*) =\log  \sum_{{\bf{x}} \in \mathcal{X}} e^{ \left \langle {\bf{x}} , \delta \theta^* \right \rangle } p({\bf{x}}|\theta_2^*)$.
 
 Also, using the Taylor expansion, we have
  \begin{align}
    \hat{\Psi}(\delta \theta^*+t {\bf{u}}) - \hat{\Psi}(\delta \theta^*) - \left\langle \nabla \hat{\Psi}(\delta \theta^*), t{\bf{u}} \right\rangle
    &= \frac{1}{2} t{\bf{u}}^T \nabla^2 \hat{\Psi}(\delta \theta^*+\tilde{{\bf{u}}}) t{\bf{u}} \nonumber \\
    & \leq \frac{1}{2} t^2 \|{\bf{u}}\|_2^2 \lambda_{\max}\left(\nabla^2 \hat{\Psi}(\delta \theta^*+\tilde{{\bf{u}}})\right) \leq t^2 \eta_0 \|{\bf{u}}\|_2^2 = t^2 \eta_1,
    \label{eq:upperbndBregmanPsiHat}
  \end{align}
  where $\frac{1}{2}  \lambda_{\max} \left(\nabla^2 \hat{\Psi}(\delta \theta^*+\tilde{{\bf{u}}} \right) \leq \eta_0$ and $\eta_1=\eta_0 \|{\bf{u}}\|_2^2$.
 
 Then, given the event $ M_{\tilde{t}}$, the moment generating function of $\left \langle Z - \nabla \hat{\Psi}(\delta \theta^*), {\bf{u}} \right \rangle$ can be bounded as,
 \begin{align}
  E_{X \sim p(X|\theta_2^*)}\left[ E_{X \sim p(X|\theta_1^*)} \left[e^{ \left \langle Z - \nabla \hat{\Psi}(\delta \theta^*), t{\bf{u}} \right \rangle } \right] \big|  M_{\tilde{t}} \right]
  & = E_{X \sim p(X|\theta_2^*)}\left[ e^{\left\langle -\nabla \hat{\Psi}(\delta \theta^*), t{\bf{u}} \right\rangle} ~~ E_{X \sim p(X|\theta_1^*)} \left[e^{ \left \langle Z , t{\bf{u}} \right \rangle } \right] \big|  M_{\tilde{t}} \right] \nonumber \\
   & \stackrel{\eqref{eq:expectPZt}}{=} E_{X \sim p(X|\theta_2^*)}\left[ e^{\Psi(\delta \theta^*+t{\bf{u}})-\Psi(\delta \theta^*) - \left\langle \nabla \hat{\Psi}(\delta \theta^*), t{\bf{u}} \right\rangle} \big|  M_{\tilde{t}} \right] \nonumber \\
   & \stackrel{\eqref{eq:eventMt}}{\leq} E_{X \sim p(X|\theta_2^*)}\left[ e^{\hat{\Psi}(\delta \theta^*+t{\bf{u}}) - \hat{\Psi}(\delta \theta^*) - \left\langle \nabla \hat{\Psi}(\delta \theta^*), t{\bf{u}} \right\rangle + \sqrt{\eta_1} t} \big|  M_{\tilde{t}} \right] \nonumber \\
     & \stackrel{\eqref{eq:upperbndBregmanPsiHat}}{\leq} E_{X \sim p(X|\theta_2^*)}\left[ e^{t^2 \eta_1 + \sqrt{5 \eta_1} t + \frac{1}{2}} \big|  M_{\tilde{t}} \right] = e^{t^2 \eta_1 + \sqrt{5 \eta_1} t + \frac{1}{2}}.
  \label{eq:mgfbnd}
 \end{align}
 
 As a result, using the Chernoff bound, for any $t>0$, we have
 \begin{align}
 P \left( \left \langle Z - \nabla \hat{\Psi}(\delta \theta^*), {\bf{u}} \right \rangle \geq \epsilon \big|  M_{\tilde{t}} \right) 
  & \leq e^{-t\epsilon} E_{X \sim p(X|\theta_2^*)}\left[ E_{X \sim p(X|\theta_1^*)} \left[e^{ \left \langle Z - \nabla \hat{\Psi}(\delta \theta^*), t{\bf{u}} \right \rangle } \right] \big|  M_{\tilde{t}} \right] \nonumber \\
 & \stackrel{\eqref{eq:mgfbnd}}{\leq} \exp\{t^2 \eta_1 + \sqrt{5 \eta_1} t + \frac{1}{2} -t\epsilon \} \stackrel{(a)}{\leq} \exp\{-\frac{(\epsilon-\sqrt{5 \eta_1})^2}{4\eta_1}+\frac{1}{2}\} \nonumber \\
 & \leq c \exp\{-\frac{\epsilon^2}{4\eta_0 \|{\bf{u}}\|_2^2}\},
 \label{eq:cndprobnd}
 \end{align}
 where the inequality $(a)$ is obtained by setting $t = \frac{\epsilon-\sqrt{\eta_1}}{2\eta_1}$ to minimize it with respect to $t$, and the last inequality obtained by setting $c = \exp\{\frac{5\sqrt{5}}{2 \sqrt{\eta_1}} - \frac{3}{4}\}$ and using the fact that
  \begin{align}
 \exp\{-\frac{(\epsilon-\sqrt{5\eta_1})^2}{4\eta_1}+\frac{1}{2}\} \leq c \exp\{-\frac{\epsilon^2}{4\eta_1}\}.
 \end{align}
 
 Similarly, we have
  \begin{align}
 P \left( \left \langle Z - \nabla \hat{\Psi}(\delta \theta^*), {\bf{u}} \right \rangle \leq -\epsilon \big|  M_{\tilde{t}} \right) 
  & \leq e^{-t\epsilon} E_{X \sim p(X|\theta_2^*)}\left[ E_{X \sim p(X|\theta_1^*)} \left[e^{ \left \langle Z - \nabla \hat{\Psi}(\delta \theta^*), -t{\bf{u}} \right \rangle } \right] \big|  M_{\tilde{t}} \right] \nonumber \\
 &\leq c \exp\{-\frac{\epsilon^2}{4\eta_0 \|{\bf{u}}\|_2^2}\}.
 \end{align}
 
 This completes the proof.\qed

 \begin{lemm}
Under the smooth density ratio assumption, we have
 \begin{align}
 P \left( \big| \left \langle \nabla \cL(\delta \theta^*; \mathfrak{X}_1^{n_1}, \mathfrak{X}_2^{n_2}) , {\bf{u}}  \right \rangle \big | \geq \epsilon  \right)
 \leq c_1 \exp\{-\frac{\min(n_1,n_2) \epsilon^2}{4\eta_0 \|{\bf{u}}\|_2^2}\},
 \end{align}
 where $c_1$ is a positive constant. 
 \label{lem:gradLBnd}
 \end{lemm}
 
\proof 
Let $\tilde{t} = \sqrt{5\eta_1}t+\frac{\sqrt{5}}{2}$ and $t = \frac{\epsilon-\sqrt{\eta_1}}{2\eta_1}$. Using the result of lemma \ref{lem:momgenBnd} we have
\begin{align}
  P \left( \left \langle T({\bf{x}}_i^1) - \nabla \hat{\Psi}(\delta \theta^*), {\bf{u}} \right \rangle \geq \epsilon \big|  M_{\tilde{t}} \right) \nonumber &= P \left( \left \langle T({\bf{x}}_i^1) - \nabla \hat{\Psi}(\delta \theta^*), {\bf{u}} \right \rangle \geq {\epsilon} \big|  M_{\tilde{t}} \right) \nonumber \\
  & \leq c \exp\{-\frac{\epsilon^2}{4\eta_0 \|{\bf{u}}\|_2^2}\}.
 \end{align}
Applying Hoeffding inequality, we have
 \begin{align}
  P \left( \big| \left \langle \nabla \cL(\delta \theta^*; \mathfrak{X}_1^{n_1}, \mathfrak{X}_2^{n_2}) , {\bf{u}}  \right \rangle \big | \geq \epsilon  \big|  M_{\tilde{t}} \right) &= P \left( \left \langle \frac{1}{n_1} \sum_{i=1}^{n_1} T({\bf{x}}_i^1) - \nabla \hat{\Psi}(\delta \theta^*), {\bf{u}} \right \rangle \geq \epsilon \big|  M_{\tilde{t}} \right) \nonumber \\
  & \leq c \exp\{-\frac{n_1 \epsilon^2}{4\eta_0 \|{\bf{u}}\|_2^2}\}.
 \end{align}
 Moreover, we can obtain,
 \begin{align}
 P \left(\left \langle \nabla \cL(\delta \theta^*; \mathfrak{X}_1^{n_1}, \mathfrak{X}_2^{n_2}) , {\bf{u}} \right \rangle \leq -\epsilon  \right)
 &= P \left( \left \langle \frac{1}{n_1} \sum_{i=1}^{n_1} T({\bf{x}}_i^1) - \nabla \hat{\Psi}(\delta \theta^*), {\bf{u}} \right \rangle \geq \epsilon \right) \nonumber \\
  & \leq P \left(\left \langle \frac{1}{n_1} \sum_{i=1}^{n_1} T({\bf{x}}_i^1) - \nabla \hat{\Psi}(\delta \theta^*), {\bf{u}} \right \rangle \geq \epsilon \big|  M_{\tilde{t}} \right) P( M_{\tilde{t}}) \nonumber \\
  & \quad + P \left(\left \langle \frac{1}{n_1} \sum_{i=1}^{n_1} T({\bf{x}}_i^1) - \nabla \hat{\Psi}(\delta \theta^*), {\bf{u}} \right \rangle \geq \epsilon \big|  M_{\tilde{t}} ^c\right) P(  M_{\tilde{t}}^c) \nonumber \\
  & \leq c \exp\{\frac{-n_1 \epsilon^2}{4\eta_0 \|{\bf{u}}\|_2^2}\} + 
  4 \exp \{ \frac{-n_2\epsilon^2}{4\eta_0 \|{\bf{u}}\|_2^2}\} \nonumber \\
  & \leq c_1 \exp\{-\frac{\min(n_1,n_2) \epsilon^2}{4\eta_0 \|{\bf{u}}\|_2^2}\},
 \end{align}
 where the last inequality is obtained by using Lemma \ref{lem:eventM} as follows
 \begin{align*}
 P(  M_{\tilde{t}}^c)  \leq 4 \exp \{ -\frac{n_2}{5} \tilde{t}^2\}  &= 4 \exp \{ -\frac{n_2}{5} (\sqrt{5 \eta_1} t+ \frac{\sqrt{5}}{2})^2\} \\
 &= 4 \exp \{ -\frac{n_2}{5} (\sqrt{5 \eta_1} \frac{\epsilon'-\sqrt{\eta_1}}{2 \eta_1}+ \frac{\sqrt{5}}{2})^2\} \\ 
  &= 4 \exp \{ -n_2 \frac{\epsilon^2}{4 \eta_1}\}, 
 \end{align*}
 where $\eta_1=\eta_0\|{\bf{u}}\|_2^2$ and setting $c_1 = \max(4, c)$.
Similarly,
 \begin{align}
 P \left(\left \langle \nabla \cL(\delta \theta^*; \mathfrak{X}_1^{n_1}, \mathfrak{X}_2^{n_2}) , {\bf{u}} \right \rangle \geq \epsilon  \right)
 &= P \left(\left \langle \frac{1}{n_p} \sum_{i=1}^{n_p} Z_i^p - \nabla \hat{\Psi}(\delta \theta^*), {\bf{u}} \right \rangle \leq -\epsilon \right) \nonumber \\
  & \leq P \left(\left \langle \frac{1}{n_p} \sum_{i=1}^{n_p} Z_i^p - \nabla \hat{\Psi}(\delta \theta^*), {\bf{u}} \right \rangle \leq -\epsilon \big|  M_{\tilde{t}} \right) P(M_{t^2}) \nonumber \\
  & \quad + P \left(\left \langle \frac{1}{n_p} \sum_{i=1}^{n_p} Z_i^p - \nabla \hat{\Psi}(\delta \theta^*), {\bf{u}} \right \rangle \leq -\epsilon \big|  M_{\tilde{t}}^c \right) P(  M_{\tilde{t}}^c) \nonumber \\
 & \leq c_1 \exp\{-\frac{\min(n_1,n_2) \epsilon^2}{4\eta_0 \|{\bf{u}}\|_2^2}\}
 \end{align}
This completes the proof. \qed

\textbf{Theorem \ref{theo:lambdabnd}} \textit{
Define $\Omega_R = \{u: R(u) \leq 1\}$. Let $\phi(R) = \sup_{{\bf{u}}}  \frac{\|{\bf{u}}\|_2}{R({\bf{u}})}$. Assume that for any ${\bf{u}}$ that $\|{\bf{u}}\| \leq \| \theta^*\|$
\begin{align}
\frac{1}{2}  \lambda_{\max} \left(\nabla^2 {\cL}(\delta \theta^*+{{\bf{u}}}) \right) \leq \eta_0,
\end{align}
where $\lambda_{\max}(.)$ is the maximum eigenvalue.
Then under the smooth density ratio assumption, we have
\begin{align}
E\left[R^*(\nabla \cL(\delta \theta^*; \mathfrak{X}_1^{n_1}, \mathfrak{X}_2^{n_2}))\right] \leq \frac{2 \sqrt{\eta_0}}{\sqrt{\min(n_1,n_2)}} (c_1 w\left(\Omega_R) + \phi(R) \right).
\end{align}
and with probability at least $1-c_2 e^{-\epsilon^2}$
\begin{align}
R^*\left(\nabla \cL(\delta \theta^*; \mathfrak{X}_1^{n_1}, \mathfrak{X}_2^{n_2}) \right) \leq  \frac{1}{\sqrt{\min(n_1, n_2)}} \left( c_2(1+\epsilon) w(\Omega_R) + \tau_1 \right).
\label{eq:conbnd}
\end{align}
where $c_1$ and $c_2$ are positive constants, $\tau_1 = 2 \sqrt{\eta_0} \phi(R)$, and $w(\Omega_R)$ is the Gaussian width of set $\Omega_R$.
}

\proof
Define $\boldsymbol{\mu} = E[\nabla \cL(\delta \theta^*; \mathfrak{X}_1^{n_1}, \mathfrak{X}_2^{n_2})]$.
Using the triangle inequality, we have
\begin{align}
R^*\left(\nabla \cL(\delta \theta^*; \mathfrak{X}_1^{n_1}, \mathfrak{X}_2^{n_2}) \right) \leq
R^*\left(\nabla \cL(\delta \theta^*; \mathfrak{X}_1^{n_1}, \mathfrak{X}_2^{n_2}) - \boldsymbol{\mu} \right) + R^*\left(\boldsymbol{\mu} \right).
\end{align}
We upper bound two terms as follows. First, consider the first term.

Using the definition of dual norm, we have
\begin{align}
R^*\left(\nabla \cL(\delta \theta^*; \mathfrak{X}_1^{n_1}, \mathfrak{X}_2^{n_2}) -\boldsymbol{\mu} \right) = \sup_{R({\bf{u}}) \leq 1} \left \langle \nabla \cL(\delta \theta^*; \mathfrak{X}_1^{n_1}, \mathfrak{X}_2^{n_2}) - \boldsymbol{\mu}, {\bf{u}} \right \rangle.
\end{align}
Define stochastic process $H({\bf{s}}) = \left \langle \nabla \cL(\delta \theta^*; \mathfrak{X}_1^{n_1}, \mathfrak{X}_2^{n_2})  - \boldsymbol{\mu}, {\bf{s}} \right \rangle$ where $E[H({\bf{s}})] = 0$. Then, from Lemma \ref{lem:gradLBnd}, we have
\begin{align}
P \left( H({\bf{s}}) - H({\bf{t}}) \geq \epsilon  \right) &= 
P \left(\left \langle \nabla \cL(\delta \theta^*; \mathfrak{X}_1^{n_1}, \mathfrak{X}_2^{n_2})- \boldsymbol{\mu} , {\bf{s-t}} \right \rangle \geq \epsilon  \right) \\
& \leq  c_1 \exp\{-\frac{\min(n_1,n_2) \epsilon^2}{4\eta_0 \|{\bf{s-t}}\|_2^2} \}.
\end{align}

Consider the Gaussian process $G({\bf{u}}) = \langle {\bf{u}}, g \rangle$, indexed by the same set, i.e., ${\bf{u}} \in \Omega_R$, where $g \sim N (0, \I_{{d} \times {d}})$ is standard Gaussian vector. 
Now from definition sub-Gaussian random variables, we have
\begin{align}
\vertiii{H({\bf{s}}) - H({\bf{t}})}_{\psi_2} \leq \frac{2 \sqrt{\eta_0} \|{\bf{s-t}}\|_2}{\sqrt{\min(n_1,n_2)}} = K E_g[ \| G({\bf{s}}) - G({\bf{t}}) \|_2^2]^{1/2},
\end{align}
where $E_g[ \| G({\bf{s}}) - G({\bf{t}}) \|_2^2]^{1/2} = E_g[\| \langle {\bf{s}} - {\bf{t}} , g \rangle\|_2^2]^{1/2} = \| {\bf{s}} - {\bf{t}} \|_2$, and $K = \frac{2 \sqrt{\eta_0}}{\sqrt{\min(n_1,n_2)}}$.

Next, by applying the Fernique-Talagrand's comparison theorem \ref{theo:FertalsubGau}, we have
\begin{align}
E[\sup_{{\bf{u}} \in \Omega_R} H({\bf{u}})] &= 
E\left[ \sup_{{\bf{u}} \in \Omega_R} 
\left \langle \nabla \cL(\delta \theta^*; \mathfrak{X}_1^{n_1}, \mathfrak{X}_2^{n_2}) , {\bf{u}} \right \rangle \right]\nonumber \\
& \leq c_1 K E[ \sup_{u\in \Omega_R} G(u)] =  2 c_1 \sqrt{\eta_0} \frac{w(\Omega_R)}{\sqrt{\min(n_1,n_2)}},
\end{align}
where 
$c_1$ is a constant.
Thus,
\begin{align}
E\left[R^*\left(\nabla \cL(\delta \theta^*; \mathfrak{X}_1^{n_1}, \mathfrak{X}_2^{n_2}) - \boldsymbol{\mu} \right)\right] \leq c_1 \frac{w(\Omega_R)}{\sqrt{\min(n_1,n_2)}}.
\end{align}

To get the concentration bound, we use the direct application of Theorem 2.2.27 in \cite{tala14} and we have
\begin{align}
& P \left( \sup_{{\bf{s}}, {\bf{t}} \in \Omega_R} | H({\bf{s}}) - H({\bf{t}})| \leq c_2 (1+\epsilon) \frac{w(\Omega_R)}{\sqrt{\min(n_1,n_2)}}\right)  \geq 1- c_2 \exp \left(- \epsilon^2 \right).
\end{align}
Thus, with probability at least $1- c_2 \exp \left(- \epsilon^2 \right)$,
 \begin{align}
R^*\left(\nabla \cL(\delta \theta^*; \mathfrak{X}_1^{n_1}, \mathfrak{X}_2^{n_2})  - \boldsymbol{\mu} \right) \leq c_2(1+\epsilon) \frac{w(\Omega_R)}{\sqrt{\min(n_1,n_2)}}.
\end{align}

Next, we consider the second term. 
First note that $\vertiii{\nabla \cL(\delta \theta^*; \mathfrak{X}_1^{n_1}, \mathfrak{X}_2^{n_2})}_{\Psi_2} \leq \frac{2 \sqrt{\eta_0} \|{\bf{u}}\|_2}{\sqrt{\min(n_1,n_2)}}$. Using sub-Gaussian variables property, we have
\begin{align}
E[\cL(\delta \theta^*; \mathfrak{X}_1^{n_1}, \mathfrak{X}_2^{n_2})] \leq \vertiii{\nabla \cL(\delta \theta^*; \mathfrak{X}_1^{n_1}, \mathfrak{X}_2^{n_2})}_{\Psi_2} \leq \frac{2 \sqrt{\eta_0} \|{\bf{u}}\|_2}{\sqrt{\min(n_1,n_2)}}
\end{align}

Using the definition of the dual norm, we have
\begin{align}
R^*(\boldsymbol{\mu}) = R^*\left( E\left[\nabla \cL(\delta \theta^*; \mathfrak{X}_1^{n_1}, \mathfrak{X}_2^{n_2}) \right] \right) & = \sup_{{\bf{u}} \in \Omega_R} E\left[\left \langle \nabla \cL(\delta \theta^*; \mathfrak{X}_1^{n_1}, \mathfrak{X}_2^{n_2}), {\bf{u}} \right \rangle \right] \\
& \leq \frac{2 \sqrt{\eta_0}}{\sqrt{\min(n_1,n_2)}} \sup_{{\bf{u}}}  \frac{\|{\bf{u}}\|_2}{R({\bf{u}})} = \frac{2 \sqrt{\eta_0}}{\sqrt{\min(n_1,n_2)}} \Phi(R),
\end{align}
where $\Phi(R) = \sup_{{\bf{u}}}  \frac{\|{\bf{u}}\|_2}{R({\bf{u}})}$.

Also, we have
\begin{align}
E\left[R^*(\boldsymbol{\mu})\right] = \leq \frac{2 \sqrt{\eta_0}}{\sqrt{\min(n_1,n_2)}} \Phi(R),
\end{align}
where $\Phi(R) = \sup_{{\bf{u}}}  \frac{\|{\bf{u}}\|_2}{R({\bf{u}})}$.

This completes the proof. \qed

\section{RSC condition}
\label{sec:supRSC}
Let $r_i = r(X={\bf{x}_i^2} | \delta \theta^*)$ and $\bar{\varepsilon}$ denote the probability that $r_i$ exceeds some constant $\eta_0$: $\bar{\varepsilon} = p(r_i > \eta_0) \geq 1-  e^{-\frac{\eta_0^2}{2}}$.

\textbf{Theorem \ref{theo:rsc}} \textit{
Let $X \in \R^{n \times p}$ be a design matrix with independent isotropic sub-Gaussian rows with $\vertiii{X_i}_{\Psi_2} \leq \kappa$.
Then, for any set $A \subseteq S^{p-1}$, for suitable constants $\eta$, $c_1$, $c_2 >0$ with probability at least $
1-\exp(-\eta w^2(A))$, we have
\begin{align}
\inf_{u \in A} \partial \cL(\theta^*;u,X) \geq c_1 \underline{\rho}^2 \left ( 1 - c_2\kappa_1^2 \frac{w(A)}{\sqrt{n_2}} \right) - \tau
\end{align}
where $\kappa_1=\frac{\kappa}{\bar{\varepsilon}}$, $\underline{\rho}^2 = \inf_{{\bf{u}} \in A} \rho_{\bf{u}}^2$ with
$\rho_{\bf{u}}^2 = E\left[\left\langle {\bf{u}}, T(X_i^2) \right\rangle^2 \mathbb{I} (r_i > \eta_0)\right]$, and $\tau$ is smaller than the first term in right hand side. Thus, for $n_2 \geq c_2 w^2(A)$, with probability at least $1- \exp(-\eta w^2(A))$, we have
$\inf_{u \in A} \partial \cL(\theta^*;u,X) > 0$.
}

\proof
Define $Z=T(X)$ and ${\bf{z}}_i = T({\bf{x}}_i^2)$. Then,
\begin{align}
\cL(\delta \theta; \mathfrak{X}_1^{n_1}, \mathfrak{X}_2^{n_2}) &= \frac{-1}{n_1} \sum_{i=1}^{n_1} \langle T({\bf{x}}_i^1), \delta \theta \rangle + \log \frac{1}{n_2} \sum_{i=1}^{n_2} \exp\{ \langle T({\bf{x}}_i^2), \delta \theta \rangle \} \\
& = \frac{-1}{n_1} \sum_{i=1}^{n_1} \langle {\bf{z}_i}, \delta \theta \rangle + \log \frac{1}{n_2} \sum_{i=1}^{n_2} \exp\{ \langle {\bf{z}}_i, \delta \theta \rangle \}.
\end{align}
Through the analysis, we consider that $Z$ is centered random variable without loss of generality, since if it is not, the $E[Z]$ will show up as a constant.

Recall, RSC condition definition as
\begin{align}
\delta \cL(\delta \theta^{*},{\bf{ u}}) :=  \cL(\delta \theta^{*} + {\bf{u}}; \mathfrak{X}_1^{n_1}, \mathfrak{X}_2^{n_2}) - \cL(\delta \theta^{*}; \mathfrak{X}_1^{n_1}, \mathfrak{X}_2^{n_2}) - \langle \nabla \cL(\delta \theta^{*}; \mathfrak{X}_1^{n_1}, \mathfrak{X}_2^{n_2}), {\bf{u}} \rangle  \geq \kappa \|{\bf{u}}\|_{2}^{2}
\end{align}

Simplifying the expression and applying mean value theorem twice on the left side of RSC condition \eqref{eq:RSC}, for $\forall \gamma_i \in [0, 1]$, we have
\begin{align}
\delta \cL(\delta \theta^{*}, {\bf{u}})  &:= \cL(\delta \theta^{*} + {\bf{u}}; \mathfrak{X}_1^{n_1}, \mathfrak{X}_2^{n_2}) - \cL(\delta \theta^{*}; \mathfrak{X}_1^{n_1}, \mathfrak{X}_2^{n_2}) - \langle \nabla \cL(\delta \theta^{*}; \mathfrak{X}_1^{n_1}, \mathfrak{X}_2^{n_2}), {\bf{u}} \rangle \nonumber \\
& \geq {\bf{u}}^T \nabla^2 \cL(\delta \tilde{\theta}; \mathfrak{X}_1^{n_1}, \mathfrak{X}_2^{n_2}) {\bf{u}},
\label{eq:RSCextend}
\end{align}
where $\delta \tilde{\theta} = \delta \theta^*+\gamma_i {\bf{u}}$.
As a result, to show when the RSC condition is satisfied, it is enough to find a lower bound for the right side of the above equation.

Note that 
\begin{align}
\nabla^2 \cL(\delta \tilde{\theta}; \mathfrak{X}_1^{n_1}, \mathfrak{X}_2^{n_2}) = \nabla^2 \hat{\Psi}(\delta \tilde{\theta}),
\end{align}
where
\begin{align}
\nabla^2 \hat{\Psi}(\delta \tilde{\theta}) = \sum_{i=1}^{n_2} \sigma_i {\bf{z}}_i^T {\bf{z}}_i - \left(\sum_{j=1}^{n_2} \sigma_j {\bf{z}}_j \right)^T \left(\sum_{j=1}^{n_2} \sigma_j {\bf{z}}_j \right), \label{eq:secDerLogPart}
\end{align}
and
\beq
\sigma_i = \exp \{\langle {\bf{z}}_i, \delta \tilde{\theta} \rangle - \hat{\Psi}(\delta \tilde{\theta})\} = \frac{\exp{\langle {\bf{z}}_i, \delta \tilde{\theta} \rangle}}{\sum_{j=1}^{n_2} \exp{\langle {\bf{z}}_j, \delta \tilde{\theta} \rangle}}.
\eeq
Putting \eqref{eq:secDerLogPart} back in \eqref{eq:RSCextend}, we have
\begin{align}
\delta \cL(\delta \theta^{*}, {\bf{u}}) \geq \underbrace{\sum_{i=1}^{n_2} \sigma_i \langle {\bf{u}}, {\bf{z}}_i \rangle^2}_{A} - \underbrace{\left \langle {\bf{u}}, \sum_{j=1}^{n_2} \sigma_j {\bf{z}}_j \right \rangle^2}_{B}.
\label{eq:RSCextend2}
\end{align}
To show the RSC condition, we need to show that \eqref{eq:RSCextend2} is strictly positive. First, we obtain the sample complexity so that $A$ is far away from zero, then we show that $A$ is strictly greater than $B$. This is enough to obtain the sample complexity so that the RSC condition is satisfied.

\noindent {\bf{i. Lower bound on A:}}
Here, we explain how to get a lower bound on $\inf_{{\bf{u}} \in A} \sum_{i=1}^{n_2} \sigma_i \langle {\bf{u}}, {\bf{z}}_i \rangle^2$.

\noindent Let $r_i = r(X={\bf{x}_i^2} | \delta \theta^*)$, and $s_r = \sum_{j=1}^{n_2} r_j$, then $\sigma_i = \frac{r_i}{s_r}$. Then, we have
\begin{align}
\sum_{i=1}^{n_2} \sigma_i \langle {\bf{u}}, {\bf{z}}_i \rangle^2 = \frac{1}{s_r} \sum_{i=1}^{n_2} r_i \langle {\bf{u}}, {\bf{z}}_i \rangle^2~.
\end{align}

Then, we have
\begin{align}
p\left(\inf_{{\bf{u}} \in A} \frac{1}{s_r} \sum_{i=1}^{n_2} r_i \langle {\bf{u}}, {\bf{z}}_i \rangle^2 < ~ \frac{\eta_0}{\eta_1} \underline{\rho}^2 \left ( 1 - c\kappa_1^2 \frac{w(A)}{\sqrt{n_2}} \right) \right) &  \leq  p(\frac{1}{s_r} < \frac{1}{\eta_1 n_2}) \nonumber \\
& + p \left(\inf_{{\bf{u}} \in A} \sum_{i=1}^{n_2} r_i \langle {\bf{u}}, {\bf{z}}_i \rangle^2  ~< ~\eta_0 n_2 \underline{\rho}^2 \left ( 1 - c\kappa_1^2 \frac{w(A)}{\sqrt{n_2}} \right) \right).
\label{eq:bndA}
\end{align}

First, we give a bound for the first term.
Note that $E_{X \sim p(X|\theta_2)}[r(X | \delta \theta^*)]= 1$. From the smooth density ratio model assumption, we have
\begin{align}
p(|r_i - 1| > t) \leq 2 e^{-\frac{t^2}{2}}. \label{eq:r}
\end{align}
Applying Hoeffding inequality in \eqref{eq:r}, we have
\begin{align}
 & p(|\frac{1}{n_2} s_r -1| \geq t)  = p( |\frac{1}{n_2} \sum_{j=1}^{n_2} r_j - 1| \geq t) \leq 2 e^{-\frac{n_2 t^2}{2}}, \\
 \Rightarrow \quad &  p(s_r  \geq \eta_1 n_2) \leq  e^{-\frac{n_2 (\eta_1-1)^2}{2}}, \\
 \Rightarrow \quad &  p(\frac{1}{s_r} \leq \frac{1}{\eta_1 n_2}) = p(s_r  \geq \eta_1 n_2) \leq  e^{-\frac{n_2 (\eta_1-1)^2}{2}},
 \label{eq:probS}
\end{align}
where $\eta_1 = t+1$. 

Next, we focus on bounding the second term in \eqref{eq:bndA}. Recall that,
\begin{align}
&  \quad p(|r_i - 1| > t) \leq 2 e^{-\frac{t^2}{2}}, \label{eq:r}\\
\Rightarrow & \quad \bar{\varepsilon}_1 = p(r_i  \geq \eta_0)  \geq 1 - e^{-\frac{(1-\eta_0)^2}{2}},
\label{eq:probR}
\end{align}
where the last inequality holds for any $\eta_0 = 1-t$.

For any fixed $\eta_0$, let $\bar{W}_i = \bar{W}_i^u = \langle {\bf{u}}, {\bf{z}}_i \rangle \mathbb{I} (r_i > \eta_0)$. Then, the probability distribution over $\bar{W}_i$ can be written as:\footnote{With abuse of notation, we treat the distribution over $\bar{W}_i$ as discrete for ease of notation. A similar argument applies for the true continuous distribution, but more notation is needed.}
\beq
p(\bar{W}_i = w) = \frac{p(\langle {\bf{u}}, {\bf{z}}_i \rangle = w) \mathbb{I} (r_i > \eta_0)}{p(r_i > \eta_0)} \leq \frac{1}{\bar{\varepsilon}_1}p(\langle {\bf{u}}, {\bf{z}}_i \rangle = w)~.
\eeq
As a result, $\vertiii{\bar{W}_i}_{\psi_2} \leq \frac{\kappa}{\bar{\varepsilon}_1} = \kappa_1$. Thus, $\bar{W}_i = \bar{W}_i^u$ is a sub-Gaussian random variable for any ${\bf{u}} \in A$. Let $\rho_{\bf{u}}^2 = E[(\bar{W_i}^u)^2] > 0$.
For convenience of notation, let $Z_0$ be i.i.d. as the rows ${\bf{z}}_i,i = 1,\hdots,n$. Let $A \subseteq S^{p-1}$. Consider the following class of functions: 
\begin{align}
F = \{f_{\bf{u}}, {\bf{u}} \in A: f_{\bf{u}}(.)= \frac{1}{\rho_{\bf{u}}} \langle \cdot, {\bf{u}} \rangle \I(r(.|\delta\theta^*) \geq \eta_0)  : {\bf{u}} \in A \}.
\end{align}
Then for any $f_{\bf{u}} \in F$, $f_{\bf{u}}(Z_0) = \frac{1}{\rho_{\bf{u}}}\langle Z_0, {\bf{u}} \rangle \I(r_i \geq \eta_0)$ and, by construction, $F$ is a subset of the unit sphere, since for $f_{\bf{u}} \in F$
\begin{align}
\| f_{\bf{u}} \|_{L_2}^2 = \frac{1}{\rho_{\bf{u}}^2}E[\langle Z_0, {\bf{u}} \rangle^2 \I(r_i \geq \eta_0)] = 1.
\end{align}
 Further, $\sup_{f{\bf{u}} \in F} \vertiii{f_{\bf{u}}}_{\psi_2} \leq \kappa_1/2$.

Next, we show that for the current setting, the $\gamma_2$-functional can be upper bounded by $w(A)$, the Gaussian width of $A$. Since the process is sub-Gaussian with $\varphi_2$-norm bounded by $\kappa_1$, we have
\beq
\gamma_2(F \cap S_{L_2}, \vertiii{\cdot}_{\psi_2}) \leq \kappa_1 \gamma_2(F \cap S_{L_2},\vertiii{\cdot}_{L_2}) \leq \kappa_1 c_4 w(A)~,
\eeq
where the last inequality follows from generic chaining, in particular \cite[Theorem 2.1.1]{tala05}, for an absolute constant $c_4 > 0$.

In the context of Theorem~\ref{thm:subGREmain}, we choose
\beq
\theta = c_1 c_4 \kappa_1^2 \frac{w(A)}{\sqrt{n}} \geq c_1 \kappa_1 \frac{\gamma_2(F \cap S_{L_2},\vertiii{\cdot}_{\varphi_2})}{\sqrt{n}} ,
\eeq
so that the condition on $\theta$ is satisfied. With this choice of $\theta$, we have
\beq
\theta^2 n /\kappa_1^4 = c_1^2 c_4^2 w^2(A) ~.
\eeq
Then, from Theorem~\ref{thm:subGREmain}, it follows that with probability at least $1 - \exp(-\eta w^2(A))$, we have
\beq
\sup_{{\bf{u}} \in A} \left | \frac{1}{\rho_{\bf{u}} n_2} \sum_{i=1}^{n_2} \frac{1}{\rho_{\bf{u}}}\langle {\bf{z}}_i, {\bf{u}} \rangle^2 \I(r_i \geq \eta_0) - 1 \right | \leq c \kappa_1^2\frac{w(A)}{\sqrt{n_2}},
\eeq
where $\eta = c_2 c_1^2 c_4^2$ and $c = c_1 c_2$ are absolute constants. Thus, with probability at least $1 - \exp(-\eta w^2(A))$,
\begin{align}
& \quad \inf_{{\bf{u}} \in A} \frac{1}{n_2} \sum_{i=1}^{n_2} \langle {\bf{z}}_i, {\bf{u}} \rangle^2 \I(r_i \geq \eta_0) \geq \inf_{{\bf{u}} \in A} ~\rho_{\bf{u}}^2 \left( 1 - c  \kappa_1^2 \frac{w(A)}{\sqrt{n_2}} \right), \\
\Rightarrow & \quad \inf_{{\bf{u}} \in A} \sum_{i=1}^{n_2} \langle {\bf{z}}_i, {\bf{u}} \rangle^2 \I(r_i \geq \eta_0) \geq n_2 \underline{\rho}^2 \left( 1 - c  \kappa_1^2 \frac{w(A)}{\sqrt{n_2}} \right),
\end{align}
where $\underline{\rho}^2 = \inf_{{\bf{u}} \in A} \rho_{\bf{u}}^2$. Then, with probability at least $1 - \exp(-\eta w^2(A))$, we have
\begin{align}
\inf_{u \in A} \sum_{i=1}^{n_2} r_i \langle {\bf{z}}_i, {\bf{u}} \rangle^2 \geq &
\inf_{u \in A} \sum_{i=1}^{n_2} r_i \langle {\bf{z}}_i, {\bf{u}} \rangle^2 \I(r_i \geq \eta_0) \\
\geq & \inf_{u \in A} \eta_0 \sum_{i=1}^{n_2} \langle {\bf{z}}_i, {\bf{u}} \rangle^2 \I(r_i \geq \eta_0) \\
 \geq & \eta_0 n_2 \underline{\rho}^2 \left ( 1 - c\kappa_1^2 \frac{w(A)}{\sqrt{n_2}} \right)~.
\end{align}
Thus,
\begin{align}
& \quad p\left(\inf_{u \in A} \sum_{i=1}^{n_2} r_i \langle {\bf{z}}_i, {\bf{u}} \rangle^2 \geq \eta_0 n_2 \underline{\rho}^2 \left ( 1 - c\kappa_1^2 \frac{w(A)}{\sqrt{n_2}} \right)\right) \geq 1 - \exp(-\eta w^2(A)),  \\
\Rightarrow & \quad p\left(\inf_{u \in A} \sum_{i=1}^{n_2} r_i \langle {\bf{z}}_i, {\bf{u}} \rangle^2 < \eta_0 n_2 \underline{\rho}^2 \left ( 1 - c\kappa_1^2 \frac{w(A)}{\sqrt{n_2}} \right)\right) \leq \exp(-\eta w^2(A)).
\label{eq:2nd}
\end{align}

Putting \eqref{eq:probS} and \eqref{eq:2nd} into \eqref{eq:bndA}, for any $n_2 \geq \frac{2\eta w^2(A)}{(\eta_1-1)^2}$ we have
\begin{align}
p\left(\inf_{u \in A} \sum_{i=1}^{n_2} \sigma_i \langle {\bf{z}}_i, {\bf{u}} \rangle^2 < \frac{\eta_0}{\eta_1}  \underline{\rho}^2 \left ( 1 - c\kappa_1^2 \frac{w(A)}{\sqrt{n_2}} \right)\right) & \leq 
\exp(-\frac{n_2(\eta_1-1)^2}{2} ) + \exp(-\eta w^2(A)) \\
& \leq 2 \exp(-\eta w^2(A)),
\end{align}

{\bf{ii. A is strictly greater than B:}}
Note that, $0 \leq \sigma_i \leq 1$ for all $i$ and $\sum_{i=1}^{n_2} \sigma_i = 1$.
Define $f({\bf{z}}) = \langle {\bf{u}}, {\bf{z}} \rangle^2$, which is a convex function of ${\bf{z}}$. Using Jensen's inequality, we have
\begin{align}
f \left(\frac{\sum_{i=1}^{n_2} \sigma_i {\bf{z}}_i}{\sum_{i=1}^{n_2} \sigma_i} \right) &\leq \frac{\sum_{i=1}^{n_2} \sigma_i  f ({\bf{z}}_i)}{\sum_{i=1}^{n_2} \sigma_i} \\
\left \langle {\bf{u}}, \frac{\sum_{i=1}^{n_2} \sigma_i {\bf{z}}_i}{\sum_{i=1}^{n_2} \sigma_i} \right \rangle^2 &\leq \frac{\sum_{i=1}^{n_2} \sigma_i  \langle {\bf{u}}, {\bf{z}}_i \rangle^2}{\sum_{i=1}^{n_2} \sigma_i} \\
 \langle {\bf{u}}, \sum_{i=1}^{n_2} \sigma_j {\bf{z}}_i \rangle^2 &\leq \sum_{i=1}^{n_2} \sigma_i \langle {\bf{u}},  {\bf{z}}_i \rangle^2.
 \label{eq:jensenInqSecDerL}
\end{align}
The equality in \eqref{eq:jensenInqSecDerL} holds if ${\bf{z}}_1 = {\bf{z}}_2 = \cdots = {\bf{z}}_{n_2}$, or if both sides are zero i.e., ${\bf{u}}$ is in the null space of ${\bf{z}}_i$ for all $i$.
Since ${\bf{z}}_i$ are different with probability 1, then if we show that ${\bf{u}}$ is not in the null space of ${\bf{z}}_i$ for all $i$, then the inequality \eqref{eq:jensenInqSecDerL} is strict inequality.

This completes the proof. \qed

\section*{Acknowledgment}
The research was supported by NSF grants IIS-1447566, IIS-1447574, IIS-1422557, CCF-1451986, CNS- 1314560, IIS-0953274, IIS-1029711, NASA grant NNX12AQ39A, and gifts from Adobe, IBM, and Yahoo. F. F. acknowledges the support of IDF (2014-2015) and DDF (2015-2016) from the University of Minnesota.

\bibliographystyle{abbrv}
\bibliography{Sparsity,mutlitaskGraph,new_ref}

\end{document}